\documentclass[11pt,reqno]{amsart}
\usepackage[tt=false]{libertine}
\usepackage{amssymb}
\usepackage{mathtools}
\usepackage[varbb]{newpxmath}
\usepackage[foot]{amsaddr}

\let\savedbigtimes\bigtimes
\let\bigtimes\relax
\usepackage{mathabx}
\let\bigtimes\savedbigtimes

\usepackage[margin=1in, bottom=1in]{geometry}
\usepackage{graphicx}
\usepackage{caption}
\usepackage{subcaption}
\usepackage{enumitem}
\usepackage{comment}
\usepackage{bbm}

\usepackage[usenames,dvipsnames]{xcolor}
\usepackage[colorlinks=true,
  linkcolor=teal!60!black,
  citecolor=PineGreen,
  urlcolor=RedViolet]{hyperref}
\hypersetup{bookmarksopen=true}

\usepackage[yyyymmdd,hhmmss]{datetime}

\usepackage{tikz}
\usetikzlibrary{decorations.pathreplacing,calligraphy}
\usetikzlibrary{calc}

\newtheorem{thm}{Theorem}[section]
\newtheorem{lem}[thm]{Lemma}
\newtheorem{ppn}[thm]{Proposition}

\newtheorem{cor}[thm]{Corollary}
\newtheorem{fac}[thm]{Fact}

\theoremstyle{definition}
\newtheorem{dfn}[thm]{Definition}
\newtheorem{ex}[thm]{Example}
\newtheorem{asn}[thm]{Assumption}

\theoremstyle{remark}
\newtheorem{rmk}[thm]{Remark}

\numberwithin{equation}{section}
\def\beq#1\eeq{%
    \begin{equation}%
    #1%
    \end{equation}%
}
\def\baln#1\ealn{%
    \begin{align}%
    #1%
    \end{align}%
}
\def\balnn#1\ealnn{%
    \begin{align*}%
    #1%
    \end{align*}%
}

\def\lt{\left}
\def\rt{\right}
\def\fr{\frac}
\def\la{\langle}
\def\ra{\rangle}
\def\eps{\varepsilon}

\def\bA{{\boldsymbol{A}}}
\def\bB{{\boldsymbol{B}}}
\def\bC{{\boldsymbol{C}}}
\def\bD{{\boldsymbol{D}}}
\def\bE{{\boldsymbol{E}}}

\def\bG{{\boldsymbol{G}}}
\def\bH{{\boldsymbol{H}}}
\def\bI{{\boldsymbol{I}}}

\def\bK{{\boldsymbol{K}}}
\def\bL{{\boldsymbol{L}}}
\def\bM{{\boldsymbol{M}}}

\def\bO{{\boldsymbol{O}}}
\def\bP{{\boldsymbol{P}}}

\def\bR{{\boldsymbol{R}}}
\def\bS{{\boldsymbol{S}}}
\def\bT{{\boldsymbol{T}}}
\def\bU{{\boldsymbol{U}}}
\def\bV{{\boldsymbol{V}}}
\def\bW{{\boldsymbol{W}}}
\def\bX{{\boldsymbol{X}}}
\def\bY{{\boldsymbol{Y}}}
\def\bZ{{\boldsymbol{Z}}}

\def\bLambda{{\boldsymbol{\Lambda}}}

\def\bzero{{\boldsymbol{0}}}

\def\bc{{\boldsymbol{c}}}

\def\be{{\boldsymbol{e}}}

\def\bg{{\boldsymbol{g}}}

\def\bu{{\boldsymbol{u}}}
\def\bv{{\boldsymbol{v}}}
\def\bw{{\boldsymbol{w}}}
\def\bx{{\boldsymbol{x}}}
\def\by{{\boldsymbol{y}}}

\def\bbC{{\mathbb{C}}}

\def\bbE{{\mathbb{E}}}

\def\bbP{{\mathbb{P}}}

\def\bbR{{\mathbb{R}}}
\def\bbS{{\mathbb{S}}}

\def\bbZ{{\mathbb{Z}}}
\def\bbone{{\mathbb{1}}}

\def\cE{{\mathcal{E}}}
\def\cF{{\mathcal{F}}}
\def\cG{{\mathcal{G}}}

\def\cI{{\mathcal{I}}}

\def\cK{{\mathcal{K}}}

\def\cN{{\mathcal{N}}}

\def\cQ{{\mathcal{Q}}}

\def\adm{{\mathsf{adm}}}
\def\bip{{\mathsf{bip}}}
\def\de{{\mathsf{d}}}
\def\He{{\mathsf{He}}}
\def\hop{{\mathsf{hop}}}
\def\Is{{\mathsf{Is}}}
\def\per{{\mathsf{per}}}
\def\mix{{\mathsf{mix}}}
\def\mnr{{\mathsf{mnr}}}
\def\MP{{\mathsf{MP}}}
\def\nul{{\mathsf{nul}}}
\def\sgn{{\mathsf{sgn}}}
\def\spec{{\mathsf{spec}}}
\def\sph{{\mathsf{sph}}}
\def\supp{{\mathsf{supp}}}
\def\Tr{{\mathsf{Tr}}}
\def\unif{{\mathsf{unif}}}
\def\Var{{\mathsf{Var}}}

\def\oq{{\bar q}}

\def\bSigma{{\boldsymbol{\Sigma}}}

\title{Free energy fluctuations of symmetric perceptron models up to and including criticality}

\author{Hang Du \and Shuyang Gong \and Brice Huang \and Zhangsong Li}

\address[H.~Du]{Department of Mathematics, Massachusetts Institute of Technology}
\email[H.~Du]{hangdu@mit.edu}

\address[S.~Gong]{Institute for Foundations of Data Science, Yale University; Department of Statistics and Data Science, Yale University}
\email[S.~Gong]{shuyang.gong@yale.edu}

\address[B.~Huang]{Department of Statistics, Stanford University}
\email[B.~Huang]{bmhuang@stanford.edu}

\address[Z.~Li]{School of Mathematical Sciences, Peking University}
\email[Z.~Li]{ramblerlzs@pku.edu.cn}

\date{\today}


\begin{document}

\begin{abstract}
  We study the free energy fluctuations of spherical and Ising perceptron models with general even activations.
  These include the symmetric binary perceptron, the $U$-perceptron, and the bipartite Sherrington--Kirkpatrick model.
  Our main result is a Gaussian CLT for the free energy that holds both in the ``high-temperature'' interior of a second moment RS regime and on its ``critical'' boundary.
  At criticality, the CLT has scale $(\frac{\log N}{6})^{1/2}$, which matches the free energy fluctuations of the critical Sherrington--Kirkpatrick model and suggests this scaling is universal.
  
  The proof shows that the perceptron model's free energy fluctuations are explained, up to an additive $o(1)$ error with high probability, by those of a corresponding spherical Hopfield model, and that this relation persists even at criticality.
  This proof strategy builds on a critical reweighted moment method developed by the first and third authors \cite{du2026fluctuations} to identify the critical fluctuations of the SK free energy, which links these fluctuations to the spectral edge spacing of the disorder.
\end{abstract}

\maketitle

\setcounter{tocdepth}{1}
\tableofcontents

\section{Introduction and main results}
\label{s:intro}

Let $\phi:\bbR \rightarrow[0,+\infty)$ be an even, piecewise continuous function with finitely many pieces, satisfying
\begin{equation}\label{e:sub-quadratic-exponential-growth}
  \limsup_{x\to+\infty} \fr{(\log \phi(x))_+}{x^2} = 0\,.
\end{equation}
Here $x_+ = \max(x,0)$.
Throughout, we fix $\phi$ and $\alpha \in (0,+\infty)$ and let $N \rightarrow \infty$ be a positive integer.
Let $M = \lfloor \alpha N \rfloor$ and $\bG \in \bbR^{M\times N}$ have independent standard Gaussian entries, with rows $\bg^1,\ldots,\bg^M \in \bbR^N$.

Set $S_N=\sqrt{N}\bbS^{N-1}$, $\Sigma_N=\{-1,+1\}^N$, and let $\nu_N$ be the normalized surface measure on $S_N$.
In this paper we study perceptron models with general even activation $\phi$ and spherical or Ising spins.
These have partition functions
\balnn
  Z_{\per,N}^{\sph}(\bG)
  &= \int_{S_N}\prod_{a=1}^M
  \phi\lt(\fr{(\bg^a,\bx)}{\sqrt N}\rt)
  \,\de\nu_N(\bx)\,, \\
  Z_{\per,N}^{\Is}(\bG)
  &=\fr1{2^N}\sum_{\bx\in\Sigma_N}\prod_{a=1}^M
  \phi\lt(\fr{(\bg^a,\bx)}{\sqrt N}\rt)\,,
\ealnn
respectively.
These include the symmetric Ising perceptron and $U$-perceptron, corresponding to the activations $\phi(x)=\bbone\{|x|\leq\kappa\}$ and $\phi(x)=\bbone\{|x|\geq\kappa\}$ \cite{aubin2019storage,perkins2021frozen,abbe2022proof}.
Among non-even activations, the well-studied classical perceptron \cite{wendel1962problem,cover1965geometrical,krauth1989storage, gardner1988optimal,ding2018capacity,xu2021sharp,huang2024capacity} corresponds to $\phi(x)=\bbone\{x\geq\kappa\}$ for some $\kappa\in\bbR$.
More general activations were also studied in \cite{bolthausen2022gardner,nakajima2023sharp}.

Since the introduction of this model, a central question has been to identify the in-probability limit of the free energy density $\fr1N \log Z_{\per,N}^{\ast}$, in both settings $\ast \in \{\sph, \Is\}$.
This limit has now been established for many perceptron models, including the symmetric Ising perceptron for all $\kappa$ \cite{aubin2019storage,perkins2021frozen,abbe2022proof}.
Characterizing this limit for general $\phi$ was a significant open problem, even within the class of even $\phi$, which was resolved only recently in \cite{openai2026spherical, openai2026ising} following the Hamilton--Jacobi approach to the Parisi formula developed by \cite{mourrat2021hamilton, mourrat2022parisi, dominguez2024statistical, chen2026free}.

Another natural question is to describe the limiting distribution of the free energy $\log Z_{\per,N}^{\ast}(\bG)$.
For the symmetric Ising perceptron, Abbe, Li, and Sly \cite{abbe2022proof} showed that the recentered free energy satisfies a Gaussian CLT, with an explicit variance of scale $O(1)$ depending on $\kappa$.

We study this problem for general $\phi$.
Our main result identifies the fluctuations for $\phi, \alpha$ in a second moment RS regime, see Assumption~\ref{a:2mt} below.
The setting of \cite{abbe2022proof} lies in the interior of this regime.
In this interior, the model has ``high-temperature'' behavior, and we show the free energy continues to satisfy a Gaussian CLT, with explicit variance of scale $O(1)$ depending on $\phi$.

\textbf{The main novelty of this paper is that our main result also holds on a ``critical'' boundary of the second moment RS regime.}
The model's behavior becomes qualitatively different at criticality: we show that the recentered free energy, divided by $(\fr{\log N}{6})^{1/2}$, converges in distribution to a standard Gaussian.
This scaling matches that of the critical Sherrington--Kirkpatrick model \cite{sherrington1975solvable}, whose limiting free energy fluctuations were predicted by \cite{aspelmeier2008free} and recently confirmed by the first and third authors \cite{du2026fluctuations}.
Our result thus suggests that this behavior is universal in a broad class of critical models.

We now state our result formally.
For $q\in[-1,1]$ let $\bbE_q$ denote expectation with respect to centered Gaussian $(Z_1,Z_2)$ with unit variances and correlation $q$, and let
\baln\label{e:m-and-Phi}
  m_\phi &= \bbE[\phi(Z)]\,, & 
  \Phi_{\per}(q) &= \fr{\bbE_q[\phi(Z_1)\phi(Z_2)]}{(m_\phi)^2}\,.
\ealn
Then define the normalized logarithmic second moment rates for the two models
\baln
  \label{e:cF-sph}
  \cF_{\sph}(q)
  &=\fr12\log(1-q^2)+\alpha\log \Phi_{\per}(q)\,, \\
  \label{e:cF-Is}
  \cF_{\Is}(q)
  &=-\fr12\big((1+q)\log(1+q)+(1-q)\log(1-q)\big)
  +\alpha\log \Phi_{\per}(q)\,.
\ealn
In \eqref{e:cF-sph} we take $\log 0 = -\infty$, and in \eqref{e:cF-Is} we take $0 \log 0 = 0$.
\begin{asn}[Second moment RS regime]
  \label{a:2mt}
  For the spherical and Ising perceptron models respectively, we impose the following hypotheses.
  \begin{enumerate}[label=(\alph*)]
    \item \label{i:2mt-sphere} The function $\cF_{\sph}$ is uniquely maximized over $[-1,1]$ at $q=0$.
    \item \label{i:2mt-ising} The function $\cF_{\Is}$ is uniquely maximized over $[-1,1]$ at $q=0$.
  \end{enumerate}
\end{asn}
\noindent For $Z\sim\cN(0,1)$, define
\beq\label{e:s}
  s = \fr{\bbE[Z^2\phi(Z)]}{\bbE[\phi(Z)]} \in (0, +\infty)\,.
\eeq
Then let
\[
  \theta 
  = 1 + \cF''_{\sph}(0)
  = 1 + \cF''_{\Is}(0)
  = \alpha (s-1)^2\,.
\]
Since $\cF_{\sph}$ and $\cF_{\Is}$ are clearly even, under Assumption~\ref{a:2mt} we have $\theta \in [0,1]$.
The cases $\theta < 1$ and $\theta = 1$ correspond to the qualitatively different high-temperature and critical behaviors alluded to above.
Our main result, which addresses both cases, is the following.
\begin{thm}[Perceptron free energy CLT]\label{t:main}
  Define $\sigma > 0$ by
  \beq\label{e:sigma}
    \sigma^2
    = \begin{cases}
      -\fr12 \log (1-\theta) & \theta < 1\,, \\ 
      \fr16 \log N & \theta = 1\,.
    \end{cases}
  \eeq
  As $N\rightarrow\infty$, the following limits hold in distribution.
  \begin{enumerate}[label=(\alph*)]
    \item \label{i:main-sphere} Under Assumption~\ref{a:2mt}\ref{i:2mt-sphere}, when $s\neq 1$, then $\sigma\neq 0$ and the spherical free energy satisfies
    \[
      \sigma^{-1} \lt(
        \log Z_{\per,N}^{\sph} - M\log m_\phi
        + \fr12 \sigma^2
      \rt)
      \rightarrow \cN(0,1)\,,
    \]
    and when $s=1$, $\log Z_{\per,N}^{\sph} - M\log m_\phi \stackrel{\mathrm{P}}{\rightarrow} 0$. 
    
    \item \label{i:main-ising} Under Assumption~\ref{a:2mt}\ref{i:2mt-ising}, when $s\neq 1$, then $\sigma\neq 0$ and the Ising free energy satisfies
    \[
      \sigma^{-1} \lt(
        \log Z_{\per,N}^{\Is} - M\log m_\phi
        + \fr12 \sigma^2
      \rt)
      \rightarrow \cN(0,1)\,,
    \]
    and when $s=1$, $\log Z_{\per,N}^{\Is} - M\log m_\phi\stackrel{\mathrm{P}}{\rightarrow}  0$.
  \end{enumerate}
\end{thm}
\begin{rmk}
  More generally, we can consider a sequence of constraint densities $M/N = \alpha_N$ with $\alpha_N \to \alpha \in (0,\infty)$.
  If $\theta = 1$, the resulting CLT will also depend on the rate at which $\theta_N = \alpha_N (s-1)^2$ tends to $1$.
  For all $\theta_N \le 1 + O(N^{-1/3})$, our methods can be extended to show an analogous CLT with scale
  \[
    \sigma_N^2 = -\fr12 \log \max(1-\theta_N, N^{-1/3})
  \]
  in place of $\sigma^2$.
  For simplicity, we do not carry out this extension in this paper.
\end{rmk}

\subsection{Applications of the main result}
\label{ss:applications}

We next discuss Assumption~\ref{a:2mt} and give several examples illustrating the applications of Theorem~\ref{t:main}.
For any $\phi$, in either the spherical or Ising setting, the model lies in the high-temperature interior of the regime described by Assumption~\ref{a:2mt} when $\alpha > 0$ is sufficiently small.
As we increase $\alpha$, we eventually reach the boundary of this regime in one of two ways:
\begin{enumerate}[label=(S\arabic*)]
  \item \label{i:critical-away-from-0} There is a nonzero $q$ such that $\cF_\ast(q) = \cF_\ast(0) = 0$, for the appropriate $\ast \in \{\sph, \Is\}$.
  Increasing $\alpha$ further causes $\cF_\ast(q)$ to become larger than $\cF_\ast(0) = 0$.
  \item \label{i:critical-at-0} $\cF''_\ast(0) = 0$.
  Increasing $\alpha$ further causes $0$ to become a local minimum of $\cF_\ast$.
\end{enumerate}
\begin{figure}[htbp]
  \centering
  \begin{subfigure}[t]{0.49\textwidth}
    \centering
    \includegraphics[width=\linewidth]{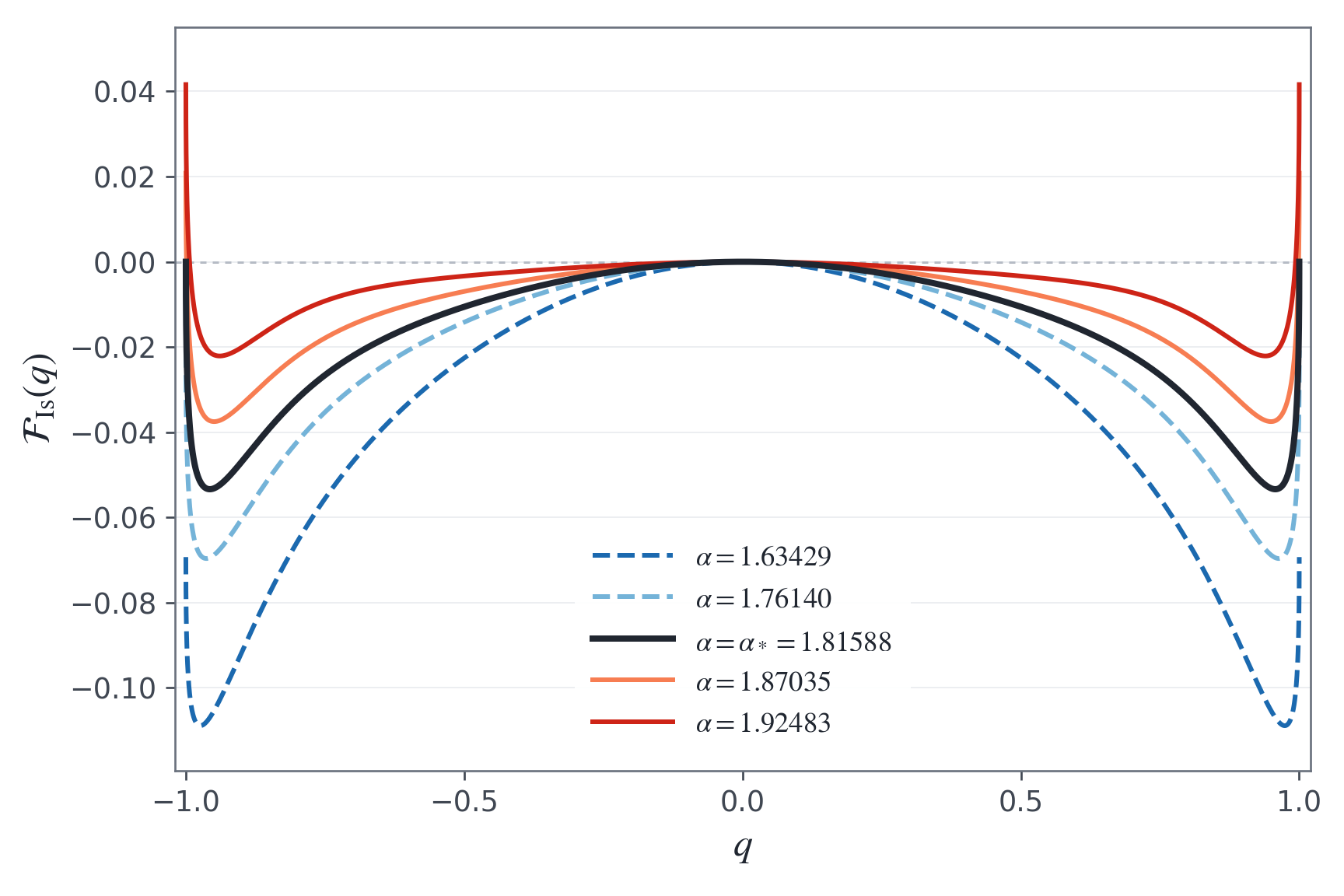}
    \caption{Symmetric Ising perceptron, $\kappa = 1$.}
    \label{f:symmetric-ising-perceptron-rate}
  \end{subfigure}\hfill
  \begin{subfigure}[t]{0.49\textwidth}
    \centering
    \includegraphics[width=\linewidth]{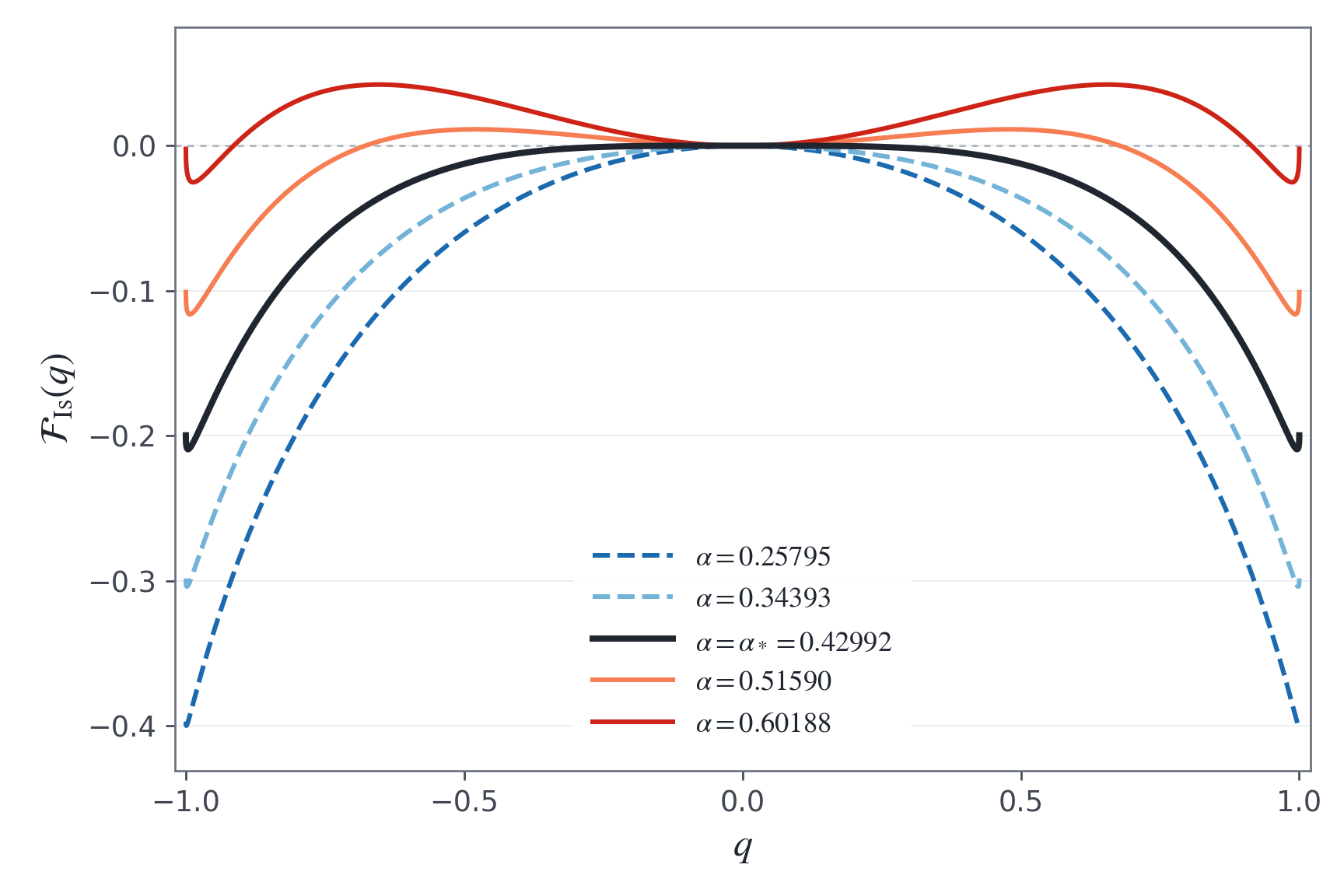}
    \caption{$U$-perceptron, $\kappa = 1$.}
    \label{f:u-perceptron-rate}
  \end{subfigure}
  \caption{
    Plots of $\cF_\Is(q)$ for two perceptron models at various constraint densities $\alpha$.
    The black curve marks the density $\alpha_*$ where the model exits the Assumption~\ref{a:2mt} regime.
    In Figure~\ref{f:symmetric-ising-perceptron-rate}, this occurs in scenario~\ref{i:critical-away-from-0}, and in Figure~\ref{f:u-perceptron-rate} it occurs in scenario~\ref{i:critical-at-0}.}
  \label{f:perceptron-rate-scenarios}
\end{figure}

\begin{figure}[htbp]
  \centering
  \begin{subfigure}[t]{0.325\textwidth}
    \centering
    \includegraphics[width=\linewidth]{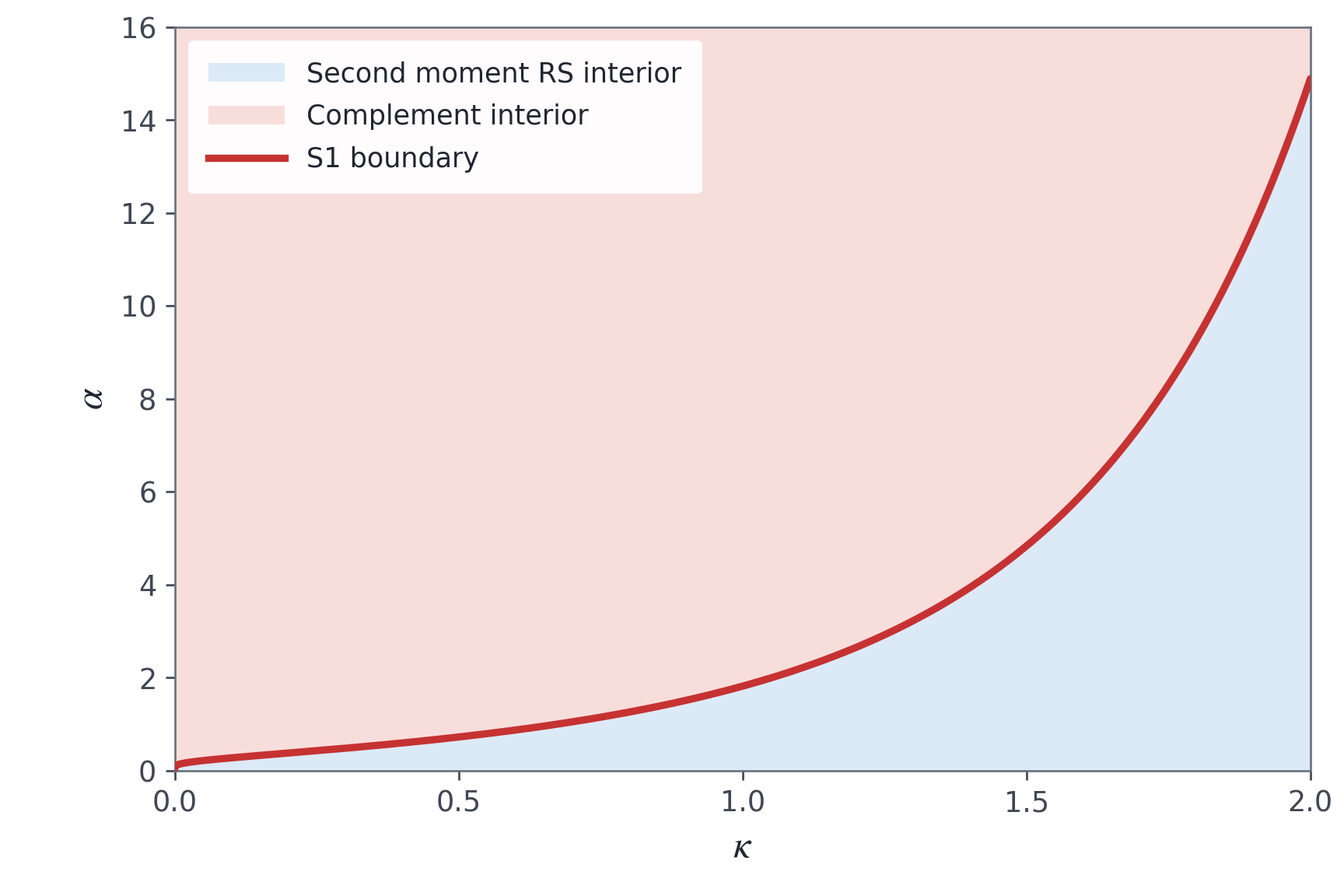}
    \caption{Symmetric Ising perceptron.}
    \label{f:symmetric-ising-perceptron-phase-diagram}
  \end{subfigure}
  \begin{subfigure}[t]{0.325\textwidth}
    \centering
    \includegraphics[width=\linewidth]{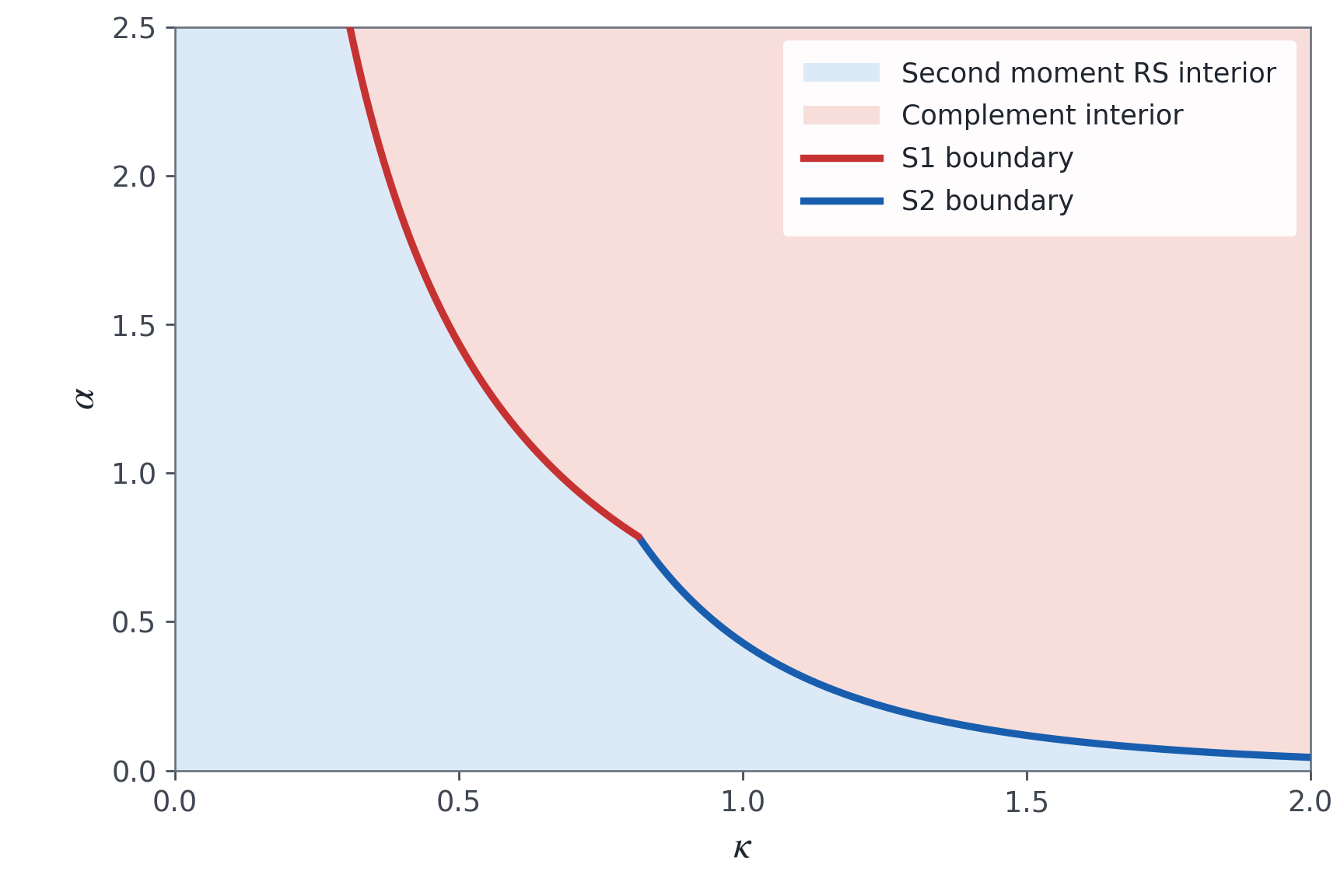}
    \caption{$U$-perceptron.}
    \label{f:u-perceptron-phase-diagram}
  \end{subfigure}
  \begin{subfigure}[t]{0.325\textwidth}
    \centering
    \includegraphics[width=\linewidth]{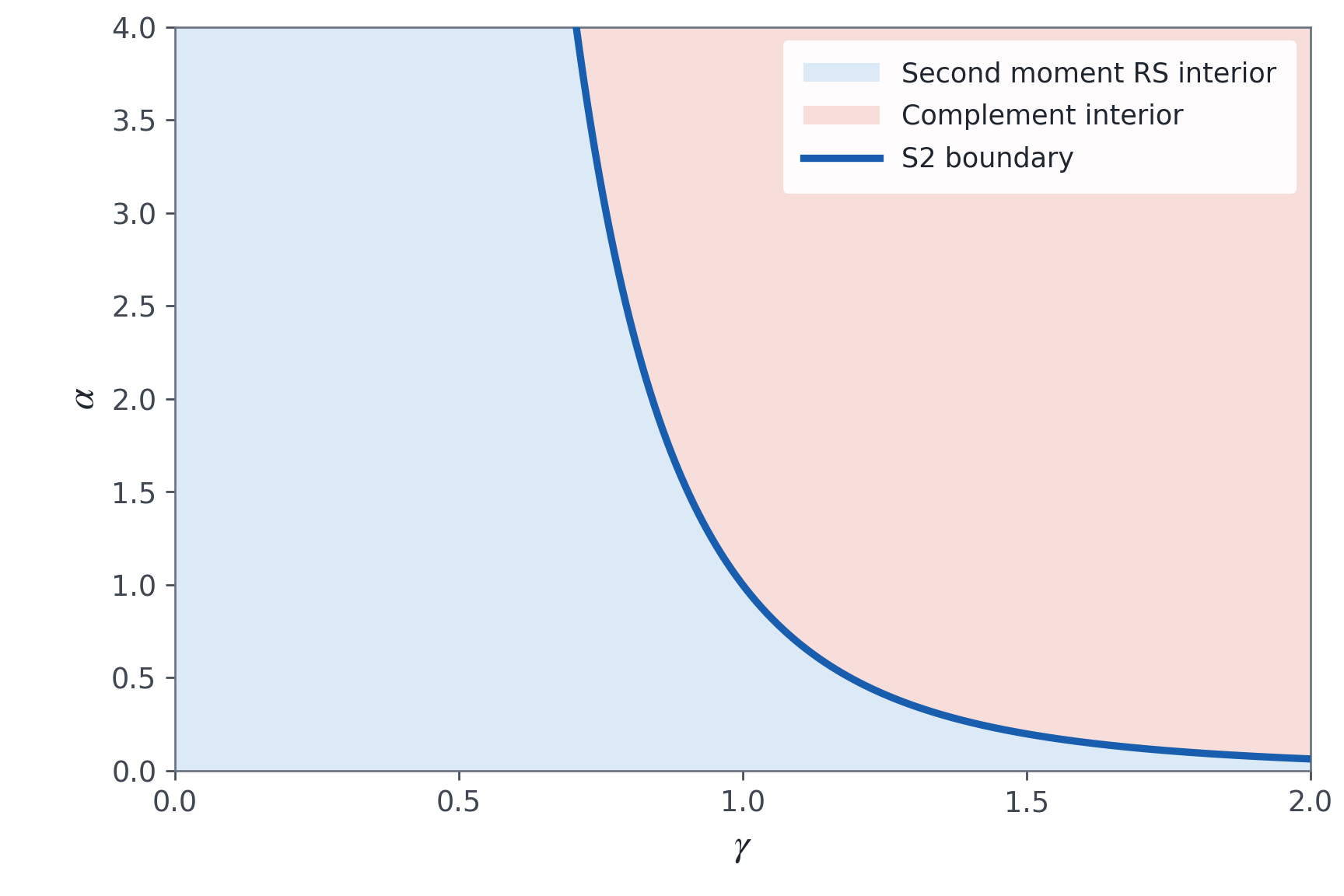}
    \caption{Bipartite SK model.}
    \label{f:bipartite-sk-phase-diagram}
  \end{subfigure}
  \caption{Phase diagrams of the three examples. The horizontal axis varies the model parameter $\kappa$ or $\gamma$, while the vertical axis varies $\alpha$. The interior of the Assumption~\ref{a:2mt} regime is shaded blue, while the interior of its complement is shaded red. The boundary is colored red if scenario \ref{i:critical-away-from-0} occurs, and blue if scenario \ref{i:critical-at-0} occurs. Theorem~\ref{t:main} holds in the blue shaded region and blue part of the boundary.}
  \label{f:phase-diagrams}
\end{figure}

The critical case of Theorem~\ref{t:main} applies when we reach the boundary through scenario \ref{i:critical-at-0}, and we view this as the most interesting case of this theorem.
We next examine which of these two scenarios occurs in several examples of random perceptron models, and plot their phase diagrams in Figure~\ref{f:phase-diagrams}.

\begin{ex}[Symmetric Ising perceptron]
  \label{e:symmetric-ising-perceptron}
  This model corresponds to the choice $\phi(x) = \bbone\{|x| \le \kappa\}$.
  As shown in \cite{abbe2022proof}, this model satisfies Assumption~\ref{a:2mt} for all $\alpha$ below the satisfiability threshold, at which point $\cF_\Is(1) = \cF_\Is(0) = 0$.
  It thus exits the regime defined by Assumption~\ref{a:2mt} through scenario \ref{i:critical-away-from-0} for all $\kappa > 0$.
\end{ex}

\begin{ex}[$U$-perceptron]
  \label{e:u-perceptron}
  This model, studied in \cite{aubin2019storage}, corresponds to $\phi(x) = \bbone\{|x| \ge \kappa\}$.
  \cite{aubin2019storage} numerically identified a threshold $\kappa_U \approx 0.8171$ such that scenario \ref{i:critical-away-from-0} holds for $\kappa \le \kappa_U$, and scenario \ref{i:critical-at-0} holds for $\kappa \ge \kappa_U$.
  The critical case of Theorem~\ref{t:main} applies for $\kappa > \kappa_U$.
\end{ex}

\begin{ex}[Bipartite SK model]
  \label{e:bipartite-sk}
  Let $\gamma > 0$ and $\phi(x) = \cosh(\gamma x)$.
  This choice of activation encodes the bipartite SK model at inverse temperature $\gamma$, which has partition function
  \[
    Z^\Is_{\bip,N}
    = \fr{1}{2^{N+M}}
    \sum_{\bx \in \Sigma_N} \sum_{\by \in \Sigma_M} 
    \exp\lt(
      \fr{\gamma}{\sqrt{N}}
      (\bG \bx, \by)
    \rt)
    = \fr{1}{2^N}
    \sum_{\bx \in \Sigma_N}
    \prod_{a=1}^M
    \cosh\lt(\fr{\gamma (\bg^a, \bx)}{\sqrt{N}}\rt).
  \]
  Scenario \ref{i:critical-at-0} holds for all $\gamma > 0$.
  Indeed, we can directly calculate 
  \balnn
    s &= 1 + \gamma^2, & 
    \Phi_\per(q) &= \cosh(\gamma^2 q).
  \ealnn
  Thus $\theta = \alpha (s-1)^2$ equals $1$ when $\alpha = \gamma^{-4}$.
  We will show that Assumption~\ref{a:2mt} holds for all $\alpha \le \gamma^{-4}$.
  Note the elementary inequalities 
  \balnn
    (1+q)\log(1+q)+(1-q)\log(1-q) &\ge q^2, &
    \log\cosh x &\le x^2/2,
  \ealnn
  with equality only at $q=0$ and $x=0$ respectively.
  Thus for all $q\neq 0$,
  \[
    \cF_\Is(q) < -\fr{q^2}{2} + \alpha \cdot \fr{(\gamma^2 q)^2}{2}
    = -\fr12 (1 - \alpha \gamma^4) q^2 
    \le 0 = \cF_\Is(0).
  \]
\end{ex}

\subsection{Proof method: comparison with Hopfield model}

We prove Theorem~\ref{t:main} by comparing the perceptron models with a spherical Hopfield model defined with the same disorder $\bG$.
Set
\beq\label{e:beta}
  \beta=\fr{s-1}{s}.
\eeq
In particular, $\beta<0$ when $s \in (0,1)$ and $\beta \in (0,1)$ when $s>1$.
The spherical Hopfield partition function at inverse temperature $\beta$ is
\beq\label{e:hopfield}
  Z_{\hop,N}^{\sph}(\bG)
  =\int_{S_N}
  \exp\lt(\fr{\beta}{2N}\|\bG\bx\|^2\rt)
  \,\de\nu_N(\bx).
\eeq
The choice \eqref{e:beta} of $\beta$ makes possible the following comparison theorem, which shows that the fluctuations of the perceptron and Hopfield models coincide in probability.
In Remark~\ref{r:beta-choice} below, we explain why this particular choice of $\beta$ enables this comparison.
\begin{thm}
  \label{t:perceptron-comparison}
  For $\beta=(s-1)/s$ and $\ast \in \{\sph, \Is\}$, define
  \[
      X_N^{\ast}
      = \fr{Z_{\per,N}^{\ast}(\bG)}{Z_{\hop,N}^{\sph}(\bG)}
      \cdot \fr{\bbE Z_{\hop,N}^{\sph}(\bG)}{\bbE Z_{\per,N}^{\ast}(\bG)}\,.
  \]
  Then, the following holds.
  \begin{enumerate}[label=(\alph*)]
    \item \label{i:spherical-comparison} Under Assumption~\ref{a:2mt}\ref{i:2mt-sphere}, $\bbE[(X_N^{\sph}-1)^2] \rightarrow 0$.
    \item \label{i:ising-comparison} Under Assumption~\ref{a:2mt}\ref{i:2mt-ising}, $\bbE[(X_N^{\Is}-1)^2] \rightarrow 0$.
  \end{enumerate}
\end{thm}
We compare to a spherical Hopfield model because its partition function has a random matrix interpretation, due to \cite{baik2016fluctuations} and developed in \cite{baik2018ferromagnetic,baik2020free,baik2021spherical,landon2022free,johnstone2024spin,collins2025fluctuations}, which makes it possible to characterize its free energy fluctuations.
This connection is very general, and applies to spherical models with quadratic interaction whose interaction matrix comes from a broad class of random matrix ensembles.
The following theorem is proved in Section~\ref{rewrite:a:hopfield-clt} by adapting ideas from this literature, and identifies the free energy fluctuations of the model we compare to.

\begin{thm}[Spherical Hopfield free energy CLT]\label{t:hopfield-clt}
  For $0<\theta\le 1$, let $\sigma$ be defined in \eqref{e:sigma}.
  When $s\neq 1$, then $\sigma\neq 0$ and for $\beta=(s-1)/s$,
  \[
    \sigma^{-1}\lt(
      \log Z_{\hop,N}^{\sph}-\fr{M}{2} \log s
      + \fr12\sigma^2
    \rt)
    \rightarrow \cN(0,1)
  \]
  in distribution. When $s=1$, $\log Z_{\hop,N}^\sph=0$. 
\end{thm}

At criticality, the scaling $\sigma^2 = \fr16 \log N$ matches the one found by Landon \cite{landon2022free} for the critical Sherrington--Kirkpatrick model. 
This scaling is universal, and comes from the $N^{-2/3}$ edge eigenvalue spacing common to a broad class of random matrix ensembles; in \cite{landon2022free} the relevant random matrix is the GOE, and in this paper it is the Wishart ensemble.

\begin{proof}[Proof of Theorem~\ref{t:main}]
  Note that the shifts $M \log m_\phi$ and $\fr{M}{2} \log s$ in Theorems~\ref{t:main} and \ref{t:hopfield-clt} are simply the annealed free energies
  \baln\label{e:annealed-free-energies}
    \log \bbE Z_{\per,N}^{\Is}(\bG)
    = \log \bbE Z_{\per,N}^{\sph}(\bG)
    &= M \log m_\phi\,, &
    \log \bbE Z_{\hop,N}^{\sph}(\bG)
    &= \fr{M}{2} \log s\,.
  \ealn
  The result is now immediate from Theorems~\ref{t:perceptron-comparison} and \ref{t:hopfield-clt}.
\end{proof}

\begin{rmk}\label{r:beta-choice}
  Let $\varphi(z) = (2\pi)^{-1/2} \exp(-z^2/2)$ denote the standard Gaussian density.
  We choose $\beta$ as in \eqref{e:beta} so that the probability measures on $\bbR$ with densities proportional to $\phi(z) \varphi(z)$ and $\exp(\beta z^2/2) \varphi(z)$ have the same second moment, namely $s = \fr{1}{1-\beta}$.
  The statement of Theorem~\ref{t:main} suggests that the fluctuation scale $\sigma$ depends on the activation $\phi$ only through the reweighted second moment $s$.
  Accordingly, the special case $\phi(z) = \exp(\beta z^2/2)$ motivates choosing $\beta$ to match these second moments, so that the perceptron and Hopfield models have the same fluctuation behavior.
\end{rmk}

\begin{rmk}\label{r:comparison-principle-general}
  Theorem~\ref{t:perceptron-comparison} builds on a comparison idea introduced by the first and third authors in \cite{du2026fluctuations}, where it was shown that free energy fluctuations of the critical (Ising) SK model are explained by those of the critical spherical SK model up to an additive $o(1)$ error.
  \textbf{The main contribution of this paper is to show that this is an instance of a much more general principle, which extends well beyond relating an Ising model to its spherical counterpart.}
  The more general comparison target is a spherical model with a quadratic random matrix interaction, which captures the original model's free energy fluctuations at both high temperature and criticality.
  Comparing to such a target is useful when the original model has either Ising or spherical spins: here, the spherical Hopfield model captures the free energy fluctuations of both versions of the perceptron.
  The comparison model's free energy fluctuations can be analyzed by its aforementioned random matrix interpretation.
  We expect that these ideas provide a way to identify the free energy fluctuations of a broad class of spin systems at high temperature and criticality, and in particular to show the universal $\fr16 \log N$ variance scale at criticality.
  Establishing this more general comparison principle also requires several innovations in the proof, see \S\ref{ss:proof-discussion} for further discussion.
\end{rmk}

\subsection{Related work}

In models with spherical spins and quadratic random matrix interaction, such as the spherical SK model, the free energy fluctuations are mostly understood.
In such models, \cite[Lemma 1.3]{baik2016fluctuations} introduced a formula expressing the free energy in terms of a contour integral depending on the eigenvalues of the disorder, through which the free energy can be studied using random matrix theory (see Lemma~\ref{l:hopfield-coarea}).
Using this connection, \cite{baik2016fluctuations} showed that the spherical SK free energy has Gaussian fluctuations of order $1$ at high temperature $\beta < 1$ and Tracy--Widom fluctuations of scale $N^{1/6}$ at low temperature $\beta > 1$.
At criticality $\beta = 1$, \cite{landon2022free} showed Gaussian fluctuations of scale $(\frac{\log N}{6})^{1/2}$, while \cite{johnstone2024spin} identified the fluctuations on the critical window $\beta = 1 + b N^{-1/3} \sqrt{\log N}$, over which the fluctuations transition from Gaussian to Tracy--Widom.
These ideas have also been used to study spherical spin glasses with ferromagnetic interaction \cite{baik2018ferromagnetic}, external field \cite{baik2021spherical}, and bipartite spins \cite{baik2020free, collins2025fluctuations}.
We expect these phenomena to be universal and adapt the above ideas to prove Theorem~\ref{t:hopfield-clt}.

At high temperature, the free energy fluctuations of a broad class of (Ising or spherical) spin systems are also understood through the \textbf{small subgraph conditioning} method.
This technique was introduced in \cite{robinson1992almost,robinson1994almost}, and has been applied in a wide variety of spin systems on random graphs \cite{cooper1996perfect,mossel2009hardness,kemkes2010chromatic,fabian2021ising,coja2026fluctuations,prodromidis2026distribution} and mean-field models \cite{abbe2022proof,bencs2025zeros}, including the aforementioned high temperature SK and symmetric Ising perceptron models.
The idea is to identify a reweighting random variable $A_N$, whose distribution we understand, such that $\bbE[(Z_N / A_N - 1)^2] = o_N(1)$.
There are by now mature heuristics to choose this $A_N$, which is usually a statistic of subgraphs of bounded or slowly growing size.
At high temperature $\theta < 1$ (but not $\theta = 1$), the CLT in Theorem~\ref{t:main} can also be derived using small subgraph conditioning; we prove this result using our Hopfield comparison approach to demonstrate the generality of this connection.

Closer to the present paper, extensive work has studied fluctuations of the Sherrington--Kirkpatrick \cite{sherrington1975solvable} free energy.
It is now known that the free energy has Gaussian fluctuations of order $1$ at high temperature $\beta < 1$ \cite{aizenman1987some,comets1995sherrington}, and fluctuations of order $N^{1/6}$ at low temperature $\beta > 1$ \cite{chatterjee2009disorder, openai2026low}.
At criticality $\beta = 1$, physics work of Aspelmeier \cite{aspelmeier2008free} predicted the free energy has variance $\fr16 \log N + O(1)$.
Following the line of work \cite{talagrand2011mean2, chen2019order, dey2026fluctuations, schertzer2026order} this was proved in \cite{du2026fluctuations} (see also \cite{chen2026moderate}) along with a Gaussian CLT for the free energy.

\subsection{Relation to critical SK model and proof discussion}
\label{ss:proof-discussion}

A key idea in \cite{du2026fluctuations} is to identify the free energy fluctuations of the critical SK model by comparing it to the critical spherical SK model.
They proved that $\bbE[(Z_N/Z_N^\sph - 1)^2] = o_N(1)$, where $Z_N$ and $Z_N^\sph$ are the partition functions of these models.
This can be understood as a critical version of the small subgraph conditioning method, where $Z_N^\sph$ is the explanatory variable.
Unlike earlier applications of small subgraph conditioning, which required the model to have free energy fluctuations of order $1$ (and therefore be at high temperature), this method works even when the fluctuations diverge.
As described in Remark~\ref{r:comparison-principle-general}, Theorem~\ref{t:perceptron-comparison} substantially generalizes this comparison principle.

At the same time, proving this more general comparison principle involves several new challenges.
In \cite{du2026fluctuations}, the exact identity $\bbE[Z_N/Z_N^\sph]=1$ meant that it sufficed to control the second moment of $Z_N/Z_N^\sph$.
No such identity is available for $\bbE[X_N^\ast]$, so we must show that $\bbE[X_N^\ast]$ and $\bbE[(X_N^\ast)^2]$ both converge to $1$.
Consequently, we need to show that the second moments of the reweighting factor $Z^\sph_{\hop,N}(\bG)^{-1}$ approximate each other under several different planted models; see Definition~\ref{d:planted-models} and \ref{i:overview-central} below.

More importantly, \cite[Lemma~3.3]{du2026fluctuations} used a log-concavity argument to sidestep many random matrix estimates.
This argument is no longer available in this paper's setting.
Instead, we develop the random matrix machinery needed to implement the critical reweighted moment method much more robustly.
In particular, we prove random matrix estimates that accurately control the most difficult term in the reweighted moment calculation, the second moment of $Z^\sph_{\hop,N}(\bG)^{-1}$ under a planted model; see Propositions~\ref{p:uniform-planted-inverse-moments}--\ref{p:bbp}.
We expect that our more robust method will also be useful for proving outstanding critical-window extensions of the earlier work on the critical SK model; see \cite[Subsection 1.4]{du2026overlap}.

\subsection*{Acknowledgments}

Part of this work was completed while BH and HD were visiting Peking University.
We are grateful to Elizabeth Collins-Woodfin, Fu-Hsuan Ho, Han Gia Le, Jason Prodromidis, and Qiang Wu for motivating conversations.
HD was partially supported by an NSF-Simons research collaboration grant (award number 2031883).
BH was supported by a Stanford Science Fellowship and NSF Mathematical Sciences Postdoctoral Fellowship. 
SG was supported by the Institute for Foundations of Data Science at Yale University.
ZL was supported in part by the NSFC Major Program (number 12595294) and the NSFC Key Program (number 12231002).

\subsection*{Statement on AI use}

The overall strategy and architecture of the proof was designed by the human authors.
This includes the plan of comparing to a spherical Hopfield model and formulation of Theorems~\ref{t:perceptron-comparison}--\ref{t:hopfield-clt}; the plan and phenomenology described in \S\ref{s:overview}; and the reduction in \S\ref{s:proof} from Theorem~\ref{t:perceptron-comparison} to the the inputs in \S\ref{ss:proof-inputs}.
These inputs were proved with AI assistance, as described below.

The proofs of Propositions~\ref{p:uniform-planted-inverse-moments}--\ref{p:bbp}, which occupy \S\ref{s:rmt}, were obtained by iterating with GPT-6 Astra.
Within these proofs, Lemmas~\ref{l:coarea-upper-bound}--\ref{l:gaussian-replacement-matrix} were proved by GPT-6 Astra.
Given these lemmas, the human authors realized Propositions~\ref{p:uniform-planted-inverse-moments}--\ref{p:bbp} could be reduced to:
\begin{itemize}
  \item A suitably sharp estimate on $\bbE |\det (z - \beta \bW)|^k$ where $z,k$ are fixed and $\bW$ is a Wishart matrix. This and related problems are well studied in the literature \cite{auffinger2013random,ben2023exponential}.
  \item Inward deviation estimates on extremal eigenvalues of Wishart matrices (possibly with $O(1)$ planted rows and columns).
\end{itemize} 
This reduction is carried out in Lemma~\ref{l:zeta-discretization}.
The reduced estimates are standard random matrix calculations, which were proved by GPT-6 Astra in Lemmas~\ref{l:characteristic-polynomial-moment} and \ref{l:wishart-inward-deviation}--\ref{l:wishart-inward-deviation-bbp}.

Lemmas~\ref{l:good-minor}--\ref{l:hopfield-overlap-concentration} were proved by GPT-6 Astra, and their proofs are presented in \S\ref{s:cov}. The human authors formulated Lemma~\ref{l:good-minor} using heuristics described in \S\ref{s:overview} and proposed the proof strategy. The proof steps were filled in by GPT-6 Astra by combining well-established random matrix estimates.
The human authors formulated Lemma~\ref{l:hopfield-overlap-concentration} after GPT-5.6 Sol found that an earlier version of our argument had an error for $s\ge 2$ (see Remark~\ref{r:hopfield-planting-domain}); its proof is by standard random matrix estimates, filled in by GPT-6 Astra.

Theorem~\ref{t:hopfield-clt} was proved by GPT-6 Astra, and its proof is presented in \S\ref{rewrite:a:hopfield-clt}.
This proof is a routine extension of the methods of \cite{landon2022free,collins2025fluctuations}.

All AI proofs in this paper were digested and rewritten by the human authors, and we take full responsibility for the results.

%\newpage

\section{Proof overview}
\label{s:overview}

In this section, we give an overview of the proof of Theorem~\ref{t:perceptron-comparison}.

\begin{itemize}
  \item In \S\ref{ss:overview-planting}, we present Lemma~\ref{l:reweight}, which provides algebraic identities relating moments of $X_N^\ast$, for $\ast \in \{\sph, \Is\}$, to expectations in three structurally similar random matrices.
  \item As a result of Lemma~\ref{l:reweight}, the target $\bbE[(X_N^\ast - 1)^2]$ in Theorem~\ref{t:perceptron-comparison} can be expressed as a difference of three similar integrals, over overlaps $q\in [-1,1]$, corresponding to these random matrices; see \eqref{e:2mt-expansion-sphere}.
  In \S\ref{ss:overview-comparison}, we describe a comparison strategy inspired by \cite{du2026fluctuations} to bound this difference.
  We will identify a central window $|q| \le q_N$ accounting for the main contribution to these integrals, and argue that the integrands approximate each other within this window, while the contribution from outside it is negligible.
  The cutoff $q_N$ depends on whether we are in the high-temperature or critical setting, and we give a heuristic explanation for how it is set.
\item In \S~\ref{s:overview-regimes}, we outline the proofs of the central comparison and tail estimates claimed in \S~\ref{ss:overview-comparison}, explaining in greater technical detail how the analysis is related to random matrix theory and differs between the high-temperature and critical regimes.
\end{itemize}

\subsection{Planted model identity for reweighted overlaps}
\label{ss:overview-planting}

In this subsection, we present Lemma~\ref{l:reweight}, which provides algebraic identities relating moments of $X_N^\sph$ and $X_N^\Is$ to expectations in the following three planted models.
We will use this in \S\ref{ss:overview-comparison} to control the quantity $\bbE[(X^\ast_N-1)^2]$ from Theorem~\ref{t:perceptron-comparison}. Denote $R(x,y)=N^{-1}\langle x,y\rangle$ for $x,y\in \mathbb{R}^N$. 
\begin{dfn}[Planted models]\label{d:planted-models}
  Consider any $q\in [-1,1]$.
  Let $\bbE_q$ denote expectation with respect to a centered jointly Gaussian pair $(Z_1,Z_2)$ with unit variances and correlation $q$, and define
  \beq\label{e:overview-planted-weights}
  \begin{aligned}
    \Phi_{\per}(q)
    &=\fr{\bbE_q[\phi(Z_1)\phi(Z_2)]}{m_\phi^2}, &
    \Phi_{\mix}(q)
    &=\fr{\bbE_q[\phi(Z_1)\exp(\beta Z_2^2/2)]}
    {m_\phi\sqrt{s}}, \\
    &&
    \Phi_{\hop}(q)
    &=\fr{\bbE_q \exp(\beta (Z_1^2 + Z_2^2) / 2)}{s}.
  \end{aligned}
  \eeq
  Let $\be_1,\be_2 \in \bbR^N$ denote the first two standard basis vectors and
  \balnn
    \bx^1_q &= \sqrt{\fr{N}{2}} \lt(\sqrt{1+q} \be_1 + \sqrt{1-q} \be_2\rt), & 
    \bx^2_q &= \sqrt{\fr{N}{2}} \lt(\sqrt{1+q} \be_1 - \sqrt{1-q} \be_2\rt),
  \ealnn
  be fixed vectors in $S_N$ with $R(\bx^1_q,\bx^2_q) = q$.
  Then let $\bbP_{\per,q}$, $\bbP_{\mix,q}$, $\bbP_{\hop,q}$ denote the probability measures on $\bbR^N$ with densities relative to $\bbP_0 = \cN(0,\bI_N)$ given by
  \beq\label{e:overview-planted-laws}
    \begin{aligned}
    \fr{\de\bbP_{\per,q}}{\de\bbP_0}(\bg)
    &=
    \fr{\phi(Z'_1)\phi(Z'_2)}{m_\phi^2\Phi_{\per}(q)}, &
    \fr{\de\bbP_{\mix,q}}{\de\bbP_0}(\bg)
    &=
    \fr{\phi(Z'_1)\exp(\beta (Z'_2)^2/2)}{m_\phi\sqrt{s}\,\Phi_{\mix}(q)}, \\
    &&
    \fr{\de\bbP_{\hop,q}}{\de\bbP_0}(\bg)
    &=
    \fr{\exp(\beta [(Z'_1)^2 + (Z'_2)^2]/2)}{s\Phi_{\hop}(q)},
    \end{aligned}
  \eeq
  where $Z'_1 = (\bx^1_q, \bg) / \sqrt{N}$ and $Z'_2 = (\bx^2_q, \bg) / \sqrt{N}$.
  For $\star \in \{\per, \mix, \hop\}$, we slightly abuse notation by letting $\bbP_{\star,q}$ also denote the law of $\bG \in \bbR^{M\times N}$ whose rows are drawn i.i.d.\ from $\bbP_{\star,q}$, and $\bbE_{\star,q}$ denote expectation with respect to this random matrix.
\end{dfn}
Note that each planted law $\bbP_{\star,q}$ modifies the distribution of $\bG$ only in the first two columns; the remaining $M \times (N-2)$ submatrix remains standard Gaussian and independent of the first two columns.
\begin{rmk}\label{r:planting-rotationally-invariant}
  The choice of $\bx^1_q,\bx^2_q$ in Definition~\ref{d:planted-models} is arbitrary among all $\bx^1_q,\bx^2_q \in S_N$ with $R(\bx^1_q,\bx^2_q) = q$, as in what follows we will only apply the planted expectations $\bbE_{\star,q}$ to functions of $\bG$ that are right-orthogonally invariant.
\end{rmk}
\begin{rmk}\label{r:hopfield-planting-domain}
  Condition \eqref{e:sub-quadratic-exponential-growth} ensures $\bbP_{\per,q}$ and $\bbP_{\mix,q}$ are well-defined for all $q\in [-1,1]$.
  However, when $\beta > 0$, the normalization $\Phi_{\hop}(q)$ of $\bbP_{\hop,q}$ is finite only when $\beta (1+|q|) < 1$.
  By \eqref{e:beta}, this holds if
  \beq\label{e:Iadm}
    q \in \cI_\adm(s) \equiv \begin{cases}
      [-1,1] & s < 2, \\
      (-\fr{1}{s-1}, \fr{1}{s-1}) & s \ge 2.
    \end{cases}
  \eeq
  The fact that $\bbP_{\hop,q}$ is sometimes not defined for $s\ge 2$ causes no essential difficulty, as in our proofs we only use $\bbP_{\hop,q}$ for $q$ in a neighborhood of zero.
  In this overview we mostly focus on the case $s<2$.
\end{rmk}
Let $\rho_N^\sph$ and $\rho_N^\Is$ denote the laws of $q = R(\bx^1,\bx^2)$ for independent $\bx^1,\bx^2 \sim \unif(S_N)$ and $\bx^1,\bx^2 \sim \unif(\Sigma_N)$, respectively.
Let $\la \cdot \ra_\bG$ denote Gibbs average with respect to the Gibbs measure associated to the Hopfield model \eqref{e:hopfield}, and $R_{1,2} = R(\bx^1,\bx^2)$ for i.i.d. $\bx^1,\bx^2$ sampled from this Gibbs measure, conditionally on $\bG$.
The following lemma follows by rotational invariance and a change-of-measure calculation.
\begin{lem}[Proved in \S\ref{ss:proof-deferred}]\label{l:reweight}
  For $\star \in \{\per, \mix\}$ and $q\in [-1,1]$, or $\star = \hop$ and $q \in \cI_\adm(s)$, define
  \beq\label{e:IstarN}
    I_{\star,N}(q) =
    \lt[s\Phi_{\star}(q)\rt]^M
    \bbE_{\star,q} \lt[
      Z^\sph_{\hop,N}(\bG)^{-2}
    \rt].
  \eeq
  For either $\ast \in \{\sph,\Is\}$,
  \baln\label{e:reweight-per-mix}
    \bbE[(X_N^\ast)^2]
    &= \int_{-1}^1
    I_{\per,N}(q) \,\rho_N^\ast(\de q), &
    \bbE[X_N^\ast]
    &= \int_{-1}^1 
    I_{\mix,N}(q) \,\rho_N^\sph(\de q).
  \ealn
  Also, for any bounded measurable $g : [-1,1] \to \bbR$ supported on a compact subset of $\cI_\adm(s)$,
  \baln
    \label{e:reweight-hop}
    \bbE[\la g(R_{1,2})\ra_\bG]
    &= \int_{\cI_\adm(s)}
    g(q) I_{\hop,N}(q) \,\rho_N^\sph(\de q).
  \ealn
\end{lem}

Before proceeding, we record two intuitions about these planted models that will be useful below.
First, the $\bbP_{\star,q}$ are designed to be spectrally similar at $q=0$, in the following sense.
Let $\bC \in \bbR^{M\times 2}$ be the first two columns of $\bG$, and form the empirical covariance matrix $\bW = \fr{1}{N} \bG^\top \bG$.
By rotational invariance, the conditional law of $\spec(\bW)$ given $\bC$ depends on $\bC$ only through the Gram matrix $\bK = \bC^\top \bC$.
As a consequence of Remark~\ref{r:beta-choice}, $\bK$ has expectation $\bbE_{\star,q=0}[\bK] = M s \cdot I_2$ under any of the planted laws, and further concentrates around this common expectation.
Thus we expect the spectra of $\bW$ to be similar under these planted laws.

Second, the criticality transition at $\theta=1$ coincides with a (shared) BBP-type phase transition in the planted models $\bbP_{\star,q=0}$.
Under these distributions, $\bW$ has bottom-right $(N-2)\times (N-2)$ minor sampled from the Wishart distribution, and two special rows and columns.
At $q=0$ and high temperature $\theta < 1$, the limiting spectrum of $\bW$ consists of a Marchenko--Pastur bulk and no outliers; as we increase $\theta$ beyond $1$, two outliers push out from the bulk.
These outliers lie above the bulk if $\beta > 0$ and below if $\beta < 0$.

\subsection{Comparison strategy}
\label{ss:overview-comparison}

We next explain our strategy for bounding the objective $\bbE[(X_N^\ast - 1)^2]$ from Theorem~\ref{t:perceptron-comparison}.
For simplicity, we mainly focus on $s<2$, $\ast = \sph$.
Taking $g = 1$ in Lemma~\ref{l:reweight} yields
\beq\label{e:2mt-expansion-sphere}
  \bbE[(X_N^\sph - 1)^2] = \int_{-1}^1
  \big( I_{\per,N}(q) - 2I_{\mix,N}(q) + I_{\hop,N}(q) \big)
  \,\rho_N^\sph(\de q).
\eeq
Define the threshold
\beq\label{e:qN}
  q_N =
  \begin{cases}
  N^{-1/2}\log N, & \theta<1,\\
  N^{-1/3}\log N, & \theta=1.
  \end{cases}
\eeq
We will estimate \eqref{e:2mt-expansion-sphere} using a comparison strategy inspired by \cite{du2026fluctuations}, consisting of two steps:
\begin{enumerate}[label=(C\arabic*)]
  \item\label{i:overview-central}
  In the central window $|q| \le q_N$, the integrands $I_{\star,N}$ approximate each other in the sense that
  \[
    \max_{\star \in \{\per, \mix\}} \int_{|q| \le q_N} \lt|{I_{\star,N}(q)}-{I_{\hop,N}(q)}\rt|\rho_N^\sph(\de q) = o_N(1).
  \]
  \item\label{i:overview-tail}
  The contributions from the tail $|q| > q_N$ are negligible, in the sense that
  \[
    \max_{\star \in \{\per, \hop\}} \int_{|q| > q_N} I_{\star,N}(q)\,\rho_N^\sph(\de q) = o_N(1).
  \]
\end{enumerate}
Recalling that $I_{\hop,N}(q)$ integrates to exactly $1$ under $\rho_N^\sph$ by \eqref{e:reweight-hop}, these together imply that \eqref{e:2mt-expansion-sphere} is $o_N(1)$.
(The tail contribution from $I_{\mix,N}$ to \eqref{e:2mt-expansion-sphere} is negative and can be discarded without argument.)

We describe the proofs of \ref{i:overview-central} and \ref{i:overview-tail} in more detail in \S\ref{s:overview-regimes}, after the following remarks.

\begin{rmk}\label{r:sidestep-explicit-calculation}
  Unlike usual applications of small subgraph conditioning, in our comparison strategy we never need to compute the $I_{\star,N}(q)$ themselves up to a multiplicative $1+o(1)$ factor.
  Indeed, the central window \ref{i:overview-central} only requires a relative estimate, while crude estimates suffice for the tail bound \ref{i:overview-tail}.
  We in fact do not know a way to directly compute the $I_{\star,N}(q)$ due to the complexity of the $Z^\sph_{\hop,N}(\bG)^{-2}$ term; being able to sidestep this computation is crucial to making our reweighted moment method work with the critical reweighting factor.
\end{rmk}
\begin{rmk}\label{r:overview-localization}
  The choice \eqref{e:qN} is due to the fact that at $\theta < 1$, the main contribution to \eqref{e:2mt-expansion-sphere} comes from $|q| \lesssim N^{-1/2}$, while at $\theta = 1$ it comes from $|q| \lesssim N^{-1/3}$.
  Indeed, Taylor expansion shows
  \[
    \Phi_{\star}(q)^M = \exp\lt(\fr{N}{2} \alpha (s-1)^2 q^2 + O(q^2 + Nq^4)\rt)
  \]
  for all $\star \in \{\per, \mix, \hop\}$; here our choice \eqref{e:beta} of $\beta$ makes these coincide to leading order.
  For $\theta < 1$, the leading term in $\Phi_{\star}(q)^M$ does not fully cancel the leading-order curvature in the density
  \[
    \rho_N^\sph(\de q) \,\,\propto\,\, \exp\lt(-\fr{N}{2} q^2 + O(q^2 + Nq^4)\rt)\,\de q.
  \]
  The remaining quadratic curvature around $0$ causes \eqref{e:2mt-expansion-sphere} to localize to $|q| \lesssim N^{-1/2}$.

  At criticality $\theta = 1$, the leading quadratic terms cancel, and the localization is instead caused by BBP soft-edge considerations.
  Recall from the discussion below Lemma~\ref{l:reweight} that in any of the planted models $\bbP_{\star,q=0}$, $\bW = \fr{1}{N} \bG^\top\bG$ has two exactly critical BBP spikes.
  As we vary $q$ away from $0$, one of these spikes becomes supercritical, and considerations similar to \cite[Remark 2.3]{du2026fluctuations}
  % \BH{new version}
  suggest that for some constant $c>0$,
  \begin{equation}\label{e:overview-reduction-to-q=0}
    \bbE_{\star,q} \lt[
      Z^\sph_{\hop,N}(\bG)^{-2}
    \rt]
    \le e^{-cN|q|^3} \cdot 
    \bbE_{\star,q=0} \lt[
      Z^\sph_{\hop,N}(\bG)^{-2}
    \rt].
  \end{equation}
  This causes \eqref{e:2mt-expansion-sphere} to localize to $|q| \lesssim N^{-1/3}$.
\end{rmk}
In the case $s\ge 2$, the only difference is that we take $g(q) = \bbone\{|q| \le q_0\}$ for a small positive $q_0 \in \cI_\adm(s)$.
We bound the large-overlap contribution $\bbE[\la \bbone\{|R_{1,2}| > q_0\} \ra_\bG]$ to \eqref{e:2mt-expansion-sphere} by a separate quenched argument (see Lemma~\ref{l:hopfield-overlap-concentration}).

For the Ising perceptron, the difference is that (only) the $I_{\per,N}(q)$ term in \eqref{e:2mt-expansion-sphere} is now integrated over $\rho_N^\Is$ instead of $\rho_N^\sph$.
This can be addressed by a discrete-to-continuous approximation argument (see Lemma~\ref{l:ising-grid} and the \hyperlink{p:t:perceptron-ising-high-temp}{proof of Theorem~\ref*{t:perceptron-comparison}\ref*{i:ising-comparison}}), showing that the central-window integrals of $I_{\per,N}(q)$ over these two measures are close.

\subsection{Proof of the two inputs: high-temperature and critical regimes}
\label{s:overview-regimes}

We now explain how to prove \ref{i:overview-central} and \ref{i:overview-tail}.
The following random variable $\zeta$ will be a key object in both proofs.
For disorder $\bG$, let $\bW = \bG^\top \bG / N$, and let $\zeta = \zeta(\bG) > \lambda_{\max}(\beta \bW)$ be the unique solution to
\beq\label{e:overview-saddle}
  \Tr (\zeta I_N - \beta \bW)^{-1} = N.
\eeq
The meaning of $\zeta$ can be explained through the following physical intuition.
Consider the Gaussian $\bY \sim \cN(\bzero, (\zeta I_N - \beta \bW)^{-1})$.
By the definition \eqref{e:overview-saddle} of $\zeta$, $\bbE \|\bY\|^2 = N$, and therefore $\|\bY\|$ concentrates around $\sqrt{N}$.
At the same time, the law of $\bY$ conditional on $\|\bY\| = \sqrt{N}$ is precisely the Gibbs measure of the Hopfield model \eqref{e:hopfield}.
We thus expect this Gibbs measure to resemble the law of $\bY$.
In particular, this heuristic suggests the Gibbs measure's covariance satisfies
\beq\label{e:overview-covariance}
  \la \bx\bx^\top \ra_\bG
  \approx (\zeta I_N - \beta \bW)^{-1}.
\eeq
This connection between the Hopfield model and $\bY$ also suggests that to leading exponential order,
\beq\label{e:overview-determinant}
  Z^\sph_{\hop,N}(\bG)
  \approx \fr{(2\pi)^{N/2}}{|S_N|} e^{N\zeta /2} \det(\zeta I_N - \beta \bW)^{-1/2}
  \approx e^{N(\zeta-1) /2} \det(\zeta I_N - \beta \bW)^{-1/2},
\eeq
where $|S_N|$ is the Hausdorff measure of $S_N \subseteq \bbR^N$.
Over the course of our proof, we will obtain rigorous versions of \eqref{e:overview-covariance}, \eqref{e:overview-determinant} (see \eqref{e:overview-covariance-rigorous}, \eqref{e:overview-determinant-rigorous} below), which reduce the difficult objects $\la \bx\bx^\top \ra_\bG$ and $Z^\sph_{\hop,N}(\bG)$ to tractable random matrix calculations.
We now outline the proofs of \ref{i:overview-central} and \ref{i:overview-tail}.

\paragraph{\textbf{Central comparison.}}
To prove \ref{i:overview-central}, we will show that the $I_{\star,N}(q)$ for $\star \in \{\per, \mix, \hop\}$ approximate each other to a multiplicative $1+o_N(1)$ factor uniformly for $|q| \le q_N$.
A direct Taylor expansion shows that the $\Phi_{\star}(q)^M$ in \eqref{e:IstarN} approximate each other in this sense; our main task is to show that the $\bbE_{\star,q}[Z^\sph_{\hop,N}(\bG)^{-2}]$ do as well.

Write $\bG = [\bC, \bA_0]$, where we recall $\bC \in \bbR^{M\times 2}$ consists of the first two columns of $\bG$ and (under any of the $\bbP_{\star,q}$) $\bA_0$ consists of i.i.d. standard Gaussians independent of $\bC$.
As noted above, the law of $Z^\sph_{\hop,N}(\bG)$ conditional on $\bC$ depends on $\bC$ only through $\bK = \bC^\top \bC$.
Thus, for fixed isometry $\bE \in \bbR^{M\times 2}$ and
\[
  f_\bK(\bA_0) = Z^\sph_{\hop,N} ([\bE \bK^{1/2}, \bA_0])^{-2},
\]
it suffices to show $\bbE_{\star,q}[f_\bK(\bA_0)]$ for $\star \in \{\per, \mix, \hop\}$ approximate each other to multiplicative $1+o_N(1)$, uniformly for $|q| \le q_N$.

We achieve this by showing that $\log f_\bK(\bA_0)$ is stable with overwhelming probability.
By standard concentration inequalities, under any $\bbP_{\star,q}$, $\bK$ concentrates in a region of the positive semidefinite cone $\bbS_+^2$ of size $N^{1/2 + o(1)}$, which is independent of $\star \in \{\per, \mix, \hop\}$.
We show (Lemma~\ref{l:good-minor}) that with overwhelming probability over $\bA_0$,
\beq\label{e:overview-stability}
  \|\nabla_\bK \log f_\bK(\bA_0)\| \le N^{-2/3 + o(1)}
\eeq
for all $\bK$ in this region.
Thus $f_{\bK}(A_0)=(1+o(1))f_{\bK'}(\bA_0)$ for typical $\bA_0$, and typical $\bK, \bK'$ from any two planted laws $\bbP_{\star,q}$; this is the stability mechanism underlying \ref{i:overview-central}.

The proof of Lemma~\ref{l:good-minor} occupies most of \S\ref{s:cov}, and we only comment on it briefly here.
The derivative in \eqref{e:overview-stability} can be expressed in terms of compressions $\bM^\top \la \bx\bx^\top \ra_\bG \bM$ of the covariance matrix where $\bM \in \bbR^{N \times 2}$ (Lemma~\ref{l:gram-response}).
We are able to control this matrix through a rigorous version of the heuristic \eqref{e:overview-covariance} (Lemma~\ref{l:hopfield-gibbs-covariance}): for a constant $C$,
\beq\label{e:overview-covariance-rigorous}
  \la \bx\bx^\top \ra_\bG
  \preceq C (\zeta I_N - \beta \bW)^{-1}.
\eeq
The relevant compression of the right-hand side of \eqref{e:overview-covariance-rigorous} can be controlled by random-matrix estimates, yielding \eqref{e:overview-stability}.

\paragraph{\textbf{Tail control.}}
Recall that
\beq\label{e:overview-Irho}
  I_{\star,N}(q)\rho_N^\sph(\de q)
  = [s\Phi_\star(q)]^M
  \bbE_{\star,q}[
    Z_{\hop,N}^{\sph}(\bG)^{-2}
  ]
  \rho_N^\sph(\de q).
\eeq
To show \ref{i:overview-tail}, we aim to show that the contribution from $|q| > q_N$ is $o_N(1)$.
We control the $Z_{\hop,N}^{\sph}(\bG)^{-2}$ term by a rigorous version of \eqref{e:overview-determinant} (Lemma~\ref{l:coarea-upper-bound}): for a constant $C>0$,
\beq\label{e:overview-determinant-rigorous}
  Z_{\hop,N}^{\sph}(\bG)^{-2}
  \le CN e^{N(1-\zeta)} \det(\zeta I_N - \beta \bW).
\eeq
The expectation of the right-hand side of \eqref{e:overview-determinant-rigorous} can be bounded by random matrix theory, with the loss of a harmless polynomial factor.
This shows (Propositions~\ref{p:uniform-planted-inverse-moments}--\ref{p:bbp}) that for constants $C,c$,
\beq\label{e:overview-moment-bound}
  \bbE_{\star,q}[
    Z_{\hop,N}^{\sph}(\bG)^{-2}
  ]
  \le \begin{cases}
    N^C s^{-M} & \theta < 1, \\
    \exp(-cN|q|^3) N^C s^{-M} & \theta = 1,
  \end{cases}
\eeq
where the factor $\exp(-cN|q|^3)$ arises by the mechanism in Remark~\ref{r:overview-localization}.
We now substitute these estimates into \eqref{e:overview-Irho}.
For $\theta < 1$, \eqref{e:overview-moment-bound} gives
\[
  I_{\star,N}(q)\rho_N^\sph(\de q)\le N^C\Phi_{\star}(q)^M\rho_N^\sph(\de q)\le N^C\exp(-cNq^2)\de q\,,
\]
where the last inequality follows from the discussion in Remark~\ref{r:overview-localization}.
At $\theta = 1$, \eqref{e:overview-moment-bound} gives
\[
  I_{\star,N}(q)\rho_N^\sph(\de q)\le N^C\exp(-cN|q|^3)\Phi_{\star}(q)^M\rho_N^\sph(\de q)\le N^C\exp(-cN|q|^3)\de q\,.
\]
In either case, the contribution from $|q| > q_N$ is $o_N(1)$, and the factor $N^C$ is immaterial.

\subsection*{Organization of the paper}

The rest of the paper is organized as follows.
\begin{itemize}
  \item The proof of Theorem~\ref{t:perceptron-comparison} occupies \S\ref{s:proof}--\ref{s:cov}.
  \S\ref{s:proof} formalizes the overall strategy of the proof, and proves this theorem assuming several random matrix inputs (Propositions~\ref{p:uniform-planted-inverse-moments}--\ref{p:bbp} and Lemmas~\ref{l:good-minor}--\ref{l:hopfield-overlap-concentration}).
  \item \S\ref{s:rmt} proves Propositions~\ref{p:uniform-planted-inverse-moments}--\ref{p:bbp}, which provide the upper bound \eqref{e:overview-moment-bound} on planted moments of $Z_{\hop,N}^{\sph}(\bG)^{-2}$.
  \item \S\ref{s:cov} proves Lemmas~\ref{l:good-minor}--\ref{l:hopfield-overlap-concentration}, which provide the high-probability stability bound \eqref{e:overview-stability} and a related bound that $R_{1,2}$ concentrates near $0$ under $\bbE \la \cdot \ra_\bG$.
  \item \S\ref{rewrite:a:hopfield-clt} proves Theorem~\ref{t:hopfield-clt} following the approach of \cite{landon2022free,collins2025fluctuations}.
\end{itemize}

%\newpage

\section{Proof of Hopfield comparison theorem}
\label{s:proof}

The purpose of this section is to record the overall architecture of the proof of Theorem~\ref{t:perceptron-comparison}, by proving this theorem assuming several deferred inputs.
This section is organized as follows.
\begin{itemize}
  \item In \S\ref{ss:proof-inputs}, we state the inputs that the proof will require.
  They belong to two main categories: some have routine or short, ad-hoc proofs, presented in \S\ref{ss:proof-deferred} immediately after the proof of Theorem~\ref{t:perceptron-comparison}.
  Others rely on more substantive random-matrix inputs, and their proofs are deferred to \S\ref{s:rmt}--\ref{s:cov}.
  \item The proof of Theorem~\ref{t:perceptron-comparison} occupies \S\ref{ss:proof-assembly-common}--\ref{ss:proof-assembly-critical}.
  The proof differs slightly between the high temperature $\theta < 1$ and critical $\theta = 1$ cases, and is organized as follows.
  \item In \S\ref{ss:proof-assembly-common} we give the parts of the proof that are shared by the two cases.
  The most important output of this section is Proposition~\ref{p:central-comparison}, which compares integrals of the $I_{\star,N}(q)$ in the central window.
  \item In \S\ref{ss:proof-assembly-high-temperature}, we complete the proof of Theorem~\ref{t:perceptron-comparison} for $\theta<1$.
  \item In \S\ref{ss:proof-assembly-critical}, we complete the proof for $\theta=1$.
  This case involves an additional random matrix input Proposition~\ref{p:bbp}, which captures the BBP soft-edge effect on the planted Hopfield moment described in Remark~\ref{r:overview-localization}.
  \item In \S\ref{ss:proof-deferred}, we prove the elementary inputs stated in \S\ref{ss:proof-inputs}.
\end{itemize}

Throughout this section we assume $\theta \le 1$.
We do not assume $\theta < 1$ or $\theta = 1$ unless otherwise stated.
We also do not assume either part of Assumption~\ref{a:2mt}, except the common implication $\theta \le 1$, unless otherwise stated.

\subsection{Inputs to the proof}
\label{ss:proof-inputs}

In this subsection, we collect several inputs to the proof of Theorem~\ref{t:perceptron-comparison}.
Of these, Lemmas~\ref{l:good-gram} and \ref{l:ising-grid} have short, elementary proofs presented in \S\ref{ss:proof-deferred}, while Propositions~\ref{p:uniform-planted-inverse-moments}--\ref{p:bbp} and Lemmas~\ref{l:good-minor}--\ref{l:hopfield-overlap-concentration} are proved in future sections using random matrix theory.

Due to the form of the integrand \eqref{e:IstarN} appearing in Lemma~\ref{l:reweight}, we will often need to upper bound expectations of $Z^\sph_{\hop,N}(\bG)^{-2}$ under the planted measures $\bbP_{\star,q}$.
The next proposition provides an upper bound within a polynomial factor of the natural exponential scale $s^{-kM}$.
\begin{ppn}[Proved in \S\ref{s:rmt}]
\label{p:uniform-planted-inverse-moments}
  There exists $C < \infty$ depending only on $\phi,\alpha$, such that for all $k\in \{1,2\}$, and $\star \in \{\per, \mix\}$ and $q\in [-1,1]$,
  \[
    \bbE_{\star,q}\lt[
      Z^\sph_{\hop,N}(\bG)^{-{2k}}
    \rt]
    \le CN^{4+7k} s^{-kM}.
  \]
  The same estimate holds for $\star = \hop$ and all $q\in J$, for any compact $J\subseteq \cI_\adm(s)$, where $C$ may also depend on $J$ in addition to $\phi,\alpha$.
\end{ppn}
In the critical case, we will need the following stronger estimate, which captures the decay speed in $q$ heuristically predicted in Remark~\ref{r:overview-localization}. 
The polynomial factor $CN^{11}$ will be harmless.
\begin{ppn}[Proved in \S\ref{s:rmt}]\label{p:bbp}
  Suppose $\theta = 1$, so that $q_N = N^{-1/3} \log N$ as in \eqref{e:qN}.
  There exist $c,C,q_0 > 0$ depending only on $\phi,\alpha$ such that for all $\star \in \{\per, \mix, \hop\}$ and all $|q| \in [q_N/2, q_0]$,
  \[
    \bbE_{\star,q}\lt[
      Z^\sph_{\hop,N}(\bG)^{-2}
    \rt] \le CN^{11} e^{-cN|q|^3} s^{-M}.
  \]
\end{ppn}
For the next two lemmas, we write $\bG \in \bbR^{M\times N}$ as $\bG = [\bC,\bA_0]$, where $\bC \in \bbR^{M\times 2}$ and $\bA_0 \in \bbR^{M\times (N-2)}$, and $\bK = \bC^\top\bC \in \bbS_+^2$, where $\bbS_+^2$ is the cone of $2\times 2$ positive semidefinite matrices.
Under any of the planted laws from Definition~\ref{d:planted-models}, $\bA_0$ is a standard Gaussian matrix independent of $\bC$.
These lemmas imply a stability estimate (Lemma~\ref{l:central-stability}) for the planted expectations of the $Z^\sph_{\hop,N}(\bG)^{-2}$ in Lemma~\ref{l:reweight}, when $q$ is in the central window.
This is the main input to the comparison of the integrals of $I_{\star,N}(q)$ in the central window (Proposition~\ref{p:central-comparison}).
Let
\beq\label{e:bSigma-hop}
  \bSigma_{\hop,q} = \bbE_{\hop,q}[\bK] = M \begin{pmatrix}
    (1 - \beta(1+q))^{-1} & 0 \\ 
    0 & (1 - \beta(1-q))^{-1} 
  \end{pmatrix}.
\eeq
Fix $\eps_0 \in (0, \fr{1}{12})$, and let $\oq_N = N^{-1/3}\log N$ be the larger of the two possible values of $q_N$ defined in \eqref{e:qN}.
For each $|q| \le \oq_N$, define
\balnn
  \cK_N(q) &= \lt\{
    \bK \in \bbS_+^2 : 
    \|\bK - \bSigma_{\hop,q}\|_F \le N^{2/3 - 2\eps_0}
  \rt\}, & 
  \cK_N &= \bigcup_{|q| \le \oq_N} \cK_N(q).
\ealnn
\begin{lem}[Proved in \S\ref{ss:proof-deferred}]
\label{l:good-gram}
  There exist $c,C>0$ depending only on $\phi,\alpha,\eps_0$ such that, uniformly for $|q|\le \oq_N$ and $\star \in \{\per, \mix, \hop\}$,
  \[
    \bbP_{\star,q}(\bK \not\in \cK_N(q)) \le C\exp(-cN^{1/3 - 4\eps_0}).
  \]
\end{lem}
\begin{lem}[Proved in \S\ref{s:cov}]
  \label{l:good-minor}
  Let $\bE \in \bbR^{M\times 2}$ be any isometry, so that $\bE^\top \bE = I_2$.
  Let $\bA_0 \in \bbR^{M\times (N-2)}$ have i.i.d. standard Gaussian entries.
  For $\bK \in \bbS_+^2$ and $C_0 > 0$, write
  \balnn
    f_\bK(\bA_0) &= Z^\sph_{\hop,N}([\bE \bK^{1/2}, \bA_0])^{-2}, & 
    \cG_\mnr(C_0) &= \lt\{
      \bA_0 :
      \sup_{\bK \in \cK_N} \|\nabla_{\bK}\log f_{\bK}(\bA_0)\|_F
      \le \fr{C_0}{N^{2/3 - \eps_0}}
    \rt\}.
  \ealnn
  There exists $C_0 < \infty$ depending only on $\phi, \alpha, \eps_0$ such that, for every $D>0$ and all sufficiently large $N$,
  \[
    \bbP_{\bA_0}\bigl(\cG_\mnr(C_0)^c\bigr)\le N^{-D}.
  \]
\end{lem}

Lemma~\ref{l:hopfield-overlap-concentration} will be used to control overlaps of the Hopfield model larger than a small constant, since as explained in Remark~\ref{r:hopfield-planting-domain} it may not be possible to study them by planting.
Finally, Lemma~\ref{l:ising-grid} relates $\rho^\Is_N$ to $\rho^\sph_N$ in the central window, which enables our comparison strategy for the Ising perceptron.

\begin{lem}[Proved in \S\ref{s:cov}]
  \label{l:hopfield-overlap-concentration}
For any $N$-independent $q_0>0$,
  \[
    \bbE\la\bbone\{|R_{1,2}|>q_0\}\ra_{\bG} = o_N(1).
  \]
\end{lem}
\begin{lem}[Proved in \S\ref{ss:proof-deferred}]
\label{l:ising-grid}
  Let the support of $\rho^\Is_N$ be denoted
  \beq\label{e:cQN}
    \cQ_N = \{-1,-1+\textstyle \fr2N,\ldots,1\}.
  \eeq
  Uniformly over $q \in \cQ_N \cap [-\oq_N, \oq_N]$, we have
  \[
    \rho_N^\Is(\{q\})
    =\lt(1+O(N^{-1}+q^2+Nq^4)\rt)
    \rho_N^\sph([q - 1/N, q + 1/N]).
  \]
\end{lem}

\subsection{Completion of the proof: common ingredients}
\label{ss:proof-assembly-common}

We will next work towards the proof of Theorem~\ref{t:perceptron-comparison}, assuming the inputs from \S\ref{ss:proof-inputs}.
In this subsection, we give several proof steps that are shared between the high temperature and critical cases.

We first work towards Proposition~\ref{p:central-comparison}, which shows that the various planted models approximate each other in the central window.
We first show this approximation for the inverse moments in the $I_{\star,N}(q)$.

\begin{lem}% [Stability in central window]
\label{l:central-stability}
  There exists $C$ depending only on $\phi,\alpha$ such that for all $\star \in \{\per, \mix, \hop\}$, $q,q' \in [-\oq_N, \oq_N]$ with $|q-q'| \le 1/N$, and any $D > 0$,
  \[
    \lt| 
      \bbE_{\star,q'} \lt[
        Z^\sph_{\hop,N}(\bG)^{-2}
      \rt]
      - \bbE_{\hop,q} \lt[
        Z^\sph_{\hop,N}(\bG)^{-2}
      \rt]
    \rt|
    \le CN^{-\eps_0} \bbE_{\hop,q} \lt[
      Z^\sph_{\hop,N}(\bG)^{-2}
    \rt] + CN^{-D} \cdot s^{-M}.
  \]
\end{lem}
\begin{proof}
  Fix an isometry $\bE \in \bbR^{M\times 2}$, sample $\bA_0 \in \bbR^{M\times (N-2)}$ with i.i.d. standard Gaussian entries, and independently sample $\bK, \bK'$ from their marginals under $\bbP_{\hop,q}$ and $\bbP_{\star,q'}$.
  Form the coupled matrices
  \balnn
    \bG &= [\bE \bK^{1/2}, \bA_0], &
    \bG' &= [\bE (\bK')^{1/2}, \bA_0]
  \ealnn
  and let $\bbE$ denote expectation with respect to the coupled law.
  By rotational invariance, the marginal distributions of $Z^\sph_{\hop,N}(\bG)$ and $Z^\sph_{\hop,N}(\bG')$ under this law are the same as if $\bG, \bG'$ were samples from $\bbP_{\hop,q}$ and $\bbP_{\star,q'}$.
  Thus, for $f_\bK$ as in Lemma~\ref{l:good-minor},
  \[
    \lt|\bbE_{\star,q'} \lt[
      Z^\sph_{\hop,N}(\bG)^{-2}
    \rt]
    - \bbE_{\hop,q} \lt[
      Z^\sph_{\hop,N}(\bG)^{-2}
    \rt]\rt|
    \le \bbE \lt[
      \lt|f_{\bK'}(\bA_0) - f_{\bK}(\bA_0)\rt|
    \rt].
  \]
  Let $C_0$ be given by Lemma~\ref{l:good-minor} and define the event
  \[
    \cE = \{\bK \in \cK_N(q)\} \cap \{\bK' \in \cK_N(q')\} \cap \{\bA_0 \in \cG_\mnr(C_0)\}.
  \]
  Recall $\bSigma_{\hop,q}$ defined in \eqref{e:bSigma-hop}.
  On the event $\cE$, the curve 
  \balnn
    \{t \bK + (1-t) \bSigma_{\hop,q} : 0\le t \le 1\}
    &\cup \{\bSigma_{\hop,t} : \min(q,q') \le t \le \max(q,q') \} \\
    &\cup \{t \bSigma_{\hop,q'} + (1-t) \bK' : 0\le t\le 1\}
  \ealnn
  lies in $\cK_N$, and has length at most $C(N^{2/3-2\eps_0} + N|q-q'|)$ for some $C = C(\phi,\alpha)$.
  Integrating the bound on $\nabla \log f_\bK$ from the definition of $\cG_\mnr(C_0)$ shows that on $\cE$,
  \[
    \lt|\log f_\bK(\bA_0)-\log f_{\bK'}(\bA_0)\rt| \le \fr{C_0}{N^{2/3-\eps_0}} C(N^{2/3-2\eps_0} + N|q-q'|)
    \le CN^{-\eps_0}
  \]
  after adjusting $C$.
  Lemmas~\ref{l:good-gram} and \ref{l:good-minor} show that
  \[
    \bbP(\cE^c)
    \le C\exp\big(-cN^{1/3-4\eps_0}\big)+N^{-D'}
  \]
  for any $D' > 0$.
  By Cauchy--Schwarz,
  \balnn
    \bbE \lt[
      \lt|f_{\bK'}(\bA_0) - f_{\bK}(\bA_0)\rt|
    \rt]
    &\le \bbE \lt[
      \bbone\{\cE\} \lt|f_{\bK'}(\bA_0) - f_{\bK}(\bA_0)\rt|
    \rt] + \bbE \lt[
      \bbone\{\cE^c\} \lt|f_{\bK'}(\bA_0) - f_{\bK}(\bA_0)\rt|
    \rt] \\
    &\le \lt[\exp(CN^{-\eps_0}) - 1\rt] \bbE\lt[f_{\bK}(\bA_0)\rt] 
    + 2^{1/2} \bbP(\cE^c)^{1/2} \bbE\lt[f_{\bK'}(\bA_0)^2 + f_{\bK}(\bA_0)^2\rt]^{1/2}.
  \ealnn
  Finally note that
  \balnn
    \bbE\lt[f_{\bK}(\bA_0)\rt] &= \bbE_{\hop,q} \lt[Z^\sph_{\hop,N}(\bG)^{-2}\rt], && \\
    \bbE\lt[f_{\bK}(\bA_0)^2\rt] &= \bbE_{\hop,q} \lt[Z^\sph_{\hop,N}(\bG)^{-4}\rt], & 
    \bbE\lt[f_{\bK'}(\bA_0)^2\rt] &= \bbE_{\star,q'} \lt[Z^\sph_{\hop,N}(\bG)^{-4}\rt],
  \ealnn
  and the last two moments are bounded by $N^C s^{-2M}$ by Proposition~\ref{p:uniform-planted-inverse-moments}.
  The conclusion follows by taking $D'$ sufficiently large depending on $C,D$.
\end{proof}

\begin{ppn}% [Central comparison]
\label{p:central-comparison}
  Uniformly over $\star\in\{\per, \mix, \hop\}$, Borel maps $\tau:[-\oq_N,\oq_N]\to[-\oq_N,\oq_N]$ satisfying $\sup_q|\tau(q)-q| \le 1/N$, and measurable $A\subseteq [-\oq_N, \oq_N]$,
  \[
    \lt|
      \int_A I_{\star,N}(\tau(q))\,\rho_N^\sph(\de q)
      -\int_A I_{\hop,N}(q)\,\rho_N^\sph(\de q)
    \rt|=o_N(1).
  \]
\end{ppn}
\begin{rmk}
  The function $\tau$ enables comparisons between spherical and Ising overlap distributions, since the latter is supported on the discrete grid $\cQ_N$.
  In such applications, $\tau$ will round $q\in [-\oq_N, \oq_N]$ to the nearest element of $\cQ_N$.
\end{rmk}
\begin{proof}
  For all $\star \in \{\per, \mix, \hop\}$, a direct Taylor expansion shows that in a fixed neighborhood of $0$,
  \beq\label{e:Phi-taylor}
    \Phi_\star(q)^M = \exp\lt(\fr{N}{2} \cdot \alpha (s-1)^2 q^2 + O(q^2 + Nq^4) \rt).
  \eeq
  Since $N\oq_N^4=o_N(1)$, we have uniformly over $\star\in\{\per, \mix, \hop\}$ and $|q|,|q'| \le \oq_N$, $|q-q'|\le1/N$ that
  \beq\label{e:Phi-comparison}
    \Phi_\star(q')^M = (1+o_N(1))\Phi_\hop(q)^M.
  \eeq
  Combining the last display with Lemma~\ref{l:central-stability}, we infer that for every $D>0$, there is a deterministic $\delta_N = o_N(1)$ such that
  \beq\label{e:central-spherical-comparison-1}
    \lt|I_{\star,N}(q') - I_{\hop,N}(q)\rt|
    \le \delta_N I_{\hop,N}(q)
    + CN^{-D}\Phi_\star(q')^M,
  \eeq
  uniformly over $\star,q,q'$ as above.
  In particular, taking $q' = \tau(q)$ and integrating over $A$ yields
  \baln\nonumber
    &\lt|
      \int_A I_{\star,N}(\tau(q))\,\rho_N^\sph(\de q)
      -\int_A I_{\hop,N}(q)\,\rho_N^\sph(\de q)
    \rt| \\
    \label{e:central-spherical-comparison-2}
    &\le \delta_N \int_A I_{\hop,N}(q)\,\rho_N^\sph(\de q)
    + CN^{-D} \int_A \Phi_\star(\tau(q))^M \,\rho_N^\sph(\de q).
  \ealn
  The first integral in \eqref{e:central-spherical-comparison-2} can be estimated using \eqref{e:reweight-hop} from Lemma~\ref{l:reweight}:
  \[
    \int_A I_{\hop,N}(q)\,\rho_N^\sph(\de q)
    = \bbE \la \bbone\{R_{1,2} \in A\} \ra_\bG
    \le 1.
  \]
  Note that the density of $\rho_N^\sph$ satisfies $\rho_N^\sph(\de q) \le C\sqrt{N} \exp(-\fr{N}{2}q^2 + O(q^2 + Nq^4))\de q$ in a fixed neighborhood of $0$.
  Combining with \eqref{e:Phi-taylor} and \eqref{e:Phi-comparison} shows
  \[
    \Phi_\star(\tau(q))^M \,\rho_N^\sph(\de q)
    \le C\sqrt{N} \exp\lt(
      - \fr{N}{2} \lt(1 - \alpha (s-1)^2\rt) q^2 
      + O(q^2 + Nq^4)
    \rt) \,\de q.
  \]
  The assumptions $\theta = \alpha (s-1)^2 \le 1$ and $|q| \le \oq_N$ imply this is bounded by $C\sqrt{N}$ after adjusting $C$.
  Thus
  \[
    \int_A \Phi_\star(\tau(q))^M \,\rho_N^\sph(\de q)
    \le \int_{-\oq_N}^{\oq_N} \Phi_\star(\tau(q))^M \,\rho_N^\sph(\de q)
    \le C\sqrt{N} \oq_N,
  \]
  and taking $D$ sufficiently large makes \eqref{e:central-spherical-comparison-2} bounded by $o_N(1)$.
\end{proof}

Finally, the next lemma controls large-overlap contributions to the integral formula \eqref{e:reweight-per-mix} for $\bbE[(X^\ast_N)^2]$, directly through Assumption~\ref{a:2mt}.

\begin{lem}%[Perceptron large-overlap tails]
  \label{l:perceptron-large-overlap}
  For any $N$-independent $q_0>0$ and $\ast \in \{\sph,\Is\}$, under the corresponding part of Assumption~\ref{a:2mt},
  \[
    \int_{|q| > q_0} I_{\per,N}(q) \rho_N^\ast(\de q) = o_N(1).
  \]
\end{lem}
\begin{proof}[Proof of Lemma~\ref{l:perceptron-large-overlap}]
  By Proposition~\ref{p:uniform-planted-inverse-moments} with $k=1$,
  \baln\nonumber
    \int_{|q| > q_0} I_{\per,N}(q) \rho_N^\ast(\de q)
    &= \int_{|q| > q_0} 
    \lt[s\Phi_{\per}(q)\rt]^M
    \bbE_{\per,q} \lt[
      Z^\sph_{\hop,N}(\bG)^{-2}
    \rt] 
    \,\rho_N^\ast(\de q) \\
    \label{e:lem-perceptron-large-overlap}
    &\le N^C \int_{|q| > q_0} 
    \Phi_{\per}(q)^M
    \,\rho_N^\ast(\de q).
  \ealn
  Recall the functions $\cF_\ast$ from \eqref{e:cF-sph}--\eqref{e:cF-Is}.
  Since Assumption~\ref{a:2mt} provides that $\cF_\ast$ is maximized uniquely at $0$, and $\cF_\sph$ tends to $-\infty$ as $|q| \to 1$, a standard compactness argument gives that
  \beq\label{e:cF-negative}
    \sup_{|q| > q_0} \cF_\ast(q) \le -c
  \eeq
  for a constant $c>0$.
  If $\ast = \sph$, the integrand in \eqref{e:lem-perceptron-large-overlap} has exponential order
  \[
    \Phi_{\per}(q)^M
    \,\rho_N^\sph(\de q)
    = \exp(N\cF_\sph(q) + o(N)) \,\de q,
  \]
  and if $\ast = \Is$, then for $q$ in the support $\cQ_N$ of $\rho_N^\Is$,
  \[
    \Phi_{\per}(q)^M \,\rho_N^\Is(\{q\})
    = \exp(N\cF_\Is(q) + o(N)).
  \]
  In either case the result follows from \eqref{e:cF-negative}.
\end{proof}

\subsection{Completion of the proof: high temperature}
\label{ss:proof-assembly-high-temperature}

In this subsection, we assume $\theta < 1$ and complete the proof of Theorem~\ref{t:perceptron-comparison} in this case.
We accordingly set $q_N = N^{-1/2} \log N$ as in \eqref{e:qN}.

\begin{lem}%[High-temperature perceptron and Hopfield localization]
  \label{l:high-temperature-localization}
  Suppose $\theta < 1$.
  There exists a sufficiently small constant $q_0 > 0$, depending on $\phi, \alpha$, such that for any $(\star,\ast) \in \{(\per,\Is), (\per, \sph), (\hop, \sph)\}$,
  \[
    \int_{q_N - \fr2N < |q| \le q_0}
    I_{\star,N}(q) \,\rho_N^\ast(\de q)
    = o_N(1).
  \]
\end{lem}
\begin{proof}
  The proof is similar to that of Lemma~\ref{l:perceptron-large-overlap}.
  By Proposition~\ref{p:uniform-planted-inverse-moments} with $k=1$,
  \balnn
    \int_{q_N - \fr2N < |q| \le q_0} I_{\star,N}(q) \rho_N^\ast(\de q)
    &= \int_{q_N - \fr2N < |q| \le q_0}
    \lt[s\Phi_{\star}(q)\rt]^M
    \bbE_{\star,q} \lt[
      Z^\sph_{\hop,N}(\bG)^{-2}
    \rt] 
    \,\rho_N^\ast(\de q) \\
    &\le N^C \int_{q_N - \fr2N < |q| \le q_0}
    \Phi_{\star}(q)^M
    \,\rho_N^\ast(\de q).
  \ealnn
  The density of $\rho_N^\sph$ and mass of $\rho_N^\Is$ are both bounded by $\exp(-Nq^2/2 + O(q^2 + Nq^4))$, up to polynomial factors in $N$.
  Recall $\Phi_\star(q)^M$ has Taylor expansion \eqref{e:Phi-taylor} around $q=0$.
  After adjusting $C$, the last display is bounded by
  \beq\label{e:high-temperature-localization}
    N^C \sup_{q_N - \fr2N < |q| \le q_0}
    \exp\lt(
      -\fr{N}{2} (1-\theta) q^2 + O(q^2 + Nq^4)
    \rt).
  \eeq
  Since $\theta < 1$, we may take $q_0$ small enough that the error term is at most $\fr{N}{4} (1-\theta) q^2$ for all $|q| \le q_0$.
  This yields the upper bound
  \[
    N^C \exp\lt(-\fr{N}{4} (1-\theta) \lt(q_N - \fr2N\rt)^2 \rt) \le N^C \exp\lt(-\fr{1-\theta}{5} \log^2 N\rt) = o_N(1). \qedhere
  \]
\end{proof}

\begin{cor}%[Spherical central masses, high temperature]
\label{c:high-temperature-spherical-central-masses}
  Suppose $\theta<1$.
  For $\star\in\{\per, \mix, \hop\}$,
  \[
    \int_{|q|\le q_N}I_{\star,N}(q)\,\rho_N^\sph(\de q)
    =1+o_N(1).
  \]
\end{cor}
\begin{proof}
  Let $q_0$ be given by Lemma~\ref{l:high-temperature-localization}.
  Using the identity \eqref{e:reweight-hop} from Lemma~\ref{l:reweight}, we write
  \[
    \int_{|q|\le q_N}I_{\hop,N}(q)\,\rho_N^\sph(\de q)
    = 1 - \int_{q_N < |q| \le q_0}I_{\hop,N}(q)\,\rho_N^\sph(\de q) - \bbE \la \bbone \{|R_{1,2}| > q_0\} \ra_\bG.
  \]
  The last two terms are $o_N(1)$ by Lemma~\ref{l:high-temperature-localization} and Lemma~\ref{l:hopfield-overlap-concentration}, which proves the result for $\star = \hop$.
  The cases $\star \in \{\per, \mix\}$ follow from Proposition~\ref{p:central-comparison}.
\end{proof}

\begin{proof}[Proof of Theorem~\ref{t:perceptron-comparison}\ref{i:spherical-comparison} for $\theta < 1$]
  By Lemmas~\ref{l:perceptron-large-overlap} and~\ref{l:high-temperature-localization}, and Corollary~\ref{c:high-temperature-spherical-central-masses},
  \balnn
    \int_{-1}^1 I_{\per,N}(q)\,\rho_N^\sph(\de q) &= 1 + o_N(1), &
    \int_{-q_N}^{q_N} I_{\mix,N}(q)\,\rho_N^\sph(\de q) &= 1 + o_N(1).
  \ealnn
  Lemma~\ref{l:reweight} implies $\bbE[(X^\sph_N)^2] = 1 + o_N(1)$ and $\bbE[X^\sph_N] \ge 1 - o_N(1)$.
  Thus $\bbE[(X_N^\sph-1)^2] \le o_N(1)$.
\end{proof}

\hypertarget{p:t:perceptron-ising-high-temp}{}
\begin{proof}[Proof of Theorem~\ref{t:perceptron-comparison}\ref{i:ising-comparison} for $\theta < 1$]
  Let $\hat q_N = q_N - 1/N$.
  It suffices to prove the following estimate, after which we can proceed as in the proof of part \ref{i:spherical-comparison}:
  \beq
    \label{e:high-temperature-ising-central}
    \int_{|q|\le \hat q_N} I_{\per,N}(q)\,\rho_N^\Is(\de q)
    =1+o_N(1).
  \eeq
  The intervals $[q - 1/N, q + 1/N]$ with $q\in\cQ_N \cap [-\hat q_N, \hat q_N]$ tile an interval $[-\tilde q_N,\tilde q_N]$ with $|\tilde q_N-\hat q_N|\le 1/N$.
  In particular $\tilde q_N \le q_N$.
  Let $\tau : [-\tilde q_N, \tilde q_N] \to [-\hat q_N, \hat q_N]$ be the map that rounds to the nearest element of $\cQ_N \cap [-\hat q_N, \hat q_N]$, breaking ties arbitrarily.
  Lemma~\ref{l:ising-grid} gives
  \balnn
    \int_{|q|\le \hat q_N} I_{\per,N}(q)\,\rho_N^\Is(\de q)
    &= (1+o_N(1))
    \sum_{q\in\cQ_N \cap [-\hat q_N, \hat q_N]}
    I_{\per,N}(q)
    \,\rho_N^\sph([q-1/N,q+1/N]) \\
    &= (1+o_N(1))
    \int_{-\tilde q_N}^{\tilde q_N}
    I_{\per,N}(\tau(q))
    \,\rho_N^\sph(\de q).
  \ealnn
  Proposition~\ref{p:central-comparison}, with this choice of $\tau$, gives
  \beq\label{e:high-temperature-ising-central-1}
    \int_{-\tilde q_N}^{\tilde q_N}
    I_{\per,N}(\tau(q))
    \,\rho_N^\sph(\de q)
    = \int_{-\tilde q_N}^{\tilde q_N}
    I_{\hop,N}(q)
    \,\rho_N^\sph(\de q) + o_N(1).
  \eeq
  Since $q_N-2/N\le\tilde q_N\le q_N$, Lemma~\ref{l:high-temperature-localization} bounds
  \[
    \int_{\tilde q_N < |q| \le q_N}
    I_{\hop,N}(q)
    \,\rho_N^\sph(\de q) = o_N(1).
  \]
  Combining the last display with Corollary~\ref{c:high-temperature-spherical-central-masses} shows the right-hand side integral of \eqref{e:high-temperature-ising-central-1} is $1+o_N(1)$.
  Then \eqref{e:high-temperature-ising-central} follows.
\end{proof}

\subsection{Completion of the proof at criticality}
\label{ss:proof-assembly-critical}

In this subsection, we complete the proof of Theorem~\ref{t:perceptron-comparison} in the case $\theta = 1$.
We set $q_N = N^{-1/3} \log N$ as in \eqref{e:qN}.

The main difference in the critical setting is that $\cF''_\sph(0) = \cF''_\Is(0) = \theta - 1 = 0$, so the quadratic curvature of $\cF_\ast$ at $0$ is not enough to force the integrals from Lemma~\ref{l:reweight} to localize.
In particular, at criticality the leading quadratic term in \eqref{e:high-temperature-localization} vanishes.
The localization will instead be due to the contribution of the planted $Z^\sph_{\hop,N}(\bG)^{-2}$ term, captured by Proposition~\ref{p:bbp}.

The following lemma is the critical counterpart to Lemma~\ref{l:high-temperature-localization}; the remaining proof is the same as above.
\begin{lem}%[Critical perceptron and Hopfield localization]
  \label{l:critical-localization}
  Suppose $\theta = 1$.
  There exists a sufficiently small constant $q_0 > 0$, depending on $\phi, \alpha$, such that for any $(\star,\ast) \in \{(\per,\Is), (\per, \sph), (\hop, \sph)\}$,
  \[
    \int_{q_N - \fr2N < |q| \le q_0}
    I_{\star,N}(q) \,\rho_N^\ast(\de q)
    = o_N(1).
  \]
\end{lem}
\begin{proof}
  We proceed as in the proof of Lemma~\ref{l:high-temperature-localization}, but using Proposition~\ref{p:bbp} in place of Proposition~\ref{p:uniform-planted-inverse-moments}.
  Thus
  \balnn
    \int_{q_N - \fr2N < |q| \le q_0} I_{\star,N}(q) \rho_N^\ast(\de q)
    &= \int_{q_N - \fr2N < |q| \le q_0}
    \lt[s\Phi_{\star}(q)\rt]^M
    \bbE_{\star,q} \lt[
      Z^\sph_{\hop,N}(\bG)^{-2}
    \rt] 
    \,\rho_N^\ast(\de q) \\
    &\le N^C \int_{q_N - \fr2N < |q| \le q_0}
    e^{-cN|q|^3}
    \Phi_{\star}(q)^M
    \,\rho_N^\ast(\de q).
  \ealnn
  The density of $\rho_N^\sph$ and mass of $\rho_N^\Is$ are both bounded by $\exp(-Nq^2/2 + O(q^2 + Nq^4))$, up to polynomial factors in $N$.
  Recall $\Phi_\star(q)^M$ has Taylor expansion \eqref{e:Phi-taylor} around $q=0$, which has leading-order quadratic term $\fr{N}{2} q^2$ because $\theta = \alpha (1-s)^2 = 1$.
  After adjusting $C$, the last display is bounded by
  \beq\label{e:critical-localization}
    N^C \sup_{q_N - \fr2N < |q| \le q_0}
    \exp\lt(
      -cN|q|^3 + O(q^2 + Nq^4)
    \rt).
  \eeq
  We may take $q_0$ small enough that the error term is at most $cN|q|^3 / 2$ for all $q_N - \fr2N < |q| \le q_0$, and the result follows.
\end{proof}

\begin{cor}%[Spherical central masses, critical temperature]
\label{c:critical-spherical-central-masses}
  Suppose $\theta = 1$.
  For $\star\in\{\per, \mix, \hop\}$,
  \[
    \int_{|q|\le q_N}I_{\star,N}(q)\,\rho_N^\sph(\de q)
    =1+o_N(1).
  \]
\end{cor}
\begin{proof}
  Identical to Corollary~\ref{c:high-temperature-spherical-central-masses}, but using Lemma~\ref{l:critical-localization} in place of Lemma~\ref{l:high-temperature-localization}.
\end{proof}

\begin{proof}[Proof of Theorem~\ref{t:perceptron-comparison} for $\theta = 1$]
  Both parts are analogous to the proofs for $\theta < 1$; we omit the details.
\end{proof}

\subsection{Deferred proofs}
\label{ss:proof-deferred}

In this subsection, we prove the inputs used above that have elementary proofs.

\begin{proof}[Proof of Lemma~\ref{l:reweight}]
  Let $\nu_N^\sph = \nu_N = \unif(S_N)$ and $\nu_N^\Is = \unif(\Sigma_N)$.
  For $\ast \in \{\sph,\Is\}$,
  \[
    \bbE[(X^\ast_N)^2]
    = \fr{\bbE[Z^\sph_{\hop,N}]^2}{\bbE[Z^\ast_{\per,N}]^2}
    \iint \bbE \lt\{
      \prod_{a=1}^M \lt[
        \phi\lt(\fr{(\bg^a,\bx^1)}{\sqrt{N}}\rt)
        \phi\lt(\fr{(\bg^a,\bx^2)}{\sqrt{N}}\rt)
      \rt]
      Z^\sph_{\hop,N}(\bG)^{-2}
    \rt\}
    \,\nu^\ast_N(\de \bx^1)
    \,\nu^\ast_N(\de \bx^2).
  \]
  By rotational invariance, the inner expectation depends on $\bx^1,\bx^2$ only through $q = R(\bx^1,\bx^2)$.
  By a standard change of measure argument, it evaluates to (for $\bbE_q$ defined in Definition~\ref{d:planted-models})
  \[
    \bbE_q[\phi(Z_1)\phi(Z_2)]^M \cdot \bbE_{\per,q} \lt[
      Z^\sph_{\hop,N}(\bG)^{-2}
    \rt]
    = \bbE[Z^\ast_{\per,N}]^2 
    \cdot \Phi_\per(q)^M \bbE_{\per,q} \lt[
      Z^\sph_{\hop,N}(\bG)^{-2}
    \rt].
  \]
  Further recall from \eqref{e:annealed-free-energies} that $\bbE[Z^\sph_{\hop,N}] = s^{M/2}$.
  Thus
  \[
    \bbE[(X^\ast_N)^2]
    = \int \lt[s\Phi_{\per}(q)\rt]^M
    \bbE_{\per,q} \lt[
      Z^\sph_{\hop,N}(\bG)^{-2}
    \rt]\,\de \rho^\ast_N(q)
    = \int I_{\per,N}(q)\,\de \rho^\ast_N(q),
  \]
  which is the first identity of \eqref{e:reweight-per-mix}.
  For the second identity of \eqref{e:reweight-per-mix}, write
  \balnn
    &\bbE[X^\ast_N]
    = \fr{\bbE[Z^\sph_{\hop,N}]}{\bbE[Z^\ast_{\per,N}]}
    \bbE\lt[\fr{Z^\ast_{\per,N} Z^\sph_{\hop,N}}{(Z^\sph_{\hop,N})^2}\rt] \\
    &= \fr{\bbE[Z^\sph_{\hop,N}]}{\bbE[Z^\ast_{\per,N}]}
    \iint \bbE \lt\{
      \prod_{a=1}^M \lt[
        \phi\lt(\fr{(\bg^a,\bx^1)}{\sqrt{N}}\rt)
        \exp\lt(\fr{\beta (\bg^a,\bx^2)^2}{2N}\rt)
      \rt]
      Z^\sph_{\hop,N}(\bG)^{-2}
    \rt\}
    \,\nu^\ast_N(\de \bx^1)
    \,\nu^\sph_N(\de \bx^2),
  \ealnn
  and change measure as above; since $\bx^2$ is integrated over $\nu_N^\sph$, $q = R(\bx^1,\bx^2)$ has law $\rho^\sph_N$.
  The identity \eqref{e:reweight-hop} is proved analogously, by introducing two copies of $Z^\sph_{\hop,N}$ in the numerator and denominator.
\end{proof}

\begin{proof}[Proof of Lemma~\ref{l:good-gram}]
  Recall the definition \eqref{e:bSigma-hop} of $\bSigma_{\hop,q}$.
  For $\star \in \{\per, \mix\}$, similarly define $\bSigma_{\star,q} = \bbE_{\star,q}[\bK]$.
  At $q=0$, the choice \eqref{e:beta} of $\beta$ in the Hopfield model ensures that for all $\star \in \{\per, \mix, \hop\}$
  \balnn
    \bSigma_{\star,0}
    &= Ms \cdot I_2, &
    \fr{\de}{\de q} \bSigma_{\star,q} \big|_{q=0}
    &= Ms(s-1) \begin{pmatrix} 1 & 0 \\ 0 & -1 \end{pmatrix},
  \ealnn
  i.e. the Taylor expansions of the $\bSigma_{\star,q}$ at $q=0$ differ only at second order.
  It follows that for some $C > 0$ depending only on $\phi,\alpha$, and all $|q| \le \oq_N$,
  \[
    \max_{\star \in \{\per, \mix\}}
    \|\bSigma_{\star,q} - \bSigma_{\hop,q}\|_F \le CMq^2.
  \]
  Thus, for any $\star \in \{\per, \mix, \hop\}$,
  \baln\nonumber
    \bbP_{\star,q}(\bK \not\in \cK_N(q))
    &\le \bbP_{\star,q}\lt(
      \|\bK - \bSigma_{\star,q}\|_F > N^{2/3 - 2\eps_0} - CMq^2
    \rt) \\
    \label{e:good-gram}
    &\le \bbP_{\star,q}\lt(
      \|\bK - \bSigma_{\star,q}\|_F > N^{2/3 - 2\eps_0} / 2
    \rt),
  \ealn
  where the last step uses that $|q| \le \oq_N$.
  We will bound this probability using Bernstein's inequality.
  Let $\bc^1,\ldots,\bc^M \in \bbR^2$ be the rows of $\bC \in \bbR^{M\times 2}$ and recall that
  \[
    \bK = \bC^\top \bC = \sum_{a=1}^M \bc^a (\bc^a)^\top.
  \]
  Under any of the $\bbP_{\star,q}$, the $\bc^a$ are i.i.d., and are uniformly subexponential in the sense that for some $\eta > 0$,
  \[
    \sup_{|q|\le \oq_N}
    \max_{\star\in\{\per,\mix,\hop\}}
    \bbE_{\star,q}
    \exp(\eta \|\bc^1\|^2)
    <\infty.
  \]
  Scalar Bernstein's inequality and a union bound over the four entries yield, for constants $C', c$ depending only on $\phi, \alpha$,
  \[
    \bbP_{\star,q}\lt(
      \|\bK-\bSigma_{\star,q}\|_F > t
    \rt)
    \le C'\exp\lt\{
      -c\min\lt(\frac{t^2}{N},t\rt)
    \rt\},
  \]
  which gives the bound $C\exp(-cN^{1/3 - 4\eps_0})$ when $t = N^{2/3 - 2\eps_0}/2$ as in \eqref{e:good-gram}.
\end{proof}

\begin{proof}[Proof of Lemma~\ref{l:ising-grid}]
  The spherical overlap law has density
  \[
    \rho_N^\sph(\de q)
    =\fr{\Gamma(N/2)}{\sqrt{\pi}\,\Gamma((N-1)/2)}
    (1-q^2)^{(N-3)/2}\,\de q.
  \]
  Taylor expanding this density at $q$ and integrating over $[q-1/N,q+1/N]$ yields
  \[
    \rho_N^\sph\lt([q-1/N,q+1/N]\rt) 
    = \lt(1 + O(N^{-1}+q^2)\rt) \fr2N \cdot 
    \fr{\Gamma(N/2)}{\sqrt{\pi}\,\Gamma((N-1)/2)}
    (1-q^2)^{(N-3)/2}.
  \]
  The error scale is $O(N^{-1}+q^2)$ because the linear Taylor term vanishes after integration by symmetry.
  Stirling's approximation gives, for $q \in \cQ_N$ in a neighborhood of $0$,
  \balnn
    \rho_N^\sph\lt([q-1/N,q+1/N]\rt) 
    &= \lt(1+O(N^{-1}+q^2)\rt)
    \sqrt{\fr{2}{\pi N}}
    \exp\lt(-\fr{N}{2} q^2 + O(q^2 + Nq^4) \rt), \\
    \rho_N^\Is(\{q\}) 
    = 2^{-N} \binom{N}{N \cdot \fr{1+q}{2}}
    &= (1 + O(N^{-1}))
    \sqrt{\fr{2}{\pi N}}
    \exp\lt(-\fr{N}{2} q^2 + O(q^2 + Nq^4) \rt).
  \ealnn
  Since $N\oq_N^4=o(1)$, comparison of the last two lines proves the result.
\end{proof}

%\newpage

\section{Inverse partition function moment bounds}
\label{s:rmt}

In this section, we prove Propositions~\ref{p:uniform-planted-inverse-moments}--\ref{p:bbp}.
In \S\ref{ss:rmt-outline}, we give the proofs of both propositions, assuming several lemmas whose proofs are deferred to \S\ref{ss:coarea-upper-bound}--\ref{ss:inward-deviations}.
The remainder of this section proves the deferred lemmas.
As in \S\ref{s:proof}, we assume $\theta \le 1$, but not $\theta < 1$ or $\theta = 1$ unless stated.

\subsection{Outline of the proof}
\label{ss:rmt-outline}

In this subsection we prove Propositions~\ref{p:uniform-planted-inverse-moments}--\ref{p:bbp} assuming Lemmas~\ref{l:coarea-upper-bound}--\ref{l:wishart-inward-deviation-bbp} below. If $\beta = 0$, then $s=1$, $\theta = 0$, and $Z^\sph_{\hop,N} \equiv 1$, causing Propositions~\ref{p:uniform-planted-inverse-moments}--\ref{p:bbp} to both become trivial.
We thus assume $\beta \neq 0$.
We begin with a pointwise upper bound for the integrand in Propositions~\ref{p:uniform-planted-inverse-moments} and \ref{p:bbp}.
This is proved through a coarea representation of $Z^\sph_{\hop,N}(\bG)$ (Lemma~\ref{l:hopfield-coarea}), which is closely related to a contour integral representation used in several earlier works \cite{baik2016fluctuations,baik2018ferromagnetic,baik2020free,baik2021spherical,landon2022free,johnstone2024spin,collins2025fluctuations}.

\begin{lem}[Proved in \S\ref{ss:coarea-upper-bound}]\label{l:coarea-upper-bound}
  Let $\bW = \fr1N \bG^\top \bG \in \bbR^{N\times N}$ have eigenvalues $\lambda_1 \ge \cdots \ge \lambda_N \ge 0$.
  Let $\mu_i = \beta \lambda_i$ and $\mu_+ = \max_i(\mu_i)$.
  Let $\zeta > \mu_+$ be the unique solution to
  \beq\label{e:zeta}
    \sum_{i=1}^N \fr{1}{\zeta-\mu_i} = N.
  \eeq
  Then, there exists a universal $C>0$ such that
  \[
    Z^\sph_{\hop,N}(\bG)^{-2} \le CN e^{N(1-\zeta)} \det (\zeta - \beta \bW).
  \]
\end{lem}

This estimate is tight up to a polynomial factor.
The main difficulty in controlling the resulting upper bound is that $\zeta$ is itself random and correlated with the disorder $\bG$.
In the proofs of Propositions~\ref{p:uniform-planted-inverse-moments} and \ref{p:bbp} below, we will use a discretization scheme to reduce to expectations where $\zeta$ is replaced by a deterministic number.
The next two lemmas control such expectations, by first comparing the planted expectation to an expectation over an i.i.d. Gaussian $\bG$, and then estimating that expectation up to a polynomial factor.

\begin{lem}[Proved in \S\ref{ss:gaussian-replacement}]\label{l:gaussian-replacement-matrix}
  Let $\bbE_\nul$ denote expectation with respect to $\bG \in \bbR^{M\times N}$ with i.i.d. standard Gaussian entries.
  Let $\bW = \fr1N \bG^\top \bG \in \bbR^{N\times N}$.
  There exists $C > 0$ depending only on $\phi,\alpha$ such that for all fixed $z\in \bbR$, $\star \in \{\per, \mix\}$, $q\in [-1,1]$, and $k\in \{1,2\}$,
  \[
    \bbE_{\star,q}[|\det (z - \bW)|^{2k}] \le CN^{4k} \bbE_\nul[|\det (z - \bW)|^{2k}].
  \]
  The same estimate holds for $\star = \hop$ and all $q\in J$, for any compact $J\subseteq \cI_\adm(s)$, where $C$ may also depend on $J$ in addition to $\phi,\alpha$.
\end{lem}
\begin{lem}[Proved in \S\ref{ss:characteristic-polynomial-moment}]\label{l:characteristic-polynomial-moment}
  Let $\nu_{\MP,\alpha}$ be the Marchenko--Pastur probability measure with continuous support $[(\alpha^{1/2} - 1)^2, (\alpha^{1/2} + 1)^2]$. Write $\ell_{\pm}=(\sqrt\alpha \pm 1)^2$; the measure includes its atom of mass $(1-\alpha)_+$ at zero.
  Let
  \beq\label{e:L-alpha}
    L_\alpha(z) = \int \log|z-x|\,\nu_{\MP,\alpha}(\de x).
  \eeq
  For any compact $J \subseteq \bbR \setminus \{0\}$, there exists $C>0$ depending only on $\phi, \alpha, J$ such that for all $z\in J$,
  \beq\label{e:characteristic-polynomial-moment}
    \bbE_\nul[|\det (z - \bW)|^4]
    \le CN^{14} \exp\lt( 4N L_\alpha(z) \rt).
  \eeq
  For any $\eta > 0$, if $\alpha \ge 1+\eta$, \eqref{e:characteristic-polynomial-moment} also holds for compact $J \subseteq \bbR$ (rather than $\bbR \setminus \{0\}$), where $C$ depends only on $\phi, \eta, J$.
\end{lem}
\begin{rmk}
  Since the limiting spectral measure of $\bW$ is $\nu_{\MP,\alpha}$, one heuristically expects $\log |\det (z-\bW)|$ to concentrate around $N \int \log|z-x|\,\nu_{\MP,\alpha}(\de x)$.
  Thus, the exponential term in \eqref{e:characteristic-polynomial-moment} gives the expected exponential order of $\bbE_\nul[|\det (z - \bW)|^4]$.
  This heuristic holds quite generally, as was shown in \cite{ben2023exponential} for determinantal moments of a broad class of random matrices, up to an error of $e^{o(N)}$.
  As we will see below, it is important in our setting that the error term is polynomial.
\end{rmk}
The remaining inputs we need are lower deviation bounds on $\zeta$ and the extremal eigenvalue $\mu_+$.
We first record their limiting values in the non-planted model:
\baln\label{e:zeta-and-mu-ast}
  \zeta_* &= 1 + \alpha(s-1), &
  \mu_* &= \max\{\beta x : x \in \supp(\nu_{\MP,\alpha})\}
  = \begin{cases}
    \beta \ell_+ & \beta \ge 0, \\
    \beta \ell_- & \beta < 0, \alpha > 1, \\
    0 & \beta < 0, \alpha \le 1.
  \end{cases}
\ealn
\begin{fac}\label{f:zeta-above-spectrum}
  We have $\zeta_* \ge \mu_*$.
  Equality holds if and only if $\theta = 1$, in which case $\zeta_* = \mu_* \neq 0$.
\end{fac}
\begin{proof}
  This is a direct algebraic check in all three cases of \eqref{e:zeta-and-mu-ast}, using $\beta = \fr{s-1}{s}$ and $\theta = \alpha (s-1)^2 \le 1$.
\end{proof}
We now record the lower deviation bounds on $\zeta$ and $\mu_+$.
Lemma~\ref{l:zeta-deviation} is proved by a standard concentration of spectral traces argument.
Lemma~\ref{l:wishart-inward-deviation} is a general result on inward deviations of the extremal eigenvalues of a non-planted Wishart matrix, which implies a lower deviation bound on $\mu_+$ by interlacing; Lemma~\ref{l:wishart-inward-deviation-bbp} is a stronger version at criticality that witnesses the $cq^2$ change in $\mu_+$ from the slightly supercritical BBP spike.
These are proved by Golub--Kahan bidiagonalization.
\begin{lem}[Proved in \S\ref{ss:inward-deviations}]\label{l:zeta-deviation}
  For every $u > 0$, there exists $c > 0$ depending only on $\phi, \alpha, u$, such that for all $\star \in \{\per, \mix\}$ and $q\in [-1,1]$, or $\star = \hop$ and $q \in \cI_\adm(s)$,
  \[
    \bbP_{\star,q}(\zeta < \zeta_* - u) \le e^{-cN^2}.
  \]
\end{lem}
\begin{lem}[Proved in \S\ref{ss:inward-deviations}]\label{l:wishart-inward-deviation}
  Let $\bA_0 \in \bbR^{M \times (N-2)}$ have i.i.d. standard Gaussian entries and $\bW_0 = \fr1N \bA_0^\top \bA_0$.
  Let $\ell_\pm = (\alpha^{1/2} \pm 1)^2$.
  There exists $c > 0$, depending only on $\alpha$, such that for all $0 \le u \le 1$,
  \balnn
    \bbP\lt(\lambda_{\max}(\bW_0) \le \ell_+ - u\rt) &\le 2e^{-cN^2u^3}, &
    \bbP\lt(\lambda_{\min}(\bW_0) \ge \ell_- + u\rt) &\le 2e^{-cN^2u^3}.
  \ealnn
\end{lem}
\begin{lem}[Proved in \S\ref{ss:inward-deviations}]\label{l:wishart-inward-deviation-bbp}
  Suppose $\theta = 1$, so that $q_N = N^{-1/3} \log N$.
  There exist $c,q_0 > 0$ depending only on $\phi,\alpha$ such that for all $\star \in \{\per, \mix, \hop\}$ and all $|q| \in [q_N/2, q_0]$, $u\in [0,1]$,
  \[
    \bbP_{\star,q}\lt(
      \mu_+ \le \mu_* + cq^2 - u
    \rt)
    \le 2\exp\lt(-cN \fr{(u+q^2)^2}{|q|}\rt).
  \]
\end{lem}

The remainder of this subsection presents the proof of Propositions~\ref{p:uniform-planted-inverse-moments}--\ref{p:bbp}.
The proofs will share a common strategy and differ only in the final step.
In the rest of this subsection, $C,c > 0$ denote constants uniform in $\star \in \{\per, \mix\}$, $q\in [-1,1]$ or $\star = \hop$, $q\in J$ where $J \subseteq \cI_\adm(s)$ is compact, which may depend arbitrarily on $\phi, \alpha, J$ and any other parameters stated.

The following lemma shows that if $\nu_{\MP,\alpha}$ has an atom at $0$, causing Lemma~\ref{l:characteristic-polynomial-moment} to only apply uniformly away from $0$, then $\zeta_*$ is away from $0$.
\begin{lem}\label{l:saddle-away-from-0}
  There exists $\eta_0 > 0$, depending only on $\phi$, such that for all $\alpha \le 1 + \eta_0$, we have $\zeta_* \ge 2\eta_0$.
\end{lem}
\begin{proof}
  Take $\eta_0 = \fr13\min(1,s)$.
  If $s\ge 1$, then $\zeta_* \ge 1 \ge 2\eta_0$.
  If $s\in (0,1)$ and $\alpha \le 1+\eta_0$, then
  \[
    \zeta_*
    = s - (\alpha-1)(1-s)
    \ge s - \eta_0(1-s)
    = \fr{s(2+s)}{3}
    \ge 2\eta_0. \qedhere
  \]
\end{proof}
For the rest of this subsection, let $\eta \in (0,\eta_0)$ be small and $r > 0$ be large, depending only on $\phi, \alpha, J$.
We will set $r$ in Lemma~\ref{l:zeta-discretization-upper-tail} and $\eta$ in the proofs of Propositions~\ref{p:uniform-planted-inverse-moments}--\ref{p:bbp} below.
\begin{lem}\label{l:zeta-discretization}
  There exists $C$, uniform over $\star, q$ in the sense described just below Lemma~\ref{l:wishart-inward-deviation-bbp} and depending only on $\phi, \alpha, J$, such that for all $k\in \{1,2\}$,
  \beq\label{e:zeta-discretization}
    \begin{aligned}
    \bbE_{\star,q}\lt[
      Z^\sph_{\hop,N}(\bG)^{-2k}
    \rt]
    &\le CN^{4+7k} |\beta|^{Nk} \sup \lt\{
      e^{Nk(1-z + L_\alpha(z / \beta))} 
      \bbP_{\star,q}(\zeta \le z)^{1/2}
      : \zeta_* - \eta \le z \le \zeta_* + r
    \rt\} \\
    &\quad + CN^k \bigg[ \bbE_{\star,q} \lt[
      \bbone\{\zeta < \zeta_* - \eta\}
      e^{Nk(1-\zeta)} \det (\zeta - \beta \bW)^k
    \rt] \\
    &\quad + \bbE_{\star,q} \lt[
      \bbone\{\zeta > \zeta_* + r\}
      e^{Nk(1-\zeta)} \det (\zeta - \beta \bW)^k
    \rt]
    \bigg]
    \end{aligned}
  \eeq
\end{lem}
\begin{proof}
  Throughout this proof $C$ will change from line to line but will be uniform in the above sense.
  Let $K = \lceil (\eta+r) N^4 \rceil$ and $\zeta_* - \eta = z_0 < z_1 < \cdots < z_K = \zeta_* + r$ be such that $|z_j - z_{j-1}| \le N^{-4}$ for all $1\le j\le K$.
  Then define $I_j = [z_{j-1}, z_j]$ and $I_0 = (-\infty, \zeta_* - \eta)$, $I_{K+1} = (\zeta_* + r, +\infty)$.
  By Lemma~\ref{l:coarea-upper-bound},
  \balnn
    \bbE_{\star,q}\lt[
      Z^\sph_{\hop,N}(\bG)^{-2k}
    \rt]
    &\le CN^k \bbE_{\star,q}\lt[
      e^{Nk(1-\zeta)} \det (\zeta - \beta \bW)^k
    \rt] \\
    &\le CN^k \sum_{j=0}^{K+1} \bbE_{\star,q}\lt[
      \bbone\{\zeta \in I_j\}
      e^{Nk(1-\zeta)} \det (\zeta - \beta \bW)^k
    \rt].
  \ealnn
  The terms $j=0$ and $j=K+1$ are precisely the last two terms in \eqref{e:zeta-discretization}.
  For $1\le j\le K$, we have by Cauchy--Schwarz and Lemma~\ref{l:gaussian-replacement-matrix}
  \balnn
    \bbE_{\star,q}\lt[
      \bbone\{\zeta \in I_j\}
      e^{Nk(1-\zeta)} \det (\zeta - \beta \bW)^k
    \rt]
    &\le e^{Nk(1-z_{j-1})}
    \bbE_{\star,q}\lt[
      \bbone\{\zeta \le z_j\}
      |\det (z_j - \beta \bW)|^k
    \rt] \\
    &\le
    e^{Nk(1-z_{j-1})}
    \bbP_{\star,q}(\zeta \le z_j)^{1/2}
    \bbE_{\star,q}\lt[
      |\det (z_j - \beta \bW)|^{2k}
    \rt]^{1/2} \\
    &\le 
    CN^{2k} e^{Nk(1-z_j)}
    \bbP_{\star,q}(\zeta \le z_j)^{1/2}
    \bbE_\nul \lt[
      |\det (z_j - \beta \bW)|^{2k}
    \rt]^{1/2}.
  \ealnn
  In the last step we have also absorbed the factor $e^{Nk(z_j-z_{j-1})} \le e^{kN^{-3}}$ into $C$.
  We next claim that
  \beq\label{e:characteristic-polynomial-moment-application}
    \bbE_\nul \lt[
      |\det (z_j - \beta \bW)|^{2k}
    \rt]
    \le CN^{7k} |\beta|^{2Nk} \exp(2Nk L_\alpha(z_j / \beta)),
  \eeq
  where $C$ is uniform over all $1\le j\le K$.

  Note that the $z_j$ lie in the compact set $I = [\zeta_* - \eta, \zeta_* + r]$.
  By Lemma~\ref{l:saddle-away-from-0}, if $\alpha \le 1 + \eta$, $I$ does not include $0$.
  Since both $\eta$ and $\beta = (s-1)/s$ depend only on $\phi, \alpha$, the $z_j / \beta$ lie in a compact set that depends only on $\phi,\alpha,r$, which does not include $0$ if $\alpha \le 1 + \eta$. If instead $\alpha>1+\eta$, the second assertion of Lemma~4.3 applies without requiring the compact set to avoid zero.
  We may thus apply Lemma~\ref{l:characteristic-polynomial-moment} to $\fr{z_j}{\beta} - \bW$ to obtain the case $k=2$ of \eqref{e:characteristic-polynomial-moment-application}.
  The case $k=1$ follows from $k=2$ by Cauchy--Schwarz.
  Altogether we have shown
  \[
    \bbE_{\star,q}\lt[
      \bbone\{\zeta \in I_j\}
      e^{Nk(1-\zeta)} \det (\zeta - \beta \bW)^k
    \rt]
    \le
    CN^{6k} e^{Nk(1-z_j + L_\alpha(z_j / \beta))} |\beta|^{Nk}
    \bbP_{\star,q}(\zeta \le z_j)^{1/2}.
  \]
  The sum of this over $1\le j\le K$ is bounded by the first term in \eqref{e:zeta-discretization} because $K = O(N^4)$ and all the $z_j$ satisfy $\zeta_* - \eta \le z_j \le \zeta_* + r$.
\end{proof}
The next two lemmas show the last two terms in Lemma~\ref{l:zeta-discretization} are negligible.
\begin{lem}\label{l:zeta-discretization-lower-tail}
  There exists $c > 0$, uniform over $\star, q$ in the sense described just below Lemma~\ref{l:wishart-inward-deviation-bbp}, such that for all $k\in \{1,2\}$,
  \beq\label{e:zeta-discretization-lower-tail}
    \bbE_{\star,q} \lt[
      \bbone\{\zeta < \zeta_* - \eta\}
      e^{Nk(1-\zeta)} \det (\zeta - \beta \bW)^k
    \rt]
    \le e^{-cN^2} s^{-kM}.
  \eeq
\end{lem}
\begin{proof}
  By the inequality $\log x \le x-1$, we have the pointwise bound
  \[
    N(1-\zeta) + \log\det (\zeta - \beta \bW)
    \le -\beta \Tr(\bW)
    = -\beta \|\bG\|_F^2 / N.
  \]
  By this estimate and Cauchy--Schwarz,
  \balnn
    \bbE_{\star,q} \lt[
      \bbone\{\zeta < \zeta_* - \eta\}
      e^{Nk(1-\zeta)} \det (\zeta - \beta \bW)^k
    \rt]
    &\le \bbE_{\star,q} \lt[
      \bbone\{\zeta < \zeta_* - \eta\}
      e^{-\beta k \|\bG\|_F^2 / N}
    \rt] \\
    &\le \bbP_{\star,q}(\zeta < \zeta_* - \eta)^{1/2}
    \bbE_{\star,q} \lt[
      e^{-2\beta k \|\bG\|_F^2 / N}
    \rt]^{1/2}.
  \ealnn
  Lemma~\ref{l:zeta-deviation} bounds $\bbP_{\star,q}(\zeta < \zeta_* - \eta)^{1/2}$, while a direct calculation shows that the last expectation is bounded by $e^{CN}$ for some $C$ uniform in the sense described just below Lemma~\ref{l:wishart-inward-deviation-bbp}.
\end{proof}
\begin{lem}\label{l:zeta-discretization-upper-tail}
  There exists $r > 0$, depending only on $\phi,\alpha$, such that uniformly over $\star,q$ in the sense described just below Lemma~\ref{l:wishart-inward-deviation-bbp}, for all $k\in\{1,2\}$,
  \beq\label{e:zeta-discretization-upper-tail}
    \bbE_{\star,q}\lt[
      \bbone\{\zeta > \zeta_* + r\}
      e^{Nk(1-\zeta)} \det(\zeta-\beta\bW)^k
    \rt]
    \le e^{-N}s^{-kM}.
  \eeq
\end{lem}
\begin{proof}
  First, suppose $\beta \ge 0$.
  Then we have $\det(\zeta - \beta\bW) \le \zeta^N$, so
  \[
    \bbE_{\star,q}\lt[
      \bbone\{\zeta > \zeta_* + r\}
      e^{Nk(1-\zeta)} \det(\zeta-\beta\bW)^k
    \rt]
    \le \sup\lt\{
      z e^{1-z} : 
      z > \zeta_* + r
    \rt\}^{kN}.
  \]
  The conclusion follows by taking $r$ sufficiently large.
  If $\beta < 0$, we claim that $\zeta \le 1$ deterministically, so that for large enough $r$ the left-hand side of \eqref{e:zeta-discretization-upper-tail} is $0$.
  Indeed, we have $\mu_i \le 0$ for all $i$.
  If $\zeta > 1$, then
  \[
    \sum_{i=1}^N \fr{1}{\zeta - \mu_i} \le \fr{N}{\zeta} < N,
  \]
  contradicting \eqref{e:zeta}.
\end{proof}

\begin{proof}[Proof of Proposition~\ref{p:uniform-planted-inverse-moments}]
  Fix $r$ as in Lemma~\ref{l:zeta-discretization-upper-tail}, and define
  \[
    \Phi(z)=1-z+\log|\beta|+L_\alpha(z/\beta).
  \]
  Lemmas~\ref{l:zeta-discretization-lower-tail}--\ref{l:zeta-discretization-upper-tail} control the last two terms in Lemma~\ref{l:zeta-discretization}.
  Thus it suffices to show that for $k\in \{1,2\}$,
  \beq\label{e:uniform-planted-inverse-moments-goal}
    \sup\lt\{
      e^{Nk\Phi(z)} \bbP_{\star,q}(\zeta\le z)^{1/2} : \zeta_* - \eta \le z \le \zeta_* + r
    \rt\} \le Cs^{-kM}.
  \eeq
  Since $\theta \le 1$, direct evaluation shows
  \balnn
    \Phi(\zeta_*) &= - \alpha \log s, & 
    \Phi'(\zeta_*) &= 0.
  \ealnn
  Recall $\mu_*$ defined in \eqref{e:zeta-and-mu-ast}.
  For all $z > \mu_*$,
  \[
    \Phi''(z)=-\int\fr{1}{(z-\beta x)^2}\,\nu_{\MP,\alpha}(\de x) < 0.
  \]
  Thus $\Phi(z)$ is maximized on $[\mu_*,\infty)$ at $z = \zeta_*$.
  Using the trivial bound $\bbP_{\star,q}(\zeta\le z) \le 1$ we obtain
  \beq\label{e:uniform-planted-inverse-moments-outside-spectrum}
    \sup\lt\{
      e^{Nk\Phi(z)} \bbP_{\star,q}(\zeta\le z)^{1/2} : \mu_* \le z \le \zeta_* + r
    \rt\} \le s^{-\alpha kN}.
  \eeq
  Recall from Fact~\ref{f:zeta-above-spectrum} that $\zeta_* \ge \mu_*$.
  If $\zeta_* > \mu_*$, we may take $\eta < \zeta_* - \mu_*$ (which depends only on $\phi, \alpha$) and the proof is complete.

  Otherwise, by Fact~\ref{f:zeta-above-spectrum}, $\zeta_* = \mu_* \neq 0$.
  For sufficiently small $\eta$ and all $z \in (\mu_* - \eta, \mu_*)$, $z/\beta$ is part of the continuous spectrum $(\ell_-, \ell_+)$ of $\nu_{\MP,\alpha}$.
  A direct calculation shows
  \[
    \Phi'(z) = - 1 + \fr{1}{2\beta} + \fr{1-\alpha}{2z}
  \]
  on this interval, and this formula extends continuously to $\Phi'(\mu_*) = 0$.
  Thus there exists $C$ such that $|\Phi'(z)| \le C|z-\mu_*|$ for $z \in [\mu_* - \eta, \mu_*]$, which implies
  \beq\label{e:Phi-below-edge}
    \Phi(\mu_* - u)
    \le \Phi(\mu_*) + Cu^2
    = -\alpha \log s + Cu^2
  \eeq
  for all $0\le u\le \eta$.

  Recall $\zeta > \mu_+$ and $\bW_0 = \fr1N \bA_0^\top\bA_0$.
  By interlacing, $\mu_+ \ge \lambda_{\max}(\beta \bW_0)$.
  We combine these facts with Lemma~\ref{l:wishart-inward-deviation}, using the upper edge bound when $\beta > 0$ and the lower edge bound when $\beta < 0$, to obtain
  \[
    \bbP_{\star,q}(\zeta \le \mu_*-u)
    \le\bbP\lt(\lambda_{\max}(\beta \bW_0) \le \mu_*-u\rt)
    \le 2e^{-cN^2u^3}.
  \]
  Consequently, for $z=\mu_*-u$ with $0\le u\le\eta$,
  \[
    e^{Nk\Phi(z)}\bbP_{\star,q}(\zeta\le z)^{1/2}
    \le Cs^{-\alpha Nk}e^{CNu^2-cN^2u^3}
    \le Cs^{-kM},
  \]
  since $\sup_{u\ge0}(CNu^2-cN^2u^3)=O(N^{-1})$.
  Together with \eqref{e:uniform-planted-inverse-moments-outside-spectrum} this proves \eqref{e:uniform-planted-inverse-moments-goal}.
\end{proof}
\begin{proof}[Proof of Proposition~\ref{p:bbp}]
  This is similar to the preceding proof, but using Lemma~\ref{l:wishart-inward-deviation-bbp} in place of Lemma~\ref{l:wishart-inward-deviation}.
  Fix $r$ as in Lemma~\ref{l:zeta-discretization-upper-tail}, and define $\Phi$ as in the preceding proof.
  Since $\theta=1$, Fact~\ref{f:zeta-above-spectrum} yields $\zeta_* = \mu_* \neq 0$.
  We will take $\eta, q_0$ sufficiently small depending on $\phi, \alpha$, with $q_0$ small enough that $[-q_0,q_0] \subseteq \cI_\adm(s)$ and Lemma~\ref{l:wishart-inward-deviation-bbp} applies for all $|q| \in [q_N / 2, q_0]$.

  Lemmas~\ref{l:zeta-discretization-lower-tail}--\ref{l:zeta-discretization-upper-tail} control the last two terms in Lemma~\ref{l:zeta-discretization}, so it suffices to show
  \beq\label{e:bbp-inverse-moments-goal}
    \sup\lt\{
      e^{N\Phi(z)}\bbP_{\star,q}(\zeta\le z)^{1/2} : \mu_*-\eta\le z\le\mu_*+r
    \rt\}\le Ce^{-cN|q|^3}s^{-M}.
  \eeq
  As in the preceding proof, $\Phi(\mu_*)=-\alpha\log s$, $\Phi'(\mu_*)=0$, and $\Phi'$ is decreasing on $[\mu_*,\infty)$.
  A direct calculation shows that $\Phi'(\mu_* + u) \le -c\sqrt{u}$ as $u\downarrow 0$.
  After adjusting $c$, we have
  \[
    \Phi'(\mu_* + u) \le -c\sqrt{u}
  \]
  for all $0 \le u \le r$.
  Integrating and combining with \eqref{e:Phi-below-edge} shows
  \beq\label{e:Phi-below-and-above-edge}
    \Phi(\mu_* + u)
    \le
    \begin{cases}
      -\alpha \log s - cu^{3/2} & 0\le u \le r, \\
      -\alpha \log s + Cu^2 & -\eta \le u \le 0.
    \end{cases}
  \eeq
  We turn to estimating the term $\bbP_{\star,q}(\zeta \le z)$.
  Let $c'$ be the constant $c$ in Lemma~\ref{l:wishart-inward-deviation-bbp}.
  Since $\zeta > \mu_+$, Lemma~\ref{l:wishart-inward-deviation-bbp} gives 
  \beq\label{e:wishart-inward-deviation-bbp}
    \bbP_{\star,q}(\zeta \le \mu_* + c'q^2 - u)
    \le \bbP_{\star,q}(\mu_+ \le \mu_* + c'q^2 - u)
    \le 2\exp\lt(-c'N\fr{(u+q^2)^2}{|q|}\rt),
  \eeq
  for all $|q| \in [q_N/2, q_0]$ and $0\le u\le 1$.
  Take $\eta, q_0$ small enough that $\eta + c' q_0^2 \le 1$ and $c' q_0^2 \le r$.
  We now prove \eqref{e:bbp-inverse-moments-goal} by casework on $z$.

  \textit{Case 1: $z \ge \mu_* + c' q^2$.}
  By \eqref{e:Phi-below-and-above-edge} and the trivial bound $\bbP_{\star,q}(\zeta\le z) \le 1$,
  \[
    e^{N\Phi(z)}\bbP_{\star,q}(\zeta\le z)^{1/2}
    \le s^{-\alpha N}e^{-cN(z-\mu_*)^{3/2}}
    \le Ce^{-cN|q|^3}s^{-M},
  \]
  after adjusting $C,c$.

  \textit{Case 2: $\mu_* \le z \le \mu_* + c' q^2$.}
  By \eqref{e:Phi-below-and-above-edge}, $\Phi(z)\le-\alpha\log s$.
  Taking $u=0$ in \eqref{e:wishart-inward-deviation-bbp} yields
  \[
    e^{N\Phi(z)}\bbP_{\star,q}(\zeta\le z)^{1/2}
    \le s^{-\alpha N} \cdot \sqrt{2}e^{-c'N|q|^3/2}
    \le Ce^{-cN|q|^3}s^{-M}.
  \]

  \textit{Case 3: $\mu_* - \eta \le z \le \mu_*$.}
  Write $z=\mu_*-t$, where $0\le t\le\eta$.
  We bound $\Phi(z)$ using \eqref{e:Phi-below-and-above-edge} and $\bbP_{\star,q}(\zeta\le z)$ using \eqref{e:wishart-inward-deviation-bbp}, with $u = c'q^2 + t$.
  Note that
  \[
    u \le c'q_0^2+\eta \le 1,
  \]
  so \eqref{e:wishart-inward-deviation-bbp} applies.
  Thus (adjusting $C,c$ line by line)
  \balnn
    e^{N\Phi(z)}\bbP_{\star,q}(\zeta\le z)^{1/2}
    &\le s^{-\alpha N} \exp(CNt^2) \cdot \sqrt{2} \exp\lt(-\fr{c'N(u+q^2)^2}{2|q|}\rt) \\
    &\le Ce^{-cN|q|^3} s^{-M} \cdot \exp\lt(CNt^2 - \fr{c'Nu^2}{2|q|}\rt).
  \ealnn
  If we take $q_0$ small enough that $C \le \fr{c'}{2q_0}$, the last exponential is bounded by $1$ because $u\ge t \ge 0$.
\end{proof}

\subsection{Pointwise upper bound from coarea formula}
\label{ss:coarea-upper-bound}

In this subsection we prove Lemma~\ref{l:coarea-upper-bound}.
We adopt the notation of this lemma: fix any $\bG \in \bbR^{M\times N}$, and let $\bW = \bG^\top \bG / N$ have eigenvalues $\lambda_1\ge\cdots\ge\lambda_N\ge0$. 
Let $\mu_i = \beta \lambda_i$ and $\mu_+ = \max_i \mu_i$, and define $\zeta$ as in \eqref{e:zeta}.
We begin with a coarea representation of $Z_{\hop,N}^{\sph}(\bG)$ similar to \cite[Lemma 1.3]{baik2016fluctuations}, from which Lemma~\ref{l:coarea-upper-bound} follows by upper bounding $f_\zeta(N)^{-2}$.

\begin{lem}%[Coarea representation]
  \label{l:hopfield-coarea}
  Define $K_N=2^{N/2}\Gamma(N/2)N^{1-N/2}$.
  Let $Z_1,\ldots,Z_N$ be independent standard Gaussians.
  For any $z > \mu_+$, define
  \balnn
    a_i(z) &= (z - \mu_i)^{-1}, &
    S(z) &= \sum_{i=1}^N a_i(z) Z_i^2,
  \ealnn
  and let $f_z$ be the density of $S(z)$, as a random variable measurable in $Z_1,\ldots,Z_N$ (recall $\bG$, and thus the $\mu_i$, are fixed).
  Then,
  \beq\label{e:hopfield-coarea}
    Z_{\hop,N}^{\sph}(\bG)^{-2} 
    = K_N^{-2} e^{-Nz} \det\lt( zI_N-\beta\bW \rt) f_z(N)^{-2}.
  \eeq
  In particular, this holds at $z = \zeta$.
\end{lem}
\begin{proof}
  Let $\bY \sim \cN(\bzero, (zI_N-\beta\bW)^{-1})$.
  The eigenvalues of this covariance matrix are precisely $(a_i(z) : 1\le i\le N)$, so $\|\bY\|^2$ has the same distribution as $S(z)$.
  Write
  \[
    p(\by) = (2\pi)^{-N/2}\det(zI_N-\beta\bW)^{1/2}
    \exp\lt(-\frac12 \by^\top(zI_N-\beta\bW)\by\rt)
  \]
  for the probability density of $\bY$.
  Then, expanding in polar coordinates $\by= r\bx /\sqrt N$ where $\bx \in S_N$,
  \[
    \bbP(S(z) \le t)
    = \int_{\|\by\|^2 \le t} p(\by)\,\de \by 
    = |\bbS^{N-1}| \int_0^{\sqrt{t}} \int r^{N-1} p\lt(\fr{r \bx}{\sqrt N}\rt) \,\de \nu_N(\bx) \,\de r.
  \]
  Here $|\bbS^{N-1}|=2\pi^{N/2}/\Gamma(N/2)$ is the Hausdorff measure of $\bbS^{N-1}$.
  Differentiating in $t$ at $t=N$ yields
  \[
    f_z(N) = \fr{|\bbS^{N-1}| N^{(N-2)/2}}{2} 
    \int p(\bx) \,\de \nu_N(\bx)
    = \fr{|\bbS^{N-1}| N^{(N-2)/2}}{2(2\pi)^{N/2}} 
    e^{-zN/2} \det(zI_N-\beta\bW)^{1/2} 
    Z_{\hop,N}^{\sph}(\bG).
  \]
  Rearranging gives the result.
\end{proof}

The next lemma provides a deterministic bound for the term $f_\zeta(N)^{-2}$ from \eqref{e:hopfield-coarea}.

\begin{lem}
\label{l:weighted-chi-square-density}
 Let $a_1,\ldots,a_n>0$ and $Z_1,\ldots,Z_n$ be independent $\mathcal N(0,1)$ random variables. Also, let $\widetilde Z_1^2,\ldots,\widetilde Z_n^2$ be independent $\chi^2(3)$ variables independent of $\{ Z_i:1 \le i \le n \}$. Set
\begin{equation}\label{e:V-S-Sj}
    a_\Sigma=\sum_{i=1}^n a_i, \qquad V=\sum_{i=1}^n a_i^2, \qquad S=\sum_{i=1}^n a_iZ_i^2, \qquad S^{(j)}=a_j\widetilde Z_j^2+\sum_{i\ne j}a_iZ_i^2.
\end{equation}
Then there exists a universal constant $C>0$ such that, for $f_S, f_{S^{(j)}}$ the probability densities of $S$, $S^{(j)}$,
\beq\label{e:weighted-chi-square-density}
    f_S(a_\Sigma), f_{S^{(j)}}(a_\Sigma) \le CV^{-1/2}.
\eeq
\beq\label{e:weighted-chi-square-augmented-density}
    C^{-1}V^{-1/2} \le f_S(a_\Sigma).
\eeq
\end{lem}
\begin{proof}
Without loss of generality, we assume $V=\sum_{i=1}^n a_i^2=1$ and $a_1 \ge a_2 \ge \ldots \ge a_n>0$ in the following proof. Set $T=S-a_\Sigma$ and $T^{(j)}=S^{(j)}-a_\Sigma$. Denote $\mathrm{i}$ as the imaginary unit, their characteristic functions are
\[
    \varphi(t) =\bbE e^{\mathrm{i}tT} =e^{-\mathrm{i}ta_\Sigma}\prod_{i=1}^n(1-2\mathrm{i}ta_i)^{-1/2}, \qquad \varphi_j(t) =\bbE e^{\mathrm{i}tT^{(j)}} =\frac{\varphi(t)}{1-2\mathrm{i}ta_j}.
\]
For \eqref{e:weighted-chi-square-density}, our proof splits into the following cases: $a_1>1/2$ and $a_1\leq 1/2$.

\smallskip
\noindent\emph{Proof of \eqref{e:weighted-chi-square-density} when $a_1\le 1/2$.} For the characteristic functions we have $|\varphi_j(t)|\le |\varphi(t)| =\prod_{i=1}^n(1+4a_i^2t^2)^{-1/4}$.
Since $4a_i^2\in[0,1]$, concavity of $u\mapsto\log(1+ut^2)$ on interval $[0,1]$ gives $\log(1+4a_i^2t^2)\ge 4a_i^2\log(1+t^2)$. Summing over $i$ and using $\sum_i a_i^2=1$, we obtain
\begin{equation}\label{eq:bound-abs-varphi}
    |\varphi_j(t)|\le |\varphi(t)|\le\frac{1}{1+t^2}.
\end{equation}
Thus both characteristic functions are integrable. Fourier inversion therefore implies
\[
    f_{T^{(j)}}(0), f_T(0)\le\frac{1}{2\pi}\int_{\bbR}\frac{\de t}{1+t^2} =\frac12.
\]

\smallskip
\noindent\emph{Proof of \eqref{e:weighted-chi-square-density} when $a_1>1/2$.}  For $t\ge0$, we write $\varphi(t)$ in polar form as follows
\[
    \varphi(t)=A(t)e^{\mathrm{i}\psi(t)}, \qquad A(t)=\prod_{i=1}^n(1+4a_i^2t^2)^{-1/4}, \qquad \psi(t)=\frac12\sum_{i=1}^n\arctan(2a_it)-ta_\Sigma.
\]
For $\varphi_j(t)$, we have
\[
    \varphi_j(t)=A_j(t)e^{\mathrm{i}\psi_j(t)}, \qquad A_j(t)=\frac{A(t)}{\sqrt{1+4a_j^2t^2}}, \qquad \psi_j(t)=\psi(t)+\arctan(2a_jt).
\]
Now let $(\widetilde A,\widetilde\psi)$ denote either $(A,\psi)$ or $(A_j,\psi_j)$. For $\varepsilon>0$, we define the auxiliary function
\[
    B_\varepsilon(t) =\frac{e^{-\varepsilon t^2}\widetilde A(t)}{-\widetilde\psi'(t)}.
\]
We first show
\begin{equation}\label{eq:properties-B(t)}
    B_\varepsilon(t) \text{ is nonincreasing with limit zero and satisfies }B_\varepsilon(4)\le4.
\end{equation}
 Note that $\widetilde A(t)$ is decreasing on $[0,\infty)$ with $A(0)=1$ and tends to zero. Then by taking derivatives for functions $\psi(t),\psi_j(t)$, we have
\[
    \psi'(t) =-\sum_{i=1}^n\frac{4a_i^3t^2}{1+4a_i^2t^2}, \qquad \psi_j'(t) =\psi'(t)+\frac{2a_j}{1+4a_j^2t^2}.
\]
The derivatives $\psi'$ and $\psi_j'$ are strictly decreasing on $(0,\infty)$. Moreover, $a_1>1/2$ and $1+4a_j^2t^2\ge4a_jt$ imply that $-\psi'(t),-\psi_j'(t)\ge1/4$ for $t\ge4$. This shows \eqref{eq:properties-B(t)}. Therefore, integrating by parts gives
\begin{align*}
    \left|\int_4^\infty e^{-\varepsilon t^2}\widetilde A(t) e^{i\widetilde\psi(t)}\,\de t\right| \le B_\varepsilon(4) +\int_4^\infty|B_\varepsilon'(t)|\,\de t =2B_\varepsilon(4)\le8.
\end{align*}
The integral over $[0,4]$ has absolute value at most $4$. Thus we get
\[
    \left|\frac{1}{2\pi}\int_{\bbR} e^{-\varepsilon t^2}\varphi(t)\,\de t\right|\le\frac{12}{\pi}, \qquad \left|\frac{1}{2\pi}\int_{\bbR} e^{-\varepsilon t^2}\varphi_j(t)\,\de t\right|\le\frac{12}{\pi}.
\]
These integrals equal the densities at zero of $T+\sqrt{2\varepsilon}\,Z$ and $T^{(j)}+\sqrt{2\varepsilon}\,Z$, where $Z$ is an independent standard Gaussian. Since the densities of $S$ and $S^{(j)}$ are continuous on $(0,\infty)$, $f_T$ and $f_{T^{(j)}}$ are continuous at zero. Gaussian smoothing therefore converges to their values at zero as $\varepsilon\downarrow0$, yielding
\[
f_T(0)\le\frac{12}{\pi}, \qquad f_{T^{(j)}}(0)\le\frac{12}{\pi}.
\]
This completes the proof of \eqref{e:weighted-chi-square-density}.\\
For \eqref{e:weighted-chi-square-augmented-density}, we also consider two cases $a_1\leq \frac{\sqrt{\pi}}{192}$ and $a_1>\frac{\sqrt{\pi}}{192}$. For simplicity we write $\delta=\frac{\sqrt{\pi}}{192}$.

\smallskip
\noindent\emph{Proof of \eqref{e:weighted-chi-square-augmented-density} when $a_1\leq \delta$.} 
Put $H=(16a_1)^{-1}$. For $|t|\le H$, the logarithmic series expansion gives
\[
\varphi(t)=\exp\{-t^2+E(t)\}, \qquad E(t)=\sum_{i=1}^n\sum_{m=3}^{\infty}\frac{(2\mathrm ita_i)^m}{2m}.
\]
Since $2a_i|t|\le1/8$, the series is absolutely convergent and
\[
|E(t)|\le8|t|^3\sum_{i=1}^n a_i^3 \le8\delta|t|^3\le\frac{t^2}{2}.
\]
Together with the inequality $|e^z-1|\le |z|e^{|z|}$, this implies
\[
|\varphi(t)-e^{-t^2}| \le8\delta|t|^3e^{-t^2/2}, \qquad |t|\le H.
\]
Using $\int_{\bbR}|t|^3e^{-t^2/2}\,\de t=4$, we get
\begin{equation}{\label{eq:Fourier-est-local}}
    \int_{-H}^{H}|\varphi(t)-e^{-t^2}|\,\de t\le32\delta.
\end{equation}
On the complement $|t| > H$, by \eqref{eq:bound-abs-varphi} we get that both $|\varphi(t)|$ and $e^{-t^2}$ are bounded
by $(1+t^2)^{-1}$, so
\begin{equation}{\label{eq:Fourier-est-far}}
    \int_{|t|>H}|\varphi(t)-e^{-t^2}|\,\de t \le\int_{|t|>H}\frac{2\,\de t}{1+t^2} \le\frac4H=64\delta.
\end{equation}
Fourier inversion therefore yields the explicit estimate
\[
\left|f_T(0)-\frac{1}{2\sqrt\pi}\right| \le\frac{48}{\pi}\,\delta=\frac{1}{4\sqrt{\pi}}.
\]

\smallskip
\noindent\emph{Proof of \eqref{e:weighted-chi-square-augmented-density} when $a_1\ge\delta$.}
Write
\[
    T=a_1(Z_1^2-1)+R, \qquad R=\sum_{i=2}^n a_i(Z_i^2-1).
\]
Here $R$ is independent of $Z_1$ and has mean zero and variance $\operatorname{Var}(R)=2\sum_{i=2}^n a_i^2\le2$. Cantelli's inequality (see \cite[Theorem 2.1]{IonKlaassenVanDenHeuvel2023}) gives
\[
    \bbP\left(R\le\frac{a_1}{2}\right) \ge\frac{a_1^2/4}{\operatorname{Var}(R)+a_1^2/4} \ge\frac{\delta^2}{8+\delta^2} =:p_0>0.
\]
In addition, the Chebyshev's inequality gives $\bbP(R<-2/\sqrt{p_0})\le p_0/2$, this implies
\(
\bbP\left(-\frac{2}{\sqrt{p_0}}\le R\le\frac{a_1}{2}\right)\ge\frac{p_0}{2}.
\)
Moreover, the density of $a(Z^2-1)$ is given by
\(
h_a(x)=\frac{\mathbf 1\{x>-a\}}{\sqrt{2\pi a(x+a)}}
\exp\left(-\frac{x+a}{2a}\right).
\)
On the event in the preceding display, we have
$a_1-R\in[a_1/2,a_1+2/\sqrt{p_0}]$ and $a_1\in[\delta,1]$.
It follows that
\[
h_{a_1}(-R)
\ge\frac{1}{\sqrt{2\pi(1+2/\sqrt{p_0})}}
\exp\left(-\frac{1+2/\sqrt{p_0}}{2\delta}\right).
\]
Conditioning on $R$ now gives
\[
f_T(0)=\bbE h_{a_1}(-R)
\ge\frac{p_0}{2\sqrt{2\pi(1+2/\sqrt{p_0})}}
\exp\left(-\frac{1+2/\sqrt{p_0}}{2\delta}\right)>0.\qedhere
\]
\end{proof}

Now we can finish the proof of Lemma~\ref{l:coarea-upper-bound}.
\begin{proof}[Proof of Lemma~\ref{l:coarea-upper-bound}]
  Taking $z = \zeta$ in \eqref{e:hopfield-coarea} yields
  \[
    Z_{\hop,N}^{\sph}(\bG)^{-2} 
    = K_N^{-2} e^{-N\zeta} \det\lt( \zeta I_N-\beta\bW \rt) f_\zeta(N)^{-2}.
  \]
  Recall that $f_\zeta$ is the density of $S(\zeta) = \sum_{i=1}^N a_i(\zeta) Z_i^2$, and that the definition \eqref{e:zeta} of $\zeta$ implies $\sum_{i=1}^N a_i(\zeta) = N$.
  By Lemma~\ref{l:weighted-chi-square-density},
  \[
    f_\zeta(N)^{-2}
    \le C\sum_{i=1}^N a_i(\zeta)^2
    \le C\lt(\sum_{i=1}^N a_i(\zeta)\rt)^2
    = CN^2.
  \]
  Stirling's formula gives $K_N^{-2} \le C N^{-1} e^N$.
  Combining the preceding estimates proves the claim.
\end{proof}

\subsection{Comparison of planted and null characteristic polynomial moments}
\label{ss:gaussian-replacement}

In this subsection we prove Lemma~\ref{l:gaussian-replacement-matrix}.
\begin{lem}\label{l:gaussian-replacement}
  There exists a universal constant $C$ such that if $p : \bbR^k \rightarrow \bbR$ is a polynomial of degree at most $d$, then for $g \sim \cN(0, I_k)$ and every $x \in \bbR^k$,
  \[
    |p(x)|^2 \le C (1 + k + \|x\|^2)^d \bbE[|p(g)|^2].
  \]
\end{lem}
\begin{proof}
  For $n\in \bbZ_{\ge0}$, let $\He_n$ denote the $n$-th probabilists' Hermite polynomial.
  For $\gamma \in \bbZ_{\ge0}^k$, define the normalized multivariate Hermite polynomials
  \[
    h_\gamma(x) = \prod_{i=1}^k \fr{\He_{\gamma_i}(x_i)}{\sqrt{\gamma_i!}}.
  \]
  These form an orthonormal basis for the standard Gaussian measure on $\bbR^k$.
  Write $|\gamma| = \|\gamma\|_1$.
  Expanding $p(x)$ in this basis yields
  \balnn
    p(x)
    = \sum_{|\gamma| \le d}
    a_\gamma h_\gamma(x), \qquad  
    a_\gamma &= \bbE[p(g) h_\gamma(g)].
  \ealnn
  By Cauchy--Schwarz,
  \[
    |p(x)|^2
    \le \lt(
      \sum_{|\gamma| \le d}
      a_\gamma^2
    \rt)\lt(
      \sum_{|\gamma| \le d}
      h_\gamma(x)^2
    \rt)
    = \bbE[|p(g)|^2] \sum_{|\gamma|\le d} h_\gamma(x)^2.
  \]
  The one-dimensional Mehler formula \cite[Eq.~(4.16), p.~51]{janson1997gaussian} gives, for all $t\in (0,1)$,
  \[
    \sum_{n=0}^\infty \fr{t^n}{n!} \He_n(x)^2 = \fr{1}{\sqrt{1-t^2}} \exp\lt(\fr{t}{1+t} x^2\rt).
  \]
  Tensorizing across the coordinates yields
  \[
    \sum_{\gamma \in \bbZ_{\ge0}^k}
    t^{|\gamma|}
    h_\gamma(x)^2
    = (1-t^2)^{-k/2} \exp\lt(\fr{t}{1+t} \|x\|^2\rt).
  \]
  Taking $t = (1 + k + \|x\|^2)^{-1}$ bounds the right-hand side of the last display by an absolute constant.
  Thus
  \[
    \sum_{|\gamma|\le d} h_\gamma(x)^2
    \le t^{-d} \sum_{\gamma \in \bbZ_{\ge0}^k}
    t^{|\gamma|}
    h_\gamma(x)^2
    \le C (1 + k + \|x\|^2)^d. \qedhere
  \]
\end{proof}
\begin{proof}[Proof of Lemma~\ref{l:gaussian-replacement-matrix}]
  Write $\bG = [\bC, \bA_0]$, where $\bA_0$ is standard Gaussian under either $\bbE_{\star,q}$ or $\bbE_\nul$.
  Let $\bbE'_{\star,q}$ and $\bbE'_\nul$ denote these expectations conditional on $\bA_0$.

  For any fixed $\bA_0$, $\det (z - \bW)^k$ is a polynomial in the entries of $\bC$ of degree at most $4k$.
  Lemma~\ref{l:gaussian-replacement} gives that for fixed $\bC$,
  \[
    |\det (z - \bW)|^{2k}
    \le C \lt(1 + 2M + \|\bC\|_F^2 \rt)^{4k} \bbE'_\nul[|\det (z - \bW)|^{2k}].
  \]
  Taking expectation over $\bbE'_{\star,q}$ yields that for $C$ depending only on $\phi,\alpha$, and (in the case $\star = \hop$) $J$,
  \[
    \bbE'_{\star,q}[|\det (z - \bW)|^{2k}]
    \le CN^{4k} \bbE'_\nul[|\det (z - \bW)|^{2k}].
  \]
  The conclusion follows by taking expectation over $\bA_0$.
\end{proof}

\subsection{Characteristic polynomial moment bound}
\label{ss:characteristic-polynomial-moment}

In this subsection we prove Lemma~\ref{l:characteristic-polynomial-moment}.
Recall that for $\gamma>0$, let $\nu_{\MP,\gamma}$ denote the
Marchenko--Pastur distribution given by
\[
\nu_{\MP,\gamma}(\de x)
=
\frac{\sqrt{(\ell_+(\gamma)-x)(x-\ell_-(\gamma))}}
{2\pi x}
\mathbf 1_{[\ell_-(\gamma),\ell_+(\gamma)]}(x)\de x
+(1-\gamma)_+\delta_0(\de x),
\]
where $\ell_\pm(\gamma)=(\sqrt{\gamma}\pm1)^2$.
Define its logarithmic potential by
\[
L_\gamma(z)
=\int\log|z-x|\,\nu_{\MP,\gamma}(\de x).
\]
For $w\in \mathbb C$ with $\Im w>0$, we define its Stieltjes transform as
\beq\label{e:wishart-dual-saddles}
u_\gamma(w)
=-\int\frac{\nu_{\MP,\gamma}(\de x)}{w-x}
=
\frac{\gamma-1-w+
\sqrt{(w-\ell_-(\gamma))(w-\ell_+(\gamma))}}{2w}.
\eeq
For real $z\ne0$, define $u_\gamma(z):=\lim_{\eta\downarrow0}u_\gamma(z+\mathrm i\eta)$. We obtain
\begin{equation}\label{e:s-tranform-real-z}
u_\gamma(z)=
\begin{cases}
\displaystyle
\frac{\gamma-1-z-
\sqrt{(z-\ell_-)(z-\ell_+)}}{2z},
& z<\ell_-,\\[8pt]
\displaystyle
\frac{\gamma-1-z+
\mathrm i\sqrt{(\ell_+-z)(z-\ell_-)}}{2z},
& \ell_-\le z\le\ell_+,\ z>0,\\[8pt]
\displaystyle
\frac{\gamma-1-z+
\sqrt{(z-\ell_-)(z-\ell_+)}}{2z},
& z>\ell_+.
\end{cases}
\end{equation}
When $\gamma>1$, this definition extends
continuously to $z=0$, with
$u_\gamma(0)=1/(\gamma-1)$.

\begin{lem}\label{l:log-potential-s-transform}
   For $z\in \mathbb R,z\neq 0$, the identity
\beq\label{e:wishart-saddle-potential}
L_\gamma(z)
=-1-\log|u_\gamma(z)|-z\Re u_\gamma(z)
+\gamma\log|1+u_\gamma(z)|
\eeq
holds for every $\gamma>0$ and extends to $z=0$ when $\gamma>1$. 
\end{lem}
\begin{proof}
Recall that $u_\gamma(w)$ satisfies the equation
\begin{equation}\label{e:saddle-eq}
-\frac{1}{u_\gamma(w)}-w
+\frac{\gamma}{1+u_\gamma(w)}=0.
\end{equation}
For $\Im w>0$, we define
\[
F_\gamma(w)
:=-1-\log(-u_\gamma(w))-wu_\gamma(w)
+\gamma\log(1+u_\gamma(w)),
\]
where we use the principal branches of the logarithms. Note that $\Im u_\gamma(w)>0$ because $\Im w>0$, the function $F_\gamma$
is analytic in the upper half-plane. By \eqref{e:saddle-eq} we have $F_\gamma'(w)=-u_\gamma(w)$. Moreover, note that $G_\gamma(w)=-u_\gamma(w)$ holds by definition. Thus, $F_\gamma-G_\gamma$ is constant in the upper half-plane. To determine this constant, observe that
$u_\gamma(w)=-1/w+O(w^{-2})$ as $w\to\infty$.
Consequently, as $w\to\infty$ we have $F_\gamma(w)=G_\gamma(w)=\log w+o(1).$ It follows that $F_\gamma(w)=G_\gamma(w)$
throughout the upper half-plane.
Taking real parts of $F_\gamma(w)=G_\gamma(w)$
and passing to the boundary yields
\[
L_\gamma(z)
=-1-\log|u_\gamma(z)|
-z\Re u_\gamma(z)
+\gamma\log|1+u_\gamma(z)|.
\]
Finally, when $\gamma>1$, the support of
$\nu_{\MP,\gamma}$ is bounded away from zero,
and by \eqref{e:s-tranform-real-z} $u_\gamma(z)$ extends continuously to
$z=0$, with $u_\gamma(0)=1/(\gamma-1)$.
Hence the identity also holds at $z=0$.
\end{proof}

Next we record a lemma which is useful for the proof of Lemma~\ref{l:wishart-characteristic-polynomial}.
\begin{lem}
\label{l:circular-phase}
For $\gamma\ge1$ and $z\in\bbR$, define
\[
H_{\gamma,z}(r,t)
=-\log r-zrt
+\frac{\gamma}{2}\log(1+r^2+2rt),
\qquad r>0,\quad t\in[-1,1],
\]
with the convention $\log 0=-\infty$.
Then
\[
\inf_{r>0}\max_{-1\le t\le1}
H_{\gamma,z}(r,t)
=1+L_\gamma(z).
\]
For $(\gamma,z)\ne(1,0)$, the infimum is
attained at $r=|u_\gamma(z)|$.
When $\gamma=1$ and $z=0$, the infimum
is approached as $r\to\infty$.
\end{lem}

\begin{proof}
Let $F(s):=\max_{-1\le t\le1}H_{\gamma,z}(e^s,t),s\in\bbR$. We first show that $F$ is convex. Indeed, the function
\[
h(u)=-\log|u|-z\Re u+\gamma\log|1+u|
\]
is subharmonic on $\bbC\setminus\{0\}$, Writing $u=e^{s+i\theta}$, we have
$h(u)=H_{\gamma,z}(e^s,\cos\theta)$.
Since $\cos\theta$ ranges over $[-1,1]$, thus we have $F(s)=\max_{|u|=e^s}h(u)$. For $s_1<s_2$, let $g$ be the affine function
satisfying $g(s_i)=F(s_i)$, $i=1,2$.
The function $u\mapsto g(\log|u|)$ is harmonic
on the annulus $e^{s_1}<|u|<e^{s_2}$
and dominates $h$ on its boundary.
The maximum principle therefore gives
$F(s)\le g(s)$ for $s\in[s_1,s_2]$. This proves convexity.

We next identify a stationary point of $F$.
Recall that
\[
\partial_tH_{\gamma,z}(r,t)
=-zr+\frac{\gamma r}{1+r^2+2rt},
\qquad
\partial_t^2H_{\gamma,z}(r,t)
=-\frac{2\gamma r^2}{(1+r^2+2rt)^2}<0
\]
wherever the derivatives are defined.
Thus $H_{\gamma,z}(r,\cdot)$ is strictly
concave. Suppose $(\gamma,z)\ne(1,0)$
and take $r=|u_\gamma(z)|$. If $z\in[\ell_-(\gamma),\ell_+(\gamma)]$
and $z>0$, we have $r=|u_\gamma(z)|=z^{-1/2}$ by \eqref{e:s-tranform-real-z}. This gives $\partial_tH_{\gamma,z}(r,t_*)=0$ where $t_*=\frac{\gamma-1-z}{2\sqrt z}\in[-1,1]$. 
By concavity, we have $t_*$ is the unique maximizer
of $H_{\gamma,z}(r,\cdot)$.

If $z<\ell_-(\gamma)$, then $u_\gamma(z)>0$, let $t_*=1$.
The saddle-point equation \eqref{e:saddle-eq} gives
\[
\partial_tH_{\gamma,z}(r,1)
=1-\frac{\gamma r^2}{(1+r)^2}\ge0.
\]
For $\gamma>1$, the inequality follows from the monotonicity
of $u_\gamma$ and
\(
0<r\le u_\gamma(\ell_-(\gamma))
=\frac1{\sqrt\gamma-1}.
\)
For $\gamma=1$, the inequality is immediate. Thus $t_*=1$ is the unique maximizer. If $z>\ell_+(\gamma)$, then $u_\gamma(z)<0$, thus we have $r=-u_\gamma(z)$. Similarly we can show $t_*=-1$ is the unique maximizer of $H_{\gamma,z}(r,\cdot)$ .

In all three cases, $t_*$ is the unique maximizer
of $H_{\gamma,z}(r,\cdot)$. Moreover, writing
$u=u_\gamma(z)$, we have
$r=|u|$ and $rt_*=\Re u$. Hence
\[
\begin{aligned}
r\partial_rH_{\gamma,z}(r,t_*)
&=-1-z\Re u
+\gamma\Re\left(\frac{u}{1+u}\right)=0,
\end{aligned}
\]
where the last equality follows by multiplying
\eqref{e:saddle-eq} by $u$ and taking real parts.
Consequently, we get
\(
F'(\log r)
=r\partial_rH_{\gamma,z}(r,t_*)=0.
\)
By convexity, $F$ attains its global minimum
at $\log r$. Lemma~\ref{l:log-potential-s-transform}
then gives
\[
\begin{aligned}
\inf_{s\in\bbR}F(s)
=H_{\gamma,z}(r,t_*)=-\log|u_\gamma(z)|
-z\Re u_\gamma(z)
+\gamma\log|1+u_\gamma(z)|=1+L_\gamma(z).
\end{aligned}
\]

 Finally, when $\gamma=1$ and $z=0$, we have
\(
\max_{-1\le t\le1}H_{1,0}(r,t)
=\log(1+r^{-1})\longrightarrow0
=1+L_1(0)
\)
as $r\to\infty$, completing the proof.
\end{proof}

\begin{lem}
\label{l:wishart-characteristic-polynomial}
Let $\bX\in\bbR^{m\times n}$ have independent $\mathcal{N}(0,1)$ entries,
where $m\ge n\ge1$, and set $\bW=\bX^\top\bX/n$ and $\gamma=m/n$.
There is a universal constant $C_0<\infty$ such that, for every $z\in\bbR$,
\beq\label{e:wishart-char-poly-fourth}
\bbE\bigl[|\det(zI_n-\bW)|^4\bigr]
\le C_0n^{14}\exp\{4nL_\gamma(z)\}.
\eeq
\end{lem}
\begin{proof}
For simplicity we write $D_n(z)=\det(zI_n-\bW)$. And for $\bu\in\mathbb R^4$ we let $\Delta(\bu)=\prod_{1\le j<k\le4}(u_j-u_k)$.
Recall that the eigenvalues of $\bX^\top\bX$ have joint
density proportional to
\[
\prod_{j=1}^n x_j^{(m-n-1)/2}e^{-x_j/2}
\prod_{1\le j<k\le n}|x_j-x_k|.
\]
Thus, they form a Laguerre $\beta$-ensemble with
$\beta=1$ and parameter $a=(m-n-1)/2$.
Applying \cite[Eq.~(4.1)]{ShenForrester25}
with matrix size $n$ and $p=2$ (corresponding
to the fourth moment) yields
    \begin{equation}\label{eq:integral-representation-4thmoment}
\begin{aligned}
\bbE[D_n(z)^4]
={}&n^{-4n}
\frac{C_{1,n}[x^{a+4}e^{-x/2}]}
     {C_{1,n}[x^a e^{-x/2}]\widetilde C_{4,4}[\omega]}\\
&\quad\times\int_{[-1/2,1/2]^4}
\prod_{j=1}^4\left(e^{-nz e^{2\pi i\theta_j}}
\omega(\theta_j)\right)
\prod_{j<k}|e^{2\pi i\theta_j}-e^{2\pi i\theta_k}|^4\,
\de\theta_1\cdots\de\theta_4\,,
\end{aligned}
\end{equation}
where
\[
\omega(\theta)=
e^{\pi i(m-2n)\theta} \left|1+e^{2\pi i\theta}\right|^{m}.
\]
\[
C_{\beta,N}[w]
=\int_{(0,\infty)^N}\prod_{j=1}^N w(x_j)
\prod_{1\le j<k\le N}|x_j-x_k|^\beta\,
\de x_1\cdots\de x_N.
\]
\[
\widetilde C_{4,4}[\omega] = \int_{[-1/2,1/2]^4} \prod_{j=1}^{4} \omega(\theta_j) \prod_{j<k} \left| e^{2\pi i\theta_j}- e^{2\pi i\theta_k} \right|^4 \de\theta_1\cdots\de\theta_4.
\]
We next evaluate $C_{\beta,N}[w]$ and $\widetilde C_{4,4}[\omega] $. 
By the Laguerre Selberg integral \cite{Selberg1944,ForresterWarnaar2008},
\[
C_{1,n}[x^a e^{-x/2}]
=
2^{mn/2}
\prod_{j=0}^{n-1}
\frac{
\Gamma((j+3)/2)
\Gamma((m-n+j+1)/2)
}{
\Gamma(3/2)
},
\]
\[
C_{1,n}[x^{a+4} e^{-x/2}]
=
2^{mn/2+4n}
\prod_{j=0}^{n-1}
\frac{
\Gamma((j+3)/2)
\Gamma((m-n+j+9)/2)
}{
\Gamma(3/2)
}.
\]
To evaluate the circular normalization, observe that
\(\omega(\theta)=(1+e^{2\pi i\theta})^{m-n}(1+e^{-2\pi i\theta})^n.\) The Morris integral \cite{Morris1982,ForresterWarnaar2008} with parameters $N=4$, $A=m-n$, $B=n$, and $q=2$,
then gives
\[
\widetilde C_{4,4}[\omega]
=
\prod_{j=0}^{3}
\frac{
\Gamma(m+2j+1)\Gamma(2j+3)
}{
\Gamma(m-n+2j+1)
\Gamma(n+2j+1)\Gamma(3)
}.
\]
Combining these identities, we obtain
\[
\begin{aligned}
c_n:=&n^{-4n}
\frac{C_{1,n}[x^{a+4}e^{-x/2}]}
{C_{1,n}[x^ae^{-x/2}]
\widetilde C_{4,4}[\omega]}=
n^{-4n}
\prod_{j=0}^{3}
\frac{\Gamma(n+2j+1)\Gamma(3)}
{\Gamma(2j+3)}
\le Cn^{14}e^{-4n},
\end{aligned}
\]
where the last inequality follows from
Stirling's formula. Then we set
$u_j=e^{2\pi i\theta_j}$. Note that on the unit circle,
\[
\prod_{1\leq j<k\leq 4}|u_j-u_k|^4
=\Delta(\bu)^4\prod_{j=1}^4u_j^{-6}.
\]
By change of variables $u_j=e^{2\pi i\theta_j}$ in \eqref{eq:integral-representation-4thmoment}, we then obtain
\beq
\bbE[D_n(z)^4]
=
\frac{c_n}{(2\pi i)^4}
\oint_{|u_1|=1}\!\cdots\!\oint_{|u_4|=1}
\Delta(\bu)^4
\prod_{j=1}^4
\frac{e^{-nzu_j}(1+u_j)^m}{u_j^{n+7}}
\,\de u_1\cdots\de u_4.\notag
\eeq
Note that the integrand is single-valued and
holomorphic in each variable away from $0$. By Cauchy's integral theorem we get that
\beq\label{e:wishart-fourth-dual}
\bbE[D_n(z)^4]
=
\frac{c_n}{(2\pi i)^4}
\oint_{|u_1|=r}\!\cdots\!\oint_{|u_4|=r}
\Delta(\bu)^4
\prod_{j=1}^4
\frac{e^{-nzu_j}(1+u_j)^m}{u_j^{n+7}}
\,\de u_1\cdots\de u_4.
\eeq
It remains to bound the contour integral. Recall the definition $H_{\gamma,z}(r,t)$ in Lemma~\ref{l:circular-phase}.
For $u=re^{i\theta}$, observe that
$\Re u=r\cos\theta$ and
$|1+u|^2=1+r^2+2r\cos\theta$.
Consequently,
\[
H_{\gamma,z}(r,\cos\theta)
=-\log|u|-z\Re u+\gamma\log|1+u|.
\]
Since $m=n\gamma$, it follows that
\[
\left|
\frac{e^{-nzu}(1+u)^m}{u^{n+7}}
\right|
=r^{-7}\exp\{nH_{\gamma,z}(r,\cos\theta)\}.
\]
Also, $|\Delta(\bu)|^4\le(2r)^{24}$. Plugging the estimates in \eqref{e:wishart-fourth-dual} and take infimum over $r>0$ yields
\beq\label{e:wishart-circular-bound}
\bbE[D_n(z)^4]
\le 2^{24}c_n
\exp\left\{4n\inf_{r>0}\max_{-1\le t\le1}H_{\gamma,z}(r,t)\right\}.
\eeq
By Lemma~\ref{l:circular-phase}, we have
\[
\bbE[D_n(z)^4]
\le 2^{24}c_n e^{4n(1+L_\gamma(z))}
\le C_0 n^{14}e^{4nL_\gamma(z)}.
\]
Since $D_n(z)$ is real for $z\in\bbR$, its fourth power equals
$|D_n(z)|^4$.
\end{proof}

\begin{proof}[Proof of Lemma~\ref{l:characteristic-polynomial-moment}]
Set $\alpha_N=M/N$, where $M=\lfloor\alpha N\rfloor$.
We first prove the estimate at aspect ratio $\alpha_N$.
If $M\ge N$, Lemma~\ref{l:wishart-characteristic-polynomial} gives
\beq\label{e:wishart-exact-ratio-moment}
\bbE_\nul[|\det(zI_N-\bW)|^4]
\le C_0N^{14}\exp\{4NL_{\alpha_N}(z)\}.
\eeq
If $1\le M<N$, let $\widetilde{\bW}=\bG\bG^\top/M$.
The equality of the nonzero eigenvalues of $\bG^\top\bG$ and
$\bG\bG^\top$ gives
\[
\det(zI_N-\bW)
=z^{N-M}\alpha_N^M
\det\bigl((z/\alpha_N)I_M-\widetilde{\bW}\bigr).
\]
The corresponding relation between Marchenko--Pastur measures yields,
for $z\ne0$,
\[
L_{\alpha_N}(z)
=(1-\alpha_N)\log|z|+\alpha_N\log\alpha_N
+\alpha_N L_{1/\alpha_N}(z/\alpha_N).
\]
Applying Lemma~\ref{l:wishart-characteristic-polynomial} to
$\widetilde{\bW}$ and using $M^{14}\le N^{14}$ proves
\eqref{e:wishart-exact-ratio-moment} also in this case.

It remains to replace $\alpha_N$ by $\alpha$.
Using
\(
\partial_\gamma L_\gamma(z)
=\log|1+u_\gamma(z)|
\ge-\log(1+\gamma^{-1/2})
\)
by \eqref{e:wishart-dual-saddles},
together with $0\le\alpha-\alpha_N<1/N$, we obtain
\(
L_{\alpha_N}(z)
\le L_\alpha(z)
+\frac1N\log(1+\alpha_N^{-1/2}).
\)
Consequently, \eqref{e:wishart-exact-ratio-moment} gives
\[
\bbE_\nul[|\det(zI_N-\bW)|^4]
\le C_0N^{14}(1+\alpha_N^{-1/2})^4
e^{4NL_\alpha(z)}.
\]
For fixed $\alpha>0$, the prefactor is bounded
by $C_\alpha N^{14}$. If $\alpha\ge1+\eta$,
it is bounded by $16C_0N^{14}$ for all
sufficiently large $N$, uniformly in $\alpha$.
The bound extends to $z=0$ by continuity
when $\alpha\ge1$.
\end{proof}

\subsection{Inward deviations at regular Wishart soft edges}
\label{ss:inward-deviations}

In this subsection we prove Lemmas~\ref{l:zeta-deviation}--\ref{l:wishart-inward-deviation-bbp}.
\begin{proof}[Proof of Lemma~\ref{l:zeta-deviation}]
  This is proved by a concentration of spectral traces argument.
  Let $\eps = \eps(\phi, \alpha, u) > 0$ be small, to be chosen below.
  Define
  \[
    f_\eps(x) = \fr{1}{\max(\zeta_* - u - x, \eps)}.
  \]
  Recall $\bW = \fr1N \bG^\top \bG \in \bbR^{N\times N}$, $\bW_0 = \fr1N \bA_0^\top\bA_0 \in \bbR^{(N-2) \times (N-2)}$, and define $\Tr f_\eps(\beta \bW)$ by the standard spectral calculus.
  Then $\Tr f_\eps(\beta \bW) \ge \Tr f_\eps(\beta \bW_0)$ by interlacing, so
  \beq\label{e:clipped-resolvent}
    \bbP_{\star,q}(\zeta < \zeta_* - u)
    \le \bbP_{\star,q}\lt( \Tr f_\eps(\beta \bW) < N \rt)
    \le \bbP\lt( \Tr f_\eps(\beta \bW_0) < N \rt).
  \eeq
  Since $f_\eps$ is bounded, by the Marchenko--Pastur theorem,
  \[
    \bbE \Tr f_\eps(\beta \bW_0) = N \int f_\eps(\beta x)\,\nu_{\MP,\alpha}(\de x) + o(N).
  \]
  One directly checks that $\zeta_*-\beta x>0$ for $\nu_{\MP,\alpha}$-almost every $x$.
  For such $x$, we also have
  \[
    \lim_{\eps\downarrow 0} f_\eps(\beta x) > (\zeta_*-\beta x)^{-1},
  \]
  where the left-hand side limit is increasing and may be infinite.
  By the monotone convergence theorem,
  \[
    \lim_{\eps\downarrow 0}
    \int f_\eps(\beta x)\,\nu_{\MP,\alpha}(\de x) 
    > \int (\zeta_*-\beta x)^{-1} \nu_{\MP,\alpha}(\de x)
    = 1,
  \]
  where the final equality is a direct calculation.
  We can thus pick $\eps, \delta > 0$ depending only on $\phi,\alpha,u$ such that 
  $\int f_\eps(\beta x)\,\nu_{\MP,\alpha}(\de x) \ge 1 + \delta$.
  It follows that
  \[
    \bbE \Tr f_\eps(\beta \bW_0) \ge (1+\delta/2) N.
  \]
  By \cite[Lemma~1.2(b)]{GuionnetZeitouni00concentration}, $\bA_0 \mapsto \Tr f_\eps(\beta \bW_0)$ is $C$-Lipschitz for some $C = C(\phi,\alpha,u)$.
  The result follows by Gaussian concentration.
\end{proof}
Lemmas~\ref{l:wishart-inward-deviation}--\ref{l:wishart-inward-deviation-bbp} both depend through the following representation.
\begin{lem}[{Golub--Kahan bidiagonalization, \cite[proof of Theorem~3.1]{DumitriuEdelman02}}]
  \label{l:golub-kahan}
  Let $\bX = [\bv,\bA] \in \bbR^{m\times n}$, where $\bA \in \bbR^{m\times(n-1)}$ has i.i.d. standard Gaussian entries and $\bv \in \bbR^m$ is independent of $\bA$.
  Let $\bB \in \bbR^{m\times n}$ be the upper bidiagonal matrix with $\bB_{1,1} = d_1 = \|\bv\|$ and independent entries
  \balnn
    \bB_{j,j} &= d_j \sim \chi_{m-j+1}\,, & 
    \bB_{j,j+1} &= e_j \sim \chi_{n-j}\,,
  \ealnn
  for respectively $2\le j\le \min(m,n)$ and $1\le j\le \min(m,n-1)$, all independent of $\bv$.
  All other entries of $\bB$ are zero.
  Then, there exists a joint distribution on $(\bB, \bU, \bV)$, where $\bU \in \bbR^{m\times m}$ and $\bV \in \bbR^{n\times n}$ are orthogonal, such that $\bX = \bU^\top \bB \bV$.
\end{lem}

\begin{proof}[Proof of Lemma~\ref{l:wishart-inward-deviation}]
  Below we set $c_0, C_0$ sufficiently small and large depending only on $\alpha$, such that $c_0 C_0$ is sufficiently small and $c_0^2 C_0^3$ is sufficiently large; then we set $c$ small depending only on $\alpha, C_0$.
  If $0 \le u \le C_0 N^{-2/3}$, the result follows by taking $c$ small enough that $2e^{-cN^2u^3} \ge 2e^{-cC_0^3} \ge 1$.
  From now on we assume $C_0 N^{-2/3} \le u \le 1$.

  We apply Lemma~\ref{l:golub-kahan} with $m = M$, $n = N-2$, and $\bX = \bA_0$ (so that $d_1 = \|\bv\| \sim \chi_m$).
  Let $\bB$ be defined therein.
  For any fixed nonzero $\bw, \bv \in \bbR^n$, 
  \balnn
    \lambda_{\max}(\bW_0)
    \stackrel{d}{=} \lambda_{\max}\lt(\fr1N \bB^\top \bB\rt)
    &\ge \fr{\|\bB \bw\|^2}{N \|\bw\|^2}, &
    \lambda_{\min}(\bW_0)
    \stackrel{d}{=} \lambda_{\min}\lt(\fr1N \bB^\top \bB\rt)
    &\le \fr{\|\bB \bv\|^2}{N \|\bv\|^2}.
  \ealnn
  We will fix test vectors $\bw, \bv$ and argue that for $\ell_\pm = (\alpha^{1/2} \pm 1)^2$,
  \baln\label{e:wishart-inward-deviation-goal}
    \bbP\lt(\fr{\|\bB \bw\|^2}{N \|\bw\|^2} \le \ell_+ - u\rt) &\le 2e^{-cN^2u^3}, &
    \bbP\lt(\fr{\|\bB \bv\|^2}{N \|\bv\|^2} \ge \ell_- + u\rt) &\le 2e^{-cN^2u^3}.
  \ealn
  This implies the result.
  Set $L = \lfloor c_0 uN \rfloor \le N$.
  Then define $\bw, \bv$ by
  \balnn
    w_j &= \sin \frac{\pi (j-1)}{L}, & 
    v_j &= (-1)^j w_j
  \ealnn
  for $1\le j\le L+1$, and $w_j = v_j = 0$ otherwise.
  We next estimate $\bbE \|\bB \bw\|^2$ and $\bbE \|\bB \bv\|^2$.
  By the standard estimate $\sqrt{r-1} \le \bbE [\chi_r] \le \sqrt{r}$,
  \baln
    \nonumber
    \bbE \|\bB \bw\|^2
    &= \sum_{j=1}^L \bbE[(d_j w_j + e_j w_{j+1})^2]
    \ge \sum_{j=1}^L \lt[(m-j+1) w_j^2 + (n-j) w_{j+1}^2 + 2 \sqrt{(m-j)(n-j-1)} w_j w_{j+1} \rt] \\
    \nonumber
    &\ge \sum_{j=1}^L \lt[(m-j+1) w_j^2 + (n-j) w_{j+1}^2 + \sqrt{(m-j)(n-j-1)} (w_j^2 + w_{j+1}^2) \rt] \\
    \nonumber
    &\qquad - \sum_{j=1}^L \sqrt{(m-j)(n-j-1)} (w_j - w_{j+1})^2 \\
    \label{e:wishart-upper-inward-deviation}
    &\ge \lt((\sqrt{\alpha} + 1)^2 N - C_0 L\rt) \|\bw\|^2
    - \sqrt{MN} \sum_{j=1}^L (w_j - w_{j+1})^2\,,
  \ealn
  if $C_0$ is sufficiently large depending only on $\alpha$.
  An identical argument yields
  \beq\label{e:wishart-lower-inward-deviation}
    \bbE \|\bB \bv\|^2
    \le \lt((\sqrt{\alpha} - 1)^2 N + C_0 L\rt) \|\bw\|^2
    + \sqrt{MN} \sum_{j=1}^L (w_j - w_{j+1})^2\,.
  \eeq
  Note that $\|\bw\|^2 = \|\bv\|^2$ and
  \[
    \sum_{j=1}^L (w_j - w_{j+1})^2
    = 2\lt(1-\cos\fr{\pi}{L}\rt) \|\bw\|^2
    \le \fr{\pi^2}{L^2} \|\bw\|^2.
  \]
  Thus \eqref{e:wishart-upper-inward-deviation}, \eqref{e:wishart-lower-inward-deviation} imply
  \baln\label{e:bidiagonal-expectation}
    \fr{\bbE \|\bB \bw\|^2}{N\|\bw\|^2} &\ge (\sqrt{\alpha} + 1)^2 - \fr{C_0 L}{N} - \fr{\pi^2 \sqrt{\alpha}}{L^2}, &
    \fr{\bbE \|\bB \bv\|^2}{N\|\bv\|^2} &\le (\sqrt{\alpha} - 1)^2 + \fr{C_0 L}{N} + \fr{\pi^2 \sqrt{\alpha}}{L^2}.
  \ealn
  Since we set $c_0, C_0$ such that $c_0C_0$ is sufficiently small and $c_0^2 C_0^3$ sufficiently large, and assume $u \ge C_0 N^{-2/3}$, we can ensure 
  \beq\label{e:bidiagonal-error-bound}
    \max\lt\{
      \fr{C_0 L}{N},
      \fr{\pi^2 \sqrt{\alpha}}{L^2}
    \rt\}
    \le \fr{u}{6}.
  \eeq
  Finally, if we write $d_j^2 = \sum_{k=1}^{m-j+1} \bg_{j,k}^2$ and $e_j^2 = \sum_{k=1}^{n-j} \tilde\bg_{j,k}^2$ for $\bg, \tilde \bg$ with i.i.d. standard Gaussian entries, then one easily checks the maps
  \balnn
    F_+(\bg, \tilde\bg) &= \fr{\|\bB \bw\|}{\sqrt{N}\|\bw\|}, & 
    F_-(\bg, \tilde\bg) &= \fr{\|\bB \bv\|}{\sqrt{N}\|\bv\|},
  \ealnn
  are $2/\sqrt{NL}$-Lipschitz.
  It follows that $\Var(F_{\pm}) \le 4 / NL \le u/6$.
  Together with \eqref{e:bidiagonal-expectation} and \eqref{e:bidiagonal-error-bound} we have
  \balnn
    (\bbE F_+)^2 &\ge (\sqrt{\alpha} + 1)^2 - \fr u2, & 
    (\bbE F_-)^2 &\le (\sqrt{\alpha} - 1)^2 + \fr u2.
  \ealnn
  By Gaussian concentration of measure, for $c_1$ be sufficiently small depending only on $\alpha$,
  \[
    \bbP\lt( (F_+)^2 \le (\sqrt{\alpha} + 1)^2 - u \rt)
    \le \bbP\lt( F_+ - \bbE F_+ \le -c_1 u \rt)
    \le \exp\lt(
      -\fr{NL}{8} \cdot (c_1 u)^2
    \rt)
    \le e^{-cN^2u^3}\,.
  \]
  This proves the first estimate of \eqref{e:wishart-inward-deviation-goal}, and the second estimate is analogous.
\end{proof}
\begin{proof}[Proof of Lemma~\ref{l:wishart-inward-deviation-bbp}]
  Assume without loss $q>0$, so that $q\in [q_N/2, q_0]$.
  Throughout this proof we take $q_0\le1$ sufficiently small, and all constants depend only on $\phi,\alpha$.
  We write $\bbP$, $\bbE$ for $\bbP_{\star,q}$, $\bbE_{\star,q}$.

  Form $\bX = [\bC_1,\bA_0] \in \bbR^{M \times (N-1)}$ and set $\mu'_+ = \lambda_{\max} (\beta \bX^\top \bX / N)$.
  By interlacing, $\mu_+ \ge \mu'_+$.

  Applying Lemma~\ref{l:golub-kahan} to $\bX$ with $m=M$, $n=N-1$, and $\bv=\bC_1$ provides an upper bidiagonal matrix $\bB \in \bbR^{M \times (N-1)}$, where $d_1 = \|\bC_1\|$ and the remaining nonzero entries are independent with the chi distributions stated therein.
  For any nonzero $\bw \in \bbR^M$ supported on the first $\min(M, N-1)$ coordinates,
  \[
    \mu'_+ \stackrel{d}{=} \lambda_{\max}\lt(\fr{\beta \bB^\top \bB}{N}\rt)
    \ge \fr{\beta \|\bB^\top \bw\|^2}{N\|\bw\|^2}.
  \]
  We will fix such a test vector $\bw$ and argue that for all $0\le u\le 1$,
  \beq\label{e:wishart-inward-deviation-bbp-goal}
    \bbP\lt(
      \fr{\beta \|\bB^\top \bw\|^2}{N\|\bw\|^2}
      \le \mu_* + cq^2 - u
    \rt)
    \le 2\exp\lt(-cN\fr{(u+q^2)^2}{q}\rt).
  \eeq
  Let $C$ be a sufficiently large constant depending only on $\alpha$ and $L=\lceil C/q\rceil \le \fr12 \min(M, N-1)$.
  Let $\eps = \sgn(s-1) = \sgn(\beta)$; note that $s \neq 1$ under the assumption $\theta = 1$.
  Define
  \[
    w_j = \eps^j (L+1-j)
  \]
  for $1\le j\le L$, and $w_j = 0$ for $j > L$.
  Note that
  \beq\label{e:w-norm}
    \|\bw\|^2 = \sum_{j=1}^L j^2
    \in [L^3 / 3, L^3].
  \eeq
  Also recall from the proof of Lemma~\ref{l:good-gram} the Taylor expansion
  \[
    \bbE[d_1^2] = \bbE \|\bC_1\|^2 = M \lt(s + s(s-1) q + O(q^2)\rt).
  \]
  Thus, for sufficiently small $q_0$,
  \[
    \eps(\bbE d_1^2 - Ms) \ge M \fr{s|s-1|}{2}q.
  \]
  By the standard estimate $\sqrt{r-1}\le\bbE\chi_r\le\sqrt r$, and for constant $C_1$ depending only on $\alpha$ changing from line to line,
  \balnn
    \beta\bbE \|\bB^\top\bw\|^2
    &= \beta \bbE [d_1^2] L^2
    + \beta \sum_{j=1}^L \bbE[(e_jw_j+d_{j+1}w_{j+1})^2]\\
    &\ge \beta \alpha N \lt(s + \fr{s(s-1)}{2}q\rt) L^2
    + \beta N \sum_{j=1}^L(w_j+\sqrt{\alpha} w_{j+1})^2 - C_1L\|\bw\|^2 \\
    &= \mu_* N\|\bw\|^2
    + \fr{\alpha \beta s(s-1) L^2 Nq}{2} 
    - \sqrt{\alpha} |\beta| LN
    - C_1L\|\bw\|^2,
  \ealnn
  where we have used $s = 1 + \eps/\sqrt\alpha$ and $\mu_* = \beta(\sqrt\alpha+\eps)^2$ at $\theta = 1$.
  Thus
  \[
    \fr{\beta\bbE \|\bB^\top\bw\|^2}{N\|\bw\|^2}
    \ge \mu_* 
    + \fr{\alpha \beta s(s-1) L^2 q}{2 \|\bw\|^2}
    - \fr{\sqrt{\alpha} |\beta| L}{\|\bw\|^2}
    - \fr{C_1 L}{N}.
  \]
  By the estimate \eqref{e:w-norm} and definition of $L$, we can find a small constant $c_1 > 0$ such that
  \balnn
    \fr{\alpha \beta s(s-1) L^2 q}{2 \|\bw\|^2}
    &\ge 5c_1q^2, &
    \max\lt\{
      \fr{\sqrt{\alpha} |\beta| L}{\|\bw\|^2},
      \fr{C_1 L}{N}
    \rt\} \le c_1q^2.
  \ealnn
  We have thus shown
  \beq\label{e:inward-deviation-bbp-mean}
    \fr{\beta\bbE \|\bB^\top\bw\|^2}{N\|\bw\|^2}
    \ge \mu_* + 3c_1q^2.
  \eeq
  We next bound the fluctuations $\beta\bbE \|\bB^\top\bw\|^2 / N \|\bw\|^2$.
  Write
  \balnn
    \fr{\beta\|\bB^\top\bw\|^2}{N\|\bw\|^2}
    & = \frac{\beta L^2d_1^2}{N\|\bw\|^2}+\beta F^2, & 
    F &= \frac{1}{\sqrt N\|\bw\|}
       \left(\sum_{j=1}^L(e_jw_j+d_{j+1}w_{j+1})^2\right)^{1/2}.
  \ealnn
  Recall $d_1^2 = \|\bC_1\|^2$.
  As in the proof of Lemma~\ref{l:good-gram}, the $(\bC_{1,j}^2 : 1\le j\le M)$ are independent and have uniformly bounded exponential moments for $q\le q_0$.
  Bernstein's inequality and \eqref{e:w-norm} give, for some constant $c_2$ and all $0\le t\le 1$,
  \beq\label{e:d1-conc}
    \bbP\lt(
      \fr{\beta L^2}{N\|\bw\|^2}(d_1^2-\bbE d_1^2)
      \le - t
    \rt)
    \le e^{-c_2 NL t^2}.
  \eeq
  Like in the proof of Lemma~\ref{l:wishart-inward-deviation}, we can realize the $d_j^2$ and $e_j^2$ in the definition of $F$ as sums of squares of i.i.d. standard Gaussians.
  Then $F$ is a Lipschitz function of these Gaussians, with Lipschitz constant $\sqrt{2L^2/(N\|\bw\|^2)} \le \sqrt{6/(NL)}$ by \eqref{e:w-norm}.
  It follows that
  \beq\label{e:var-F}
    |\beta| \Var(F) \le \frac{6|\beta|}{NL} \le c_1q^2,
  \eeq
  and Gaussian concentration implies that for some constant $c_3$ and all $0\le t\le 1$,
  \beq\label{e:F-conc}
    \bbP\lt(
      \beta \lt(F^2 - (\bbE F)^2\rt) \le -t
    \rt)
    \le e^{-c_3 NL t^2}.
  \eeq
  Combining \eqref{e:inward-deviation-bbp-mean} and \eqref{e:var-F}, we have
  \balnn
    \lt\{
      \fr{\beta\|\bB^\top\bw\|^2}{N\|\bw\|^2}
      \le \mu_* + c_1q^2 - u
    \rt\}
    &\subseteq \lt\{
      \fr{\beta\|\bB^\top\bw\|^2}{N\|\bw\|^2}
      - \fr{\beta\bbE \|\bB^\top\bw\|^2}{N\|\bw\|^2}
      \le -2c_1q^2 - u
    \rt\} \\
    &\subseteq \lt\{
      \fr{\beta L^2}{N\|\bw\|^2}(d_1^2-\bbE d_1^2)
      + \beta \lt(F^2 - (\bbE F)^2\rt)
      \le -c_1q^2 - u
    \rt\}\,.
  \ealnn
  If we take $t = (c_1 q^2 + u) / 2$ in \eqref{e:d1-conc}, \eqref{e:F-conc}, then we conclude
  \balnn
    \bbP\lt(
      \fr{\beta\|\bB^\top\bw\|^2}{N\|\bw\|^2}
      \le \mu_* + c_1q^2 - u
    \rt)
    &\le \bbP\lt(\fr{\beta L^2}{N\|\bw\|^2}(d_1^2-\bbE d_1^2) \le -t\rt)
    + \bbP\lt(\beta \lt(F^2 - (\bbE F)^2\rt) \le -t\rt) \\
    &\le 2\exp\lt(-cN\fr{(u+q^2)^2}{q}\rt),
  \ealnn
  for suitable $c$ depending on $c_1, c_2, c_3, C$.
  This proves \eqref{e:wishart-inward-deviation-bbp-goal}.
\end{proof}

%\newpage

\section{Stability in the plant and Hopfield overlap concentration}
\label{s:cov}

This section is devoted to proving Lemmas~\ref{l:good-minor} and~\ref{l:hopfield-overlap-concentration}. We begin with a unified bound on the Gibbs covariance matrix in terms of the $(\zeta I_N-\beta \bG^\top\bG/N)^{-1}$. The two lemmas are then proved building on this bound, with different treatments of the high-temperature and critical regimes. Throughout this section, we use ``with overwhelming probability'' to indicate with probability at least $1-N^{-D}$ for any fixed constant $D$ when $N$ is large enough. 

\subsection{Covariance bound}

For simplicity, we write $\bW_\bG=\bG^T\bG/N\in \bbR^{N\times N}$. Denote the spherical Hopfield Gibbs measure with respect to disorder $\bG$ as $\langle \cdot \rangle_\bG$. In particular, $\langle \bx\bx^\top\rangle_\bG\in \bbR^{N\times N}$ is the Gibbs covariance. Note that this implicitly depends on the temperature parameter $\beta$, but it should be clear from the context so we suppress it from the notation. 

\begin{lem}%[Hopfield Gibbs covariance]
\label{l:hopfield-gibbs-covariance}
There exists a universal constant $C$ such that at any temperature,
\[
\left\langle \bx\bx^\top\right\rangle_{\bG}
\preceq
C\lt(
\zeta I_N-\beta\bW_{\bG}
\rt)^{-1}.
\]
\end{lem}

\begin{proof}
By rotational invariance, we may work in an eigenbasis of
$\bW_{\bG}$. In this basis, the Hopfield Gibbs measure is the law of
independent Gaussian variables
\[
g_i\sim\cN(0,a_i),
\qquad
a_i=\frac{1}{\zeta-\beta\lambda_i},
\]
conditioned on (understood by coarea)
$
\sum_{i=1}^N g_i^2=N=\sum_{i=1}^N a_i.
$
The conditional law is invariant under changing the sign of any
coordinate, so
$
\langle g_i g_j\rangle_{\bG}=0$ for all $i\ne j.
$

It remains to control the diagonal entries. We claim that
$
\langle g_j^2\rangle_{\bG}\le Ca_j.
$
Weighting by $g_j^2/a_j$ changes the distribution of the $j$-th term from $a_j\chi_1^2$ to $a_j\chi_3^2$. Thus, the conditional law of $S=\sum_i a_iZ_i^2$ is changed to the law of $S^{(j)}=a_j\tilde{Z}_j^2+\sum_{i\neq j}a_iZ_i^2$ under the same conditioning (recall \eqref{e:V-S-Sj} for the relevant definitions). This yields that
\[
\langle g_j^2\rangle_{\bG}
=
a_j\,
\frac{f_{S^{(j)}}(N)}{f_{S_\zeta}(N)}.
\]
The bounds in
\eqref{e:weighted-chi-square-density} and
\eqref{e:weighted-chi-square-augmented-density} therefore imply the claim. 
Hence, in the eigenbasis of $\bW_{\bG}$,
\[
\left\langle \bx\bx^\top\right\rangle_{\bG}
\preceq
C\,\mathrm{diag}(a_1,\ldots,a_N)
=
C\lt(
\zeta I_N-\beta\bW_{\bG}
\rt)^{-1},
\]
which proves the result.
\end{proof}

\subsection{Stability in the plant}\label{ss:good-minor}

We now prove Lemma~\ref{l:good-minor}. Following the notation of Lemma~\ref{l:good-minor}, we fix an isometry $\bE\in \mathbb{R}^{M\times 2}$ and write
$
\bG_{\bK}=\bG_\bK(\bA_0)=[\bE\bK^{1/2},\bA_0].
$
The case $s=1$ is immediate, since then $\beta=0$ and
$f_{\bK}\equiv1$. We now assume $s\ne1$. Our plan is to explicitly compute the derivative $\nabla_\bK \log f_\bK(\bA_0)$ and then use Lemma~\ref{l:hopfield-gibbs-covariance} to control it by $\bE^\top(\zeta I_M-\beta\bG_\bK\bG_\bK^\top/N)^{-1}\bE$.

To express the derivative, we introduce some notations. Fix $\bA_0\in \bbR^{M\times (N-2)}$ and let $\mathcal S=\operatorname{Range}(\bA_0)+\operatorname{Range}(\bE)$ and $r=\operatorname{rank}(\mathcal S)$.
For positive definite $\bK$, we have $\operatorname{Range}(\bG_\bK)=\mathcal S$. Thus, we may take a thin singular value decomposition $\bG_\bK=\bU\boldsymbol{\Sigma} \bV^\top$, 
$\bU\in\mathbb R^{M\times r}$ and $\bV\in\mathbb R^{N\times r}$ have
orthonormal columns, and $\boldsymbol{\Sigma}\in\mathbb R^{r\times r}$ is positive diagonal.
Define
\[
 \bT_\bK:=\bU(\bV^\top \langle\bx\bx^\top\rangle_{\bG_\bK} \bV)\bU^\top\in\mathbb R^{M\times M}.
\]
It is supported on $\mathcal S$ and is independent of the choice of
singular vectors.\footnote{Since $\operatorname{Range}(\bU)=\mathcal S$,
$\bT_\bK$ maps into $\mathcal S$ and vanishes on $\mathcal S^\perp$.
Changing singular vectors replaces $(\bU,\bV)$ by $(\bU\bO,\bV\bO)$ with $\bO$
orthogonal, leaving $\bU(\bV^\top \langle\bx\bx^\top\rangle_{\bG_\bK} \bV)\bU^\top$ unchanged.}
Moreover, define 
\[
 \bB=\bA_0\bA_0^\top/N\in \bbR^{M\times M},\quad
 \bD_\bK=\bG_\bK\bG_\bK^\top/N
 =\bB+\bE\bK\bE^\top/N\in \bbR^{M\times M},
\]
Note that restricting to $\mathcal S$ removes only directions on which both
$\bB$ and $\bD_\bK$ vanish, and thus preserves all their nonzero spectral information. 

\begin{lem}\label{l:gram-response}
For every positive definite $\bK\in \bbR^{2\times 2}$,
\begin{equation}\label{e:high-temp-gram-response}
 \nabla_\bK\log f_\bK(\bA_0)
 =-\frac{\beta}{N}\bE^\top
 \bT_\bK\bE,
\end{equation}
Moreover, as operators on $\mathcal S$,
\begin{equation}\label{e:row-space-covariance}
 \bT_\bK|_{\mathcal S}
 \preceq C(\zeta_\bK I_{\mathcal S}
                   -\beta\bD_\bK|_{\mathcal S})^{-1}.
\end{equation}
Consequently,
\begin{equation}\label{e:gram-response-resolvent}
 \|\nabla_\bK\log f_\bK(\bA_0)\|_F
 \le\frac{C|\beta|}{N}
 \left\|\bE^\top(\zeta_\bK I_{\mathcal S}
                  -\beta\bD_\bK|_{\mathcal S})^{-1}\bE\right\|_{\mathrm{op}}.
\end{equation}
\footnote{Whenever the inverse on the full space exists,
$\bE^\top(\zeta_\bK I_\mathcal{S}-\beta\bD_\bK|_\mathcal{S})^{-1}\bE=\bE^\top(\zeta_\bK I_M-\beta\bD_\bK)^{-1}\bE$, since
$\operatorname{Range}(\bE)\subseteq\mathcal S$ and $\bD_\bK$ vanishes on
$\mathcal S^\perp$. We retain the restricted form for the operator norm
estimate below, avoiding the additional factor $|\zeta_\bK|^{-1}$
from the zero directions in the full-space resolvent.}Here, $\bE^\top(\zeta_\bK I_{\mathcal S}
                  -\beta\bD_\bK|_{\mathcal S})^{-1}\bE\in \bbR^{2\times 2}$ is understood as the composition of the three operators $$\bE:\bbR^2\to \mathcal S, \qquad(\zeta_\bK I_{\mathcal S}
                  -\beta\bD_\bK|_{\mathcal S})^{-1}:\mathcal S\to \mathcal S, \qquad \bE^\top:\mathcal S\to \bbR^2.$$
\end{lem}

\begin{proof}
We first establish \eqref{e:row-space-covariance} as a consequence of Lemma~\ref{l:hopfield-gibbs-covariance}. 
In a complete right singular-vector basis, $\bW_{\bG_\bK}=\bG_\bK^\top\bG_\bK^{\ }/N$ is diagonal. Thus, the Gibbs 
covariance restricting on this basis is also diagonal as argued in Lemma~\ref{l:hopfield-gibbs-covariance}.
Consequently, $\bA:=\bV^\top\langle\bx\bx^\top\rangle_{\bG_\bK}\bV$
commutes with $\boldsymbol{\Sigma}$.
Thus
\begin{equation}\label{e:commute}
 \bT_\bK\bG_\bK=\bU\bA\boldsymbol{\Sigma} \bV^\top
 =\bU\boldsymbol{\Sigma} \bA\bV^\top
 =\bG_\bK\langle\bx\bx^\top\rangle_{\bG_\bK}.
\end{equation}
Since $\bW_{\bG_\bK}=\bV\boldsymbol{\Sigma}^2\bV^\top/N$, Lemma~\ref{l:hopfield-gibbs-covariance} gives
\[
 \bV^\top\langle\bx\bx^\top\rangle_{\bG_\bK}\bV
 \preceq C \bV^\top(\zeta_\bK I_N-\beta\bW_{\bG_\bK})^{-1}\bV
 =C(\zeta_\bK I_r-\beta\boldsymbol{\Sigma}^2/N)^{-1}.
\]
Multiplying by $\bU$ and $\bU^\top$ and using
$\bD_\bK=\bU\boldsymbol{\Sigma}^2\bU^\top/N$, we obtain, as operators on $\mathcal S$,
\[
 \bT_\bK|_{\mathcal S}\preceq C \bU(\zeta_\bK I_r-\beta\boldsymbol{\Sigma}^2/N)^{-1}\bU^\top
 =C(\zeta_\bK I_{\mathcal S}-\beta\bD_\bK|_{\mathcal S})^{-1},
\]
which proves \eqref{e:row-space-covariance}.

We now compute and bound the derivative. Write $\bT=\bT_\bK$ and $\bR=\bK^{1/2}$ for simplicity. 
Recall that
\[
 f_\bK(\bA_0)
 =Z^{\mathrm{sph}}_{\mathrm{hop},N}([\bE\bR,\bA_0])^{-2}
 =\left(\int_{S_N}
   \exp\left(\frac{\beta}{2N}\|\bG_\bK\bx\|^2\right)
   \,\mathrm d\nu_N(\bx)\right)^{-2}.
\]
Fix a symmetric $\bH\in \bbR^{2\times 2}$. Differentiating under the integral gives
\[
 \partial_\bH\log f_\bK(\bA_0)=\frac{\de} {\de t}\log f_{\bK+t\bH}(\bA_0)\Big|_{t=0}
 =-2\left\langle
 \partial_\bH\left(\frac{\beta}{2N}\|\bG_\bK\bx\|^2\right)
 \right\rangle_{\bG_\bK}.
\]

To proceed with the calculation, write $\bx=(\bu,\bv)$ with $\bu\in\mathbb R^2,\bv\in \bbR^{N-2}$,
then $\bG_\bK\bx=\bE\bR\bu+\bA_0\bv$. Denote $\bX=\partial_\bH\bR=\frac{\de (\bK+t\bH)^{1/2}}{\de t}\big|_{t=0}\in \bbR^{2\times 2}$. 
 We have
\[
 \partial_\bH(\bG_\bK\bx)=\bE\bX\bu,
 \qquad
 \partial_\bH\left(\frac{\beta}{2N}\|\bG_\bK\bx\|^2\right)
 =\frac{\beta}{N}(\bG_\bK\bx)^\top\bE\bX\bu=\frac{\beta}{N}\operatorname{Tr}\big(u(\bG_\bK\bx)^\top\bE\bX\big).
\]
Consequently,
\[
 \partial_\bH\log f_\bK(\bA_0)
 =-\frac{2\beta}{N}\left\langle\operatorname{Tr}\big(u(\bG_\bK\bx)^\top\bE\bX\big)\right\rangle_{\bG_\bK}=-\frac{2\beta}{N}\operatorname{Tr}\left(\left\langle u(\bG_\bK\bx)^\top\right\rangle_{\bG_\bK}\bE\bX\right).
\]
Taking transpose of 
$\bG_\bK\langle\bx\bx^\top\rangle_{\bG_\bK}=\bT\bG_\bK$ in \eqref{e:commute} yields $\langle \bx\bx^\top\rangle \bG_\bK^\top=\bG_\bK^\top \bT$. The first two rows give
$$\langle\bu\bx^\top\rangle_{\bG_\bK}\bG_\bK^\top
=\langle\bu(\bG_\bK^\top x)^\top\rangle=\bR\bE^\top\bT.$$ 
This implies
\[
 \begin{aligned}
 \partial_\bH\log f_\bK(\bA_0)
 &=-\frac{2\beta}{N}\operatorname{Tr}(\bR\bE^\top \bT\bE\bX).
 \end{aligned}
\]
Since $\bX,\bR,\bE^\top\bT\bE$ are all symmetric, we have
\[
\operatorname{Tr}(\bR\bE^\top \bT\bE\bX)=\operatorname{Tr}(\bX\bR\bE^\top \bT\bE)=\operatorname{Tr}((\bX\bR\bE^\top \bT\bE)^\top)=\operatorname{Tr}(\bE^\top \bT\bE\bR\bX)=\operatorname{Tr}(\bR\bX\bE^\top \bT\bE).
\]
Moreover, differentiating $(\bK+t\bH)^{1/2}(\bK+t\bH)^{1/2}=\bK+t\bH$ gives
\[
 \bH=\left.\frac{\mathrm d}{\mathrm dt}
 \left((\bK+t\bH)^{1/2}(\bK+t\bH)^{1/2}\right)\right|_{t=0}
 =\bX\bR+\bR\bX.
\]
Consequently, we have
\[
 \begin{aligned}
 \partial_\bH\log f_\bK(\bA_0)
 &=-\frac{\beta}{N}\operatorname{Tr}((\bX\bR+\bR\bX)\bE^\top \bT\bE)
 =-\frac{\beta}{N}\operatorname{Tr}(\bH\bE^\top \bT\bE).
 \end{aligned}
\]
This proves \eqref{e:high-temp-gram-response}; combining it with
\eqref{e:row-space-covariance} and
$\|\bA\|_F\le\sqrt2\|\bA\|_{\mathrm{op}}$ for $2\times2$ matrices
proves \eqref{e:gram-response-resolvent}.
\end{proof}

For $\bK\in\mathcal K_N$, define
\[
\mu_{\bK}
:=
\lambda_{\max}(\beta\bW_{\bG_{\bK}}),
\]
and, for $z>\mu_{\bK}$ and $z\neq 0$, set
\[
\boldsymbol{\mathcal R}_{\bK}(z)
:=
\bigl(zI_{\mathcal S}-\beta\bD_{\bK}|_{\mathcal S}\bigr)^{-1}.
\]
In light of Lemma~\ref{l:gram-response}, Lemma~\ref{l:good-minor} follows from the following Proposition. 
\begin{ppn}\label{p:projected-resolvent-bound}
For $\bA_0\in \bbR^{M\times(N-2)}$ with i.i.d. standard Gaussian entries, with overwhelming probability
\begin{equation}
\label{e:projected-resolvent-bound}
\sup_{\bK\in\mathcal K_N}
\left\|
\bE^\top\boldsymbol{\mathcal R}_{\bK}(\zeta_{\bK})\bE
\right\|_{\mathrm{op}}
\le CN^{1/3+\varepsilon_0}.
\end{equation}
\end{ppn}

Since the high-temperature and critical regimes have quantitatively
different behavior, we prove Proposition~\ref{p:projected-resolvent-bound}
separately in these two cases.

\begin{proof}[Proof of Proposition~\ref{p:projected-resolvent-bound}, the high-temperature case]
Suppose that $\theta<1$. Since $\bE:\mathbb R^2\to\mathcal S$
is an isometry, we have
\begin{equation}\label{e:bound-operator-norm}
\sup_{\bK\in\mathcal K_N}
\left\|
\bE^\top\boldsymbol{\mathcal R}_{\bK}(\zeta_{\bK})\bE
\right\|_{\mathrm{op}}
\le
\sup_{\bK\in\mathcal K_N}
\left\|
\boldsymbol{\mathcal R}_{\bK}(\zeta_{\bK})
\right\|_{\mathrm{op}}.
\end{equation}
We will show that the latter is $O(1)$ with overwhelming probability,
which gives a stronger form of
\eqref{e:projected-resolvent-bound}.
To this end, we will establish that with overwhelming probability, $\zeta_\bK-\mu_\bK\ge \Omega(1)$ uniformly over $K\in \mathcal K_N$. Let $\bK_0=MsI_2$. We first show $\zeta_{\bK_0}-\mu_{\bK_0}\ge \Omega(1)$ with overwhelming probability, and then
extend it to all $\bK\in\mathcal K_N$.

By Fact~\ref{f:zeta-above-spectrum} we have $\zeta_\star>\mu_\star$. Fix $0<\delta<(\zeta_\star-\mu_\star)/4$. On the one hand, 
since $\alpha<(s-1)^{-2}$, the finite-rank perturbation is in the BBP subcritical phase.
Thus, when $\beta>0$, \cite[Theorem~2.3 and Remark~2.4]{bloemendal2016principal} gives
$\mu_{\bK_0}\le\mu_\star+\delta$ with overwhelming probability.\footnote{The cited result applies to the Gaussian spiked model
$[\sqrt{s}\,\bZ,\bA_0]$, which corresponds to the spectrum of $\bD_{s\bZ^\top\bZ}$.
Gaussian concentration gives
$N^{-1}\|s\bZ^\top\bZ-MsI_2\|_{\mathrm{op}}=o(1)$ with overwhelming
probability, so the same estimate transfers to $\bK_0$.}
When $\beta<0$ and $\alpha\le1$, the bound follows directly from
$\mu_{\bK_0}\le0=\mu_\star$.
On the other hand, interlacing and concentration of resolvent trace to its deterministic analogue, as in the proof of
Lemma~\ref{l:zeta-deviation}, give that with overwhelming probability
\[
 \frac{1}{N}\operatorname{Tr}\left((\zeta_{\bK_0}-\delta)I_N-\beta \bW_{\bK_0}^{-1}\right)\ge\frac1N\operatorname{Tr}
 \left((\zeta_{\bK_0}-\delta)I_{N-2}-\frac{\beta}{N}\bA_0^\top\bA_0\right)^{-1}
 =\int \frac{\nu_{\MP,\alpha}(\de x)}{\zeta_\star-\delta-x}+o(1)>1.
\]
Here the last inequality follows from the fact that $\int \frac{\nu_{\MP,\alpha}(\de x)}{\zeta_\star-x}=1$. 
Monotonicity then yields
$\zeta_{\bK_0}>\zeta_\star-\delta$. with overwhelming probability. 
Altogether we conclude that with overwhelming probability
\begin{equation}\label{e:saddle-separation-K0}
 \zeta_{\bK_0}-\mu_{\bK_0}\ge 2\delta. 
\end{equation}

We next upgrade \eqref{e:saddle-separation-K0} to all
$\bK\in\mathcal K_N$. By the definition of $\mathcal K_N$,
\[
r_N:=
\frac{|\beta|}{N}
\sup_{\bK\in\mathcal K_N}\|\bK-\bK_0\|_{\mathrm{op}}
=O(N^{-1/3}\log N)=o(1).
\]
Since
\[
 \beta\bD_\bK-\beta\bD_{\bK_0}
 =\frac{\beta}{N}\bE(\bK-\bK_0)\bE^\top,
 \qquad
 \sup_{\bK\in\mathcal K_N}
 \|\beta\bD_\bK-\beta\bD_{\bK_0}\|_{\mathrm{op}}\le r_N,
\]
Let $\mu_1(\bK)\ge\cdots\ge\mu_N(\bK)$ be the eigenvalues of
$\beta\bW_{\bG_\bK}$. Weyl's inequality gives
\[
\mu_i(\bK)\le\mu_i(\bK_0)+r_N,
 \qquad 1\le i\le N,\qquad \text{for all }\bK\in\mathcal K_N.
\]
Consequently, provided that $\zeta_{\bK_0}-r_N>\mu_\bK$, 
\[
 \begin{gathered}
 \frac{1}{N}\sum_{i=1}^N\frac{1}{\zeta_{\bK_0}-r_N-\mu_i(\bK)}
 \ge \frac1N\sum_{i=1}^N\frac1{\zeta_{\bK_0}-\mu_i(\bK_0)}=1.
 \end{gathered}
\]
Monotonicity then gives $\zeta_\bK\ge \zeta_{\bK_0}-r_N$. 
Hence, on the same overwhelming-probability event,
\[
\inf_{\bK\in\mathcal K_N}
\{\zeta_{\bK}-\mu_{\bK}\}
\ge 2\delta-2r_N\ge \delta.
\]
Finally, since
$\lambda_{\max}(\beta\bD_{\bK}|_{\mathcal S})\le\mu_{\bK}$,
we obtain
\[
\sup_{\bK\in\mathcal K_N}
\|\boldsymbol{\mathcal R}_{\bK}(\zeta_{\bK})\|_{\mathrm{op}}
=
\sup_{\bK\in\mathcal K_N}
\frac{1}{
\zeta_{\bK}-\lambda_{\max}(\beta\bD_{\bK}|_{\mathcal S})
}
\le
\sup_{\bK\in\mathcal K_N}
\frac{1}{\zeta_{\bK}-\mu_{\bK}}
\le \delta^{-1},
\]
which proves the desired bound.
\end{proof}

We next handle the critical case. Assume henceforth that $\theta=1$. Since $\zeta$ lies on the soft-edge
scale, the operator norm estimate \eqref{e:bound-operator-norm}
gives only $O(N^{2/3+o(1)})$, which is insufficient for
\eqref{e:projected-resolvent-bound}. 
However, we expect a better bound for $E^\top\mathcal R_K(\zeta_K)E$, since the eigenvectors corresponding to the largest eigenvalues of $\mathcal R_K(\zeta_K)$ should have small projections onto the fixed directions in $E$.

To make this precise, we write $\bD_\bK=\bB+\bE\bK\bE^\top/N$ and
apply the Woodbury identity to express $\bE^\top\mathcal R_\bK(\zeta_\bK)\bE$ in
terms of $\bB$ (see \eqref{eq:critical-woodbury} below). Conditional on its spectrum, $\bB$ has a Haar-distributed
eigenbasis, allowing us to use Haar concentration. One complication is that 
$\zeta_\bK$ remains random under this conditioning. We handle this by
introducing an interval $\mathcal I=\mathcal I_\bB$, determined only by
the spectrum of $\bB$, that contains all $\zeta_\bK$ with
overwhelming probability, and prove the required estimates uniformly
on this interval.

We now define this interval. In what follows, we let $\varepsilon_0$ be as in Lemma~\ref{l:good-minor} and denote $\varepsilon=\varepsilon_0/10$. Recall from Fact~\ref{f:zeta-above-spectrum} that
$\zeta_\star=\mu_\star\ne0$. Write as before $\mu_\bK=\lambda_{\max}(\beta \bW_{\bK})=\lambda_{\max}(\beta \bD_\bK\mid_\mathcal{S})$ and
\[
\nu_\bB:=
\begin{cases}
\beta\lambda_{\max}(\bB), & \beta>0,\\[2mm]
\beta\lambda_{\min}^{+}(\bB), & \beta<0,
\end{cases}
\]
where $\lambda_{\min}^{+}(\bB)$ denotes the smallest positive
eigenvalue of $\bB$.
Note that by interlacing, $\mu_\bK\ge \nu_\bB$ for all $\bK\in \mathcal K_N$.
Define
\begin{equation}\label{eq:critical-admissible-region}
\mathcal I=\mathcal I_{\bB}
=
\left[
\nu_\bB+N^{-2/3-\varepsilon},\nu_\bB+N^{-2/3+\varepsilon}
\right].
\end{equation}
This interval depends only on the eigenvalues of $\bB$. Since $\nu_\bB=\mu_\star+o(1)$ with overwhelming probability, on the same event there exists a constant $c>0$ such that every $z\in\mathcal I$ satisfies
$|z|\ge c$ and $\beta z>0$.

We next apply the Woodbury identity and introduce some simplifying notations. In the following calculations, $\bB$ and $\bD_{\bK}$ act
on the full space $\mathbb R^M$. Write $\bR_0(z)=(zI_M-\beta\bB)^{-1}$. 
Whenever the inverses exist, the Woodbury identity gives
\[
 \begin{aligned}
 (zI_M-\beta\bD_\bK)^{-1}
 &=\left(\bR_0(z)^{-1}-\frac{\beta}{N}\bE\bK\bE^\top\right)^{-1}
 =\bR_0(z)+\bR_0(z)\bE
 \left(\frac{N}{\beta}\bK^{-1}-\bE^\top\bR_0(z)\bE\right)^{-1}
 \bE^\top\bR_0(z).
 \end{aligned}
\]
Denote
\[
\bT_0(z)=\bE^\top\bR_0(z)\bE,
\qquad
\bT_2(z)=\bE^\top\bR_0(z)^2\bE,
\qquad
\bS_{\bK}(z)=\frac{N}{\beta}\bK^{-1}-\bE^\top\bR_0(z)\bE.
\]
Multiplying the Woodbury equation on the left by $\bE^\top$ and on the right by $\bE$ yields
\begin{equation}\label{eq:critical-woodbury}
\bE^\top\boldsymbol{\mathcal R}_\bK(z)\bE=\bE^\top(zI_M-\beta\bD_{\bK})^{-1}\bE
=
\bT_0(z)+\bT_0(z)\bS_{\bK}(z)^{-1}\bT_0(z),
\end{equation}
where the first identity holds since $\bD_\bK=0$ on $\mathcal S^\perp$ and the columns of $\bE$ belong to $\mathcal S$.
Also,
\begin{equation}\label{eq:critical-S-derivative}
\bS_{\bK}'(z)=\bT_2(z).
\end{equation}
Our goal is to show that with overwhelming probability, for every $\bK\in \mathcal K_N$ and $z=\zeta_\bK$, $\|\bT_0(z)\|_{\mathrm{op}}=O(1)$ and $\|\bS_\bK(z)^{-1}\|_{\mathrm{op}}=O(N^{1/3+o(1)})$. Given this is true, \eqref{e:projected-resolvent-bound} follows immediately from \eqref{eq:critical-woodbury}. 

Write $\bB=\bU\bLambda\bU^\top$, where $\bLambda$ is the diagonal
matrix of eigenvalues and $\bU$ is Haar orthogonal, independent of
$\bLambda$. As outlined above, we will combine uniform localization of $\zeta_\bK$ in $\mathcal I$ with conditional Haar estimates for
$\bT_0(z)$ and $\bT_2(z)$, uniform over $z\in\mathcal I$.
The following two lemmas provide these inputs; their proofs are postponed.

\begin{lem}
\label{l:critical-saddle}
For every $D>0$ and all sufficiently large $N$, with probability
at least $1-N^{-D}$, simultaneously for every $\bK\in\mathcal K_N$,
\begin{equation}\label{eq:critical-saddle-location}
[\zeta_\bK-N^{-2/3-\varepsilon},\zeta_{\bK}]\subset \mathcal I,
\qquad
\zeta_{\bK}-\mu_{\bK}\ge 2N^{-2/3-\varepsilon}.
\end{equation}
\end{lem}

\begin{lem}
\label{l:critical-haar}
There exist constants $c,C>0$ and an event $\mathcal E_\Lambda$,
measurable with respect to the eigenvalues of $\bB$,
with overwhelming probability, such that the following holds.
Conditionally on any $\bLambda\in\mathcal E_\Lambda$, with overwhelming probability, simultaneously for every $z\in \mathcal I$,
\begin{equation}\label{eq:critical-haar-estimates}
\|\bT_0(z)\|_{\mathrm{op}}\le C,
\qquad
\bT_2(z)\succeq c N^{1/3-\varepsilon}I_2. 
\end{equation}
\end{lem}

We now complete the proof of Proposition~\ref{p:projected-resolvent-bound} at criticality, which concludes Lemma~\ref{l:good-minor}. 

\begin{proof}[Proof of Proposition~\ref{p:projected-resolvent-bound}, the critical case]
Working on the intersection of the events in Lemmas~\ref{l:critical-saddle} and~\ref{l:critical-haar} which holds with overwhelming probability, we may assume both \eqref{eq:critical-saddle-location} and \eqref{eq:critical-haar-estimates} hold. We show that \eqref{e:projected-resolvent-bound} then follows deterministically. 

Fix $\bK\in\mathcal K_N$. For every
$u\in[\zeta_\bK-N^{-2/3-\varepsilon},\zeta_\bK]$, we claim that $\bS_\bK(u)\succ 0$. To see this,  Lemma~\ref{l:critical-saddle} gives $u\in\mathcal I$ and $u>\mu_{\bK}$, which implies $\beta u>0$ and
\[
uI_N-\frac{\beta}{N}\bG_{\bK}^\top\bG_{\bK}
=
\begin{pmatrix}
uI_2-\dfrac{\beta}{N}\bK
&
-\dfrac{\beta}{N}\bK^{1/2}\bE^\top\bA_0
\\[2mm]
-\dfrac{\beta}{N}\bA_0^\top\bE\bK^{1/2}
&
uI_{N-2}-\dfrac{\beta}{N}\bA_0^\top\bA_0
\end{pmatrix}
\succ0.
\]
Applying the Schur-complement criterion to the
lower-right $(N-2)\times(N-2)$ block, we obtain
\[
\begin{aligned}
0\prec{}&
uI_2-\frac{\beta}{N}\bK-\frac{\beta^2}{N^2}
\bK^{1/2}\bE^\top\bA_0
\left(uI_{N-2}-\frac{\beta}{N}\bA_0^\top\bA_0\right)^{-1}
\bA_0^\top\bE\bK^{1/2}\\
={}&
uI_2-\frac{\beta}{N}\bK
-\frac{\beta^2}{N}
\bK^{1/2}\bE^\top\bB\bR_0(u)\bE\bK^{1/2}=
\frac{\beta u}{N}
\bK^{1/2}\bS_{\bK}(u)\bK^{1/2}.
\end{aligned}
\]
The second equality uses the identity
\[
 (uI_M-\beta\bB)\bA_0
 =\bA_0\left(uI_{N-2}-\frac{\beta}{N}\bA_0^\top\bA_0\right),
\]
which implies
\[
 \frac1N\bA_0
 \left(uI_{N-2}-\frac{\beta}{N}\bA_0^\top\bA_0\right)^{-1}
 \bA_0^\top
 =\bR_0(u)\bB=\bB\bR_0(u).
\]
The last equality uses $\bE^\top\bE=I_2$,
$\bK^{1/2}\bK^{-1}\bK^{1/2}=I_2$, and
\[
 I_M+\beta\bB\bR_0(u)
 =(uI_M-\beta\bB)\bR_0(u)+\beta\bB\bR_0(u)
 =u\bR_0(u),
\]
together with the definition of $\bS_\bK(u)$.
Since $\beta u>0$ and $\bK\succ0$, the claim follows. 

The claim together with
\eqref{eq:critical-S-derivative} and
\eqref{eq:critical-haar-estimates} imply
\[
\begin{aligned}
\bS_{\bK}(\zeta_\bK)
&=
\bS_{\bK}(\zeta_\bK-N^{-2/3-\varepsilon})+\int_{\zeta_\bK-N^{-2/3-\varepsilon}}^{\zeta_\bK}\bT_2(u)\,\mathrm du
\succeq c N^{-1/3-2\varepsilon}I_2.
\end{aligned}
\]
Consequently,
$
\|\bS_{\bK}(\zeta_{\bK})^{-1}\|_{\mathrm{op}}
\le CN^{1/3+2\varepsilon}.
$
Combining this with \eqref{eq:critical-woodbury} and \eqref{eq:critical-haar-estimates} gives
\[
\left\|
\bE^\top(\zeta_{\bK}I_M-\beta\bD_{\bK})^{-1}\bE
\right\|_{\mathrm{op}}
\le \|\bT_0(z)\|_{\mathrm{op}}+\|\bT_0(z)\|_{\mathrm{op}}^2\|\bS_\bK(z)^{-1}\|_{\mathrm{op}}\le CN^{1/3+2\varepsilon}\le CN^{1/3+\varepsilon_0}.
\]
This holds simultaneously for all $\bK\in \mathcal K_N$ and thus the result follows. 
\end{proof}

\begin{proof}[Proof of Lemma~\ref{l:critical-saddle}]
Let $\mu_1(\bK)\ge\cdots\ge\mu_N(\bK)$ be the eigenvalues of
$\beta\bW_{\bG_\bK}=\beta \bG_\bK^\top\bG_\bK^{\ }/N$, and $\nu_1\ge\cdots\ge\nu_{N-2}$ be the nonzero eigenvalues of
$\beta\bB$, counted with multiplicity and supplemented with zeros
to form a list of length $N-2$. Thus
$\mu_1(\bK)=\mu_{\bK}$ and $\nu_1=\nu_{\bB}$.
Since $\nu_j,1\le j\le N-2$ are also the eigenvalues of the lower-right
principal submatrix $\beta\bA_0^\top\bA_0/N$, interlacing gives
\begin{equation}\label{eq:critical-interlacing}
\mu_j(\bK)\ge\nu_j\ge\mu_{j+2}(\bK),
\qquad 1\le j\le N-2,
\end{equation}
simultaneously for every $\bK\in\mathcal K_N$. In what follows, we denote for brevity
\[
h_N=N^{-2/3+\varepsilon},\quad v=\mu_\star+h_N/2,\quad w=\mu_\star+h_N/4.
\]

We first provide some random matrix estimate. To begin, we claim that with overwhelming probability, $\mu_\bK<w$ for all $K\in \mathcal K_N$. To this end, we reduce to the standard population-spiked Wishart model studied in \cite{bloemendal2016principal}. 
Let $L>0$ be a sufficiently large constant and set
\[
s_N=s+\operatorname{sgn}(\beta)LN^{-1/3}\log N.
\]
Let $\bZ\in\mathbb R^{M\times 2}$ be a standard Gaussian matrix,
independent of $\bA_0$, and put
$
\widehat{\bK}=s_N\bZ^\top\bZ.
$
As in Lemma~\ref{l:good-gram}, Gaussian concentration gives
\[
\|\widehat{\bK}-Ms_NI_2\|_{\mathrm{op}}
\le C N^{1/2}\log N
\]
with overwhelming probability. On the other hand, by the definition
of $\mathcal K_N$,
\[
\sup_{\bK\in\mathcal K_N}
\|\bK-MsI_2\|_{\mathrm{op}}
\le C N^{2/3}\log N.
\]
Hence, taking $L$ sufficiently large, with overwhelming probability,
\[
\operatorname{sgn}(\beta)(\widehat{\bK}-\bK)
\succeq \left(M|s_N-s|-\|\widehat{\bK}-Ms_NI_2\|_{\mathrm{op}}-\|\bK-MsI_2\|_{\mathrm{op}}\right)I_2\succ 0,
\]
and thus $\mathrm{sgn}(\beta)\widehat{\bK}\succ \mathrm{sgn}(\beta)\bK$ simultaneously for all $\bK\in \mathcal K_N$. Consequently, on the same event, 
\begin{equation}\label{e:mu_K-reduction}
 \beta\bD_{\widehat{\bK}}-\beta\bD_\bK
 =\frac{\beta}{N}\bE(\widehat{\bK}-\bK)\bE^\top\succeq0,
 \qquad
 \mu_\bK\le\mu_{\widehat{\bK}},
 \qquad \text{for all }\bK\in\mathcal K_N.
\end{equation}
It therefore suffices to bound $\mu_{\widehat{\bK}}$.

By Gaussian polar decomposition and rotational invariance,
$\mu_{\widehat{\bK}}$ has the same law as the corresponding edge of
$\beta[\sqrt{s_N}\bZ,\bA_0]^\top[\sqrt{s_N}\bZ,\bA_0]/N$.
Its two population spikes satisfy
$s_N=s+O(N^{-1/3}\log N)$, where
$s=1+\operatorname{sgn}(\beta)\alpha^{-1/2}$ is the BBP threshold.
Thus \cite[Theorem~2.3 and Remark~2.4]{bloemendal2016principal}
gives, with overwhelming probability,
\[
\mu_{\widehat{\bK}}\le\mu_*+O(N^{-2/3+\varepsilon/2})<w.
\]
Combining with \eqref{e:mu_K-reduction}, this verifies the claim that $\mu_\bK<w$ simultaneously for all $\bK\in \mathcal K_N$ with overwhelming probability.

We next record a resolvent trace estimate. Recall from Fact~\ref{f:zeta-above-spectrum} that $\zeta_\star=\mu_\star$ and $\gamma_N=M/(N-2)$. The exterior local-law \cite[Theorems~2.5 and~2.10, Remark~2.7]{BEKYY14} applying to $\bA_0^\top\bA_0/N$ implies that on an overwhelming-probability event, for $t\in\{h_N/4,h_N/2\}$, 
\begin{align*}
\frac1N\sum_{j=1}^{N-2}
\frac1{\zeta_\star+t-\nu_j}
&=\int
\frac{\nu_{\mathrm{MP},\gamma_N}(\mathrm dx)}
{\mu_\star+t-\beta x}+O(N^{-1/2+\varepsilon/2}h_N^{-1/4})
\end{align*}
Since $\gamma_N=\alpha+O(N^{-1})$, and the derivative of the Marchenko--Pastur Stieltjes transform with respect to the aspect ratio is of order
$O(h_N^{-1/2})$ near a regular edge,
\begin{align*}
\int
\frac{\nu_{\mathrm{MP},\gamma_N}(\mathrm dx)}
{\mu_\star+t-\beta x}=\int
\frac{\nu_{\mathrm{MP},\alpha}(\mathrm dx)}
{\mu_\star+t-\beta x}+O(N^{-1}h_N^{-1/2})\,.
\end{align*}
Finally, the explicit Marchenko--Pastur transform gives, for $a=1/\sqrt s>0$,
\[
\int
\frac{\nu_{\mathrm{MP},\alpha}(\mathrm dx)}
{\mu_\star+t-\beta x}
=
1-a\sqrt t+O(t).
\]
Combining these estimates, we obtain with overwhelming probability that for $t\in \{h_N/4,h_N/2\}$, 
\begin{equation}\label{eq:critical-test-trace}
\frac1N\sum_{j=1}^{N-2}
\frac1{\zeta_\star+t-\nu_j}=1-a\sqrt t+o(\sqrt{h_N})\,.
\end{equation}
Here, we used the estimates 
\[
\max\{N^{-1/2+\varepsilon/2}h_N^{-1/4},N^{-1}h_N^{-1/2}\}=o(\sqrt{h_N})\iff \max\{N^{-\varepsilon/4},N^{-1/3-\varepsilon}\}=o(1).
\]
We then prove $\zeta_{\bK}<v$ with overwhelming probability. \eqref{eq:critical-test-trace} at $t=h_N/2$ yields that for some constant $c>0$ and all large enough $N$, 
\begin{equation}\label{e:saddle-monotonicity}
N-CN\sqrt{h_N}\le \sum_{j=1}^{N-2}\frac1{v-\nu_j}
\le
N-cN\sqrt{h_N}.
\end{equation}
By the above claim $\mu_\bK<w$ for all $\bK\in \mathcal K_N$ on an overwhelming-probability event. Assuming this and \eqref{e:saddle-monotonicity} are true, 
combining with \eqref{eq:critical-interlacing}, we obtain
\[
\begin{aligned}
\sum_{i=1}^N\frac1{v-\mu_i(\bK)}
&\le
\sum_{j=1}^{N-2}\frac1{v-\nu_j}
+\frac2{v-\mu_{\bK}}\le
N-cN\sqrt{h_N}+\frac8{h_N}
<N
\end{aligned}
\]
for sufficiently large $N$, since $Nh_N^{3/2}\to\infty$.
The definition \eqref{e:zeta} and monotonicity imply
$\zeta_{\bK}<v$, simultaneously for every $\bK\in\mathcal K_N$.

Finally, we prove $\zeta_\bK-\mu_\bK\ge 2N^{-2/3-\varepsilon}$ holds for all $\bK\in \mathcal K_N$ with overwhelming probability. Assuming that $\sup_{\bK\in \mathcal K_N}\zeta_\bK<v$ as established above. Equation \eqref{e:zeta}, \eqref{eq:critical-interlacing} and \eqref{e:saddle-monotonicity} give
\[
\begin{aligned}
\frac1{\zeta_{\bK}-\mu_{\bK}}
&=
N-\sum_{i=2}^N
\frac1{\zeta_{\bK}-\mu_i(\bK)}\le
N-\sum_{j=2}^{N-2}\frac1{v-\nu_j}
\\&=
N-\sum_{j=1}^{N-2}\frac1{v-\nu_j}
+\frac1{v-\nu_{\bB}}\le
CN\sqrt{h_N}+\frac{C}{h_N}
\le CN\sqrt{h_N}.
\end{aligned}
\]
Hence, uniformly in $\bK$,
\begin{equation}\label{e:saddle-separation-lower-bound}
\zeta_{\bK}-\mu_{\bK}
\ge
\frac{c}{N\sqrt{h_N}}
=
cN^{-2/3-\varepsilon/2}
\ge2N^{-2/3-\varepsilon}
\end{equation}
for all sufficiently large $N$.

Since $\mu_{\bK}\ge\nu_{\bB}$ by interlacing, \eqref{e:saddle-separation-lower-bound} implies
$
\zeta_{\bK}-N^{-2/3-\varepsilon}\ge \nu_{\bB}+N^{-2/3-\varepsilon}.
$
Moreover, rigidity gives that for sufficiently large $N$, with overwhelming probability,
$|\nu_{\bB}-\zeta_\star|\le h_N/4$, 
$
\zeta_{\bK}<v
=\zeta_\star+\frac{h_N}{2}
\le\nu_{\bB}+\frac{3h_N}{4}
<\nu_{\bB}+h_N.
$
Therefore, on the intersection of these overwhelming probability events, 
\[
[\zeta_{\bK}-N^{-2/3-\varepsilon},\zeta_{\bK}]
\subset[\nu_{\bB}+N^{-2/3-\varepsilon},\nu_{\bB}+N^{-2/3+\varepsilon}]
=\mathcal I
\]
simultaneously for every $\bK\in\mathcal K_N$.
This proves the lemma.
\end{proof}

\begin{proof}[Proof of Lemma~\ref{l:critical-haar}]
Let $\mathcal E_\Lambda$ be the rigidity event in \cite[Theorem~2.10]{BEKYY14} with slack $\varepsilon/20$, together with
$\|\bB\|_{\mathrm{op}}\le C$.
This event depends only on $\bLambda$ and has overwhelming
probability.

Let $r_\bB=\operatorname{rank}\bB$, so that
$\nu_1\ge\cdots\ge\nu_{r_\bB}$ are all the nonzero eigenvalues of
$\beta\bB$, and $r_\bB\asymp N$.
Set $J_0=\lceil N^\varepsilon\rceil$.
On $\mathcal E_\Lambda$, rigidity and the square-root edge density
give
\[
\nu_1-\nu_j
\ge cN^{-2/3}j^{2/3},
\qquad J_0\le j\le r_\bB.
\]
Since for $z\in \mathcal I$, $z-\nu_1\ge N^{-2/3-\varepsilon}$ and $|z|\ge c$,
\begin{equation}
\begin{aligned}
\sup_{z\in\mathcal I}\|\bR_0(z)\|_{\mathrm{op}}
&\le CN^{2/3+\varepsilon},\\
\sup_{z\in\mathcal I}\|\bR_0(z)\|_F^2
&\le
CJ_0N^{4/3+2\varepsilon}
+CN^{4/3}\sum_{j>J_0}j^{-4/3}+CN
\le CN^{4/3+3\varepsilon},
\label{eq:critical-minor-norms}
\end{aligned}
\end{equation}
where in the Frobenius norm estimate, zero eigenvalues contribute at most $CN$.
The same summation with first powers gives
\begin{equation}\label{eq:critical-minor-trace-bound}
\sup_{z\in\mathcal I}
\frac{|\operatorname{Tr}\bR_0(z)|}{M}
\le
C\left(
\frac{J_0}{N^{1/3-\varepsilon}}
+N^{-1/3}\sum_{j>J_0}j^{-2/3}+1
\right)
=O(1).
\end{equation}

To prove the desired bounds on $\bT_0$ and $\bT_2$, we use the following Haar concentration inequality from \cite[Theorem~2.3]{Ada15}. For every
deterministic symmetric $M\times M$ matrix $\bA$ and $t\ge1$,
\begin{equation}\label{eq:critical-haar-bound}
\mathbb P_{\bU}\left(
\left\|
\bE^\top\bU\bA\bU^\top\bE
-\frac{\operatorname{Tr}\bA}{M}I_2
\right\|_{\mathrm{op}}
>
C\left(
\frac{\sqrt t}{M}\|\bA\|_F
+\frac{t}{M}\|\bA\|_{\mathrm{op}}
\right)
\right)
\le Ce^{-ct}.
\end{equation}
Condition on any $\bLambda\in\mathcal E_\Lambda$.
For the bound on $\bT_0$, apply \eqref{eq:critical-haar-bound} with
$
\bA=(zI_M-\beta\bLambda)^{-1}$ and $t=(\log N)^2$
on an $N^{-4}$-net of $\mathcal I$. By
\eqref{eq:critical-minor-norms}, the resulting error is at most
\[
CN^{-1/3+3\varepsilon/2}(\log N)^2=o(1).
\]
The identity $\bR_0(z)-\bR_0(w)=(w-z)\bR_0(z)\bR_0(w)$ gives
\[
\|\bR_0(z)-\bR_0(w)\|_{\mathrm{op}}
\le |z-w|\,\|\bR_0(z)\|_{\mathrm{op}}\|\bR_0(w)\|_{\mathrm{op}}
\le CN^{4/3+2\varepsilon}|z-w|,
\qquad z,w\in\mathcal I,
\]
so the estimate extends from the net to the entire interval.
The conditional failure probability is at most
$CN^4e^{-c(\log N)^2}$.
Together with \eqref{eq:critical-minor-trace-bound}, this proves
the bound on $\bT_0$.

For the bound on $\bT_2$, let $\bP_{J_0}$ be the
spectral projection onto the first $J_0$ eigenvectors in the
ordering above. Rigidity gives, simultaneously for
$z\in\mathcal I$ and $j\le J_0$,
$
0<z-\nu_j\le CN^{-2/3+\varepsilon}.
$
Consequently,
\[
\bT_2(z)
\succeq cN^{4/3-2\varepsilon}\bE^\top\bP_{J_0}\bE,
\qquad z\in\mathcal I.
\]
Conditionally on $\bLambda$,
\[
\bP_{J_0}=\bU\begin{pmatrix}I_{J_0}&0\\0&0_{M-J_0}\end{pmatrix}\bU^\top,
\qquad \bU\sim\operatorname{Haar}(O(M)).
\]
Applying \eqref{eq:critical-haar-bound} to this projection with
$t=c_2J_0$, for sufficiently small $c_2>0$, gives
\[
\bE^\top\bP_{J_0}\bE\succeq\frac{J_0}{2M}I_2
\]
outside an event of probability at most $Ce^{-cJ_0}$.
It follows that, simultaneously for every $z\in\mathcal I$,
\[
\bT_2(z)
\succeq cN^{4/3-2\varepsilon}\frac{J_0}{M} I_2
\succeq cN^{1/3-\varepsilon}I_2.
\]
Since $J_0\asymp N^\varepsilon$, this conditional failure
probability is overwhelming. Both conditional probability bounds
are uniform over $\bLambda\in\mathcal E_\Lambda$, completing the
proof.
\end{proof}

\subsection{Hopfield Gibbs overlap}

We now prove Lemma~\ref{l:hopfield-overlap-concentration}. Let $\bG\in\mathbb R^{M\times N}$ have
i.i.d.\ standard Gaussian entries, and let $\zeta$ defined by \eqref{e:zeta}. Throughout this section we write
$\bW=\frac{\bG^\top\bG}{N},\mu_+=\lambda_{\max}(\beta\bW).$

We first record the gap between $\zeta$ and $\mu_+$ for the unplanted Gaussian
matrix. 
\begin{lem}\label{l:separation-unplanted}
For any $\theta\le 1$ and $\varepsilon>0$, with
overwhelming probability,
$
\zeta-\mu_+\ge
N^{-2/3-\varepsilon}
$.
\end{lem}
\begin{proof}
For the high temperature case $\theta<1$, the argument proving \eqref{e:saddle-separation-K0}, applied to the unplanted Gaussian matrix, gives $\zeta-\mu_+\ge \Omega(1)$ with overwhelming probability. This is stronger than the claim. At criticality, we use the argument proving \eqref{e:saddle-separation-lower-bound}, now applied to the unplanted matrix. In this case the interlacing comparison with the \((N-2)\times(N-2)\) minor is unnecessary. The same estimates yield \(\zeta-\mu_+\ge N^{-2/3-\varepsilon}\) with overwhelming probability.
We omit the
details.
\end{proof}
\begin{proof}[Proof of Lemma~\ref{l:hopfield-overlap-concentration}]
Let $\bx^1,\bx^2$ be two independent
replicas from the spherical Hopfield Gibbs measure and write
$
R_{1,2}
=
N^{-1}{\langle\bx^1,\bx^2\rangle}.
$
By Lemmas~\ref{l:hopfield-gibbs-covariance} and \ref{l:separation-unplanted}, 
for any $\theta\le 1$ and $0<\varepsilon<1/3$, there exists an event $\mathcal E_{\mathrm op}$ with overwhelming probability, on which
\[
\left\|
\left\langle\bx\bx^\top\right\rangle_{\bG}
\right\|_{\mathrm{op}}
\le
C\left\|
(\zeta I_N-\beta\bW)^{-1}
\right\|_{\mathrm{op}}
\le
CN^{2/3+\varepsilon}. 
\]
Conditionally on $\bG\in \mathcal E_{\mathrm{op}}$, independence of the two replicas gives (note that $\operatorname{Tr}\left\langle\bx\bx^\top\right\rangle_{\bG}=\langle \|\bx\|^2\rangle_\bG=N$)
\[
\begin{aligned}
\left\langle R_{1,2}^2\right\rangle_{\bG}
&=
\frac1{N^2}
\operatorname{Tr}\left(
\left\langle\bx\bx^\top\right\rangle_{\bG}^2
\right)\le
\frac1{N^2}
\left\|
\left\langle\bx\bx^\top\right\rangle_{\bG}
\right\|_{\mathrm{op}}
\operatorname{Tr}
\left\langle\bx\bx^\top\right\rangle_{\bG}=\frac{1}{N}\left\|
\left\langle\bx\bx^\top\right\rangle_{\bG}
\right\|_{\mathrm{op}}=o(1).
\end{aligned}
\]
Finally, Markov's inequality yields that for any $q_0>0$,
\[
\mathbb E
\left\langle
\mathbf{1}\left\{|R_{1,2}|>q_0\right\}
\right\rangle_{\bG}
\le \bbP[\mathcal E_{\mathrm{op}}^c]+
q_0^{-2}
\mathbb E\left[\mathbf{1}\{\mathcal E_{\mathrm{op}}\}
\left\langle R_{1,2}^2\right\rangle_{\bG}\right]
=o(1),
\]
which proves Lemma~\ref{l:hopfield-overlap-concentration}.
\end{proof}

%\newpage

\section{Proof of the spherical Hopfield free energy CLT} \label{rewrite:a:hopfield-clt}

In this section, we prove Theorem~\ref{t:hopfield-clt}. We will follow the notation of \S\ref{ss:coarea-upper-bound}. In particular, 
\[
 \beta=\frac{s-1}{s},\qquad \theta=\alpha(s-1)^2,\qquad \zeta_\star=1+\alpha(s-1).
\]
Let $\mu_i=\beta\lambda_i$ be the eigenvalues
of $\beta\bW$, where $\bW=\bG^\top\bG/N$, and put
$\mu_+=\max_i\mu_i$. Recall that the outer saddle $\zeta$ satisfies
$\sum_i(\zeta-\mu_i)^{-1}=N$ and $\zeta>\mu_+$.
For $z>\mu_+$, define
\begin{equation}\label{e:Phi-V}
 \Phi_N(z)=\frac N2(z-1)-\frac12\sum_{i=1}^N\log(z-\mu_i),\qquad
 V_N=\sum_{i=1}^N(\zeta-\mu_i)^{-2}.
\end{equation}
The coarea identity in Lemma~\ref{l:hopfield-coarea} and Stirling's
formula give
\begin{equation}\label{rewrite:e:app-coarea}
 \log Z_{\hop,N}^{\sph}
 =\Phi_N(\zeta)+\frac12\log(4\pi N)+\log f_\zeta(N)+O(N^{-1}),
 \qquad \Phi_N'(\zeta)=0.
\end{equation}
Here $f_\zeta$ is the density, conditional on $\bG$, of
$\sum_i(\zeta-\mu_i)^{-1}Z_i^2$, evaluated at its mean $N$.

The major inputs for this section are random matrix results that establishes estimates on $\log f_\zeta(N)$ and limiting distribution of $\Phi_N(\zeta)$. To this end, we introduce the deterministic analogue of $\Phi_N$ and the Stieltjes transform
\[
 \overline\Phi_N(z)=\frac N2\left(z-1-
   \int\log(z-\beta x)\,\nu_{\MP,\alpha}(\de x)\right),\qquad
 g(z)=\int\frac{\nu_{\MP,\alpha}(\de x)}{z-\beta x}.
\]
For $z$ outside the support of $\nu_{\MP,\alpha}$, it satisfies $z=g(z)^{-1}+\alpha\beta/(1-\beta g(z))$.
For $\theta<1$, direct calculation gives
\begin{equation}\label{rewrite:e:app-deterministic-saddle}
 g(\zeta_\star)=1,\qquad
 \overline\Phi_N(\zeta_\star)=\frac{\alpha N}{2}\log s,\qquad
 -g'(\zeta_\star)=\frac1{1-\theta}.
\end{equation}
At $\theta=1$, the first two identities extends by continuity. As in \S\ref{ss:rmt-outline}, $\zeta_\star$ stays a positive
distance from the full deterministic spectrum when $\theta<1$ and
equals its right endpoint when $\theta=1$.

Let $n=\min(M,N)$, $m=\max(M,N)$, and
$\omega=(\log\log n)^3$. We use the deterministic reference point
\begin{equation}\label{rewrite:e:reference-point}
 z_N=\begin{cases}
 \zeta_\star,&\theta<1,\\[2mm]
 \zeta_\star+\omega n^{-2/3},&\theta=1.
 \end{cases}
\end{equation}
At criticality, the shift therefore exceeds the edge scale, while the
error from replacing $\zeta$ by $z_N$ remains negligible on the fluctuation scale.

The proof follows the same three steps in both regimes. First we
estimate the density term in \eqref{rewrite:e:app-coarea}. Next we handle $\Phi_N(\zeta)$ by
replacing $\zeta$ by $z_N$, with an error negligible
on the fluctuation scale. Finally we apply a central limit theorem
to $\Phi_N(z_N)$, a logarithmic linear statistic evaluated at a
deterministic point. The following three lemmas make these steps precise.
Their proofs are deferred to the subsequent subsections.

\begin{lem}%[Density at the $\zeta$]
\label{rewrite:l:density}
The following estimates hold as $N\to\infty$.
\begin{enumerate}
\item If $0<\theta<1$, then
\[
 f_\zeta(N)=(1+o_{\bbP}(1))\sqrt{\frac{1-\theta}{4\pi N}},
\]
and consequently
\[
 \log f_\zeta(N)=-\frac12\log N
 +\frac12\bigl(\log(1-\theta)-\log(4\pi)\bigr)+o_{\bbP}(1),
\]
\item If $\theta=1$, then
\[
 f_\zeta(N)=O_{\bbP}(N^{-2/3}),\qquad
 f_\zeta(N)^{-1}=O_{\bbP}(N^{2/3}).
\]
and consequently
\[
 \log f_\zeta(N)=-\frac23\log N+O_{\bbP}(1).
\]
\end{enumerate}
\end{lem}

\begin{lem}%[Replacement of $\zeta$ by $z_N$]
\label{rewrite:l:comparison}
For $z_N$ defined in \eqref{rewrite:e:reference-point},
\[
 \Phi_N(\zeta)-\Phi_N(z_N)=
 \begin{cases}
 o_{\bbP}(1),&0<\theta<1,\\
 o_{\bbP}(\sqrt{\log N}),&\theta=1.
 \end{cases}
\]
\end{lem}

\begin{lem}%[CLT at $z_N$]
\label{rewrite:l:action-clt}
The following limits hold as $N\to\infty$.
\begin{enumerate}
\item If $0<\theta<1$, then
\[
 \Phi_N(z_N)-\frac M2\log s
 \ \Longrightarrow\ 
 \mathcal N\left(-\frac14\log(1-\theta),
                 -\frac12\log(1-\theta)\right).
\]
\item If $\theta=1$, then
\[
 \frac{\Phi_N(z_N)-\frac M2\log s-\frac1{12}\log N}
 {\sqrt{\frac16\log N}}
 \ \Longrightarrow\ \mathcal N(0,1).
\]
\end{enumerate}
\end{lem}

\begin{proof}[Proof of Theorem~\ref{t:hopfield-clt}]
If $s=1$, then $\beta=0$ and $Z_{\hop,N}^{\sph}=1$, so the conclusion
is immediate. Assume $s\ne1$.

In the high-temperature regime, Lemmas~\ref{rewrite:l:density}
and~\ref{rewrite:l:comparison}, together with
\eqref{rewrite:e:app-coarea}, give
\[
 \log Z_{\hop,N}^{\sph}-\frac M2\log s
 =\Phi_N(z_N)-\frac M2\log s
   +\frac12\log(1-\theta)+o_{\bbP}(1).
\]
By
Lemma~\ref{rewrite:l:action-clt}, the right-hand side converges to
$\mathcal N(\frac14\log(1-\theta),-\frac12\log(1-\theta))$, proving the high-temperature assertion.

At criticality, the same two lemmas yield
\[
 \log Z_{\hop,N}^{\sph}-\frac M2\log s
 =\Phi_N(z_N)-\frac M2\log s-\frac16\log N
   +o_{\bbP}(\sqrt{\log N}).
\]
Here the deterministic constant
$\frac12\log(4\pi)$ and the $O_{\bbP}(1)$ density error are
absorbed in the remainder. Thus
\[
 \frac{\log Z_{\hop,N}^{\sph}-\frac M2\log s+\frac1{12}\log N}
 {\sqrt{\frac16\log N}}
 =\frac{\Phi_N(z_N)-\frac M2\log s-\frac1{12}\log N}
 {\sqrt{\frac16\log N}}+o_{\bbP}(1).
\]
Lemma~\ref{rewrite:l:action-clt} then completes the proof.
\end{proof}

\subsection{Proof of Lemma~\ref{rewrite:l:density}}
\label{rewrite:ss:density-proof}

We first estimate $V_N$ from \eqref{e:Phi-V}. We record these
estimates separately, since they are also needed for the comparison
in Lemma~\ref{rewrite:l:comparison}.

\begin{lem}%[Squared resolvent weights]
\label{rewrite:l:variance}
If $0<\theta<1$, then
\[
 V_N=\frac{N}{1-\theta}(1+o_{\bbP}(1)),\qquad
 |\zeta_\star-\zeta|=O_{\bbP}(N^{-1}).
\]
If $\theta=1$, then
\[
 \begin{gathered}
 V_N=O_{\bbP}(N^{4/3}),\qquad V_N^{-1}=O_{\bbP}(N^{-4/3}),\\
 |\zeta-\zeta_\star|=O_{\bbP}(N^{-2/3}),\qquad
 (\zeta-\mu_+)^{-1}=O_{\bbP}(N^{2/3}).
 \end{gathered}
\]
\end{lem}

\begin{proof}
\noindent\emph{High temperature.}
By \cite[Theorem~1]{Johnstone2001Largest} and \cite[Theorem~2]{Ma2012Wishart}, $\mu_+$ converges in probability to the relevant spectral edge $\mu_\star$, which is strictly less than $\zeta_\star$ by Fact~\ref{f:zeta-above-spectrum}. Thus, the spectrum is separated
from $\zeta_\star$ with probability tending to one. The analytic linear
statistics theorem \cite[Theorem~1.1]{BaiSilverstein2004}, applied to
the functions $(\zeta_\star-x)^{-1},(\zeta_\star-x)^{-2}$, gives
\[
 \sum_{i=1}^N(\zeta_\star-\mu_i)^{-1}=N+O_{\bbP}(1),\qquad
 \sum_{i=1}^N(\zeta_\star-\mu_i)^{-2}=\frac{N}{1-\theta}+O_{\bbP}(1),
\]
where the $N$-term coefficient values
follow from \eqref{rewrite:e:app-deterministic-saddle}.\footnote{More precisely, the cited theorem gives that for $f(x)=(\zeta_\star-x)^{-j},j=1,2$, $\sum_{i}f(\mu_i)-N\int f(\beta x)\nu_{\MP,M/N}(\de x)=O_\bbP(1)$. We further compare the integration with respect to $\nu_{\MP,M/N}$ to that of $\nu_{\MP,\alpha}$, and the error is $O(N^{-1})$ since $M/N-\alpha=O(N^{-1})$. This type of comparison between $\nu_{\MP,M/N}$ and $\nu_{\MP,\alpha}$ will be used in several places below. All of the errors arising this way can be absorbed in the corresponding error bounds, so in what follows we omit these steps.}
The first sum is decreasing, with derivative of magnitude comparable
to $N$ in a fixed neighborhood of $\zeta_\star$. The saddle equation therefore
gives $\zeta-\zeta_\star=O_{\bbP}(N^{-1})$. The derivative of the second sum is
$O_{\bbP}(N)$ in the same neighborhood, so
\[
 V_N=\sum_i(\zeta_\star-\mu_i)^{-2}+O_{\bbP}(N\cdot N^{-1})
     =\frac{N}{1-\theta}+O_{\bbP}(1),
\]
which proves the high-temperature estimate.

\noindent\emph{Criticality.} For $\theta=1$, $\zeta_\star$ coincides with the relevant spectral edge. Set $\rho=N^{-2/3}$. 
The edge limits of \cite[Theorem~1]{Johnstone2001Largest} and
\cite[Theorem~2]{Ma2012Wishart} give $\mu_+-\zeta_\star=O_{\bbP}(\rho)$.
We use the sharp edge resolvent estimates from
\cite[Corollary~2.7, Lemma~2.8 and Appendix~A.2]{collins2025fluctuations}
for $\beta>0$, and
\cite[Theorem~7.9 and Remark~7.10]{HuangZhang2026Airy}
at the lower soft edge for $\beta<0$. Writing
$r_N(z)=N^{-1}\sum_{i=1}^N(z-\mu_i)^{-1}$, these imply
\begin{equation}\label{rewrite:e:variance-edge-resolvent}
 \begin{gathered}
 r_N(\zeta_\star+\mathrm i\rho)=1+O_{\bbP}(N^{-1/3}),\qquad
 \lim_{L\to\infty}\limsup_{N\to\infty}
 \bbP\bigl(\Re r_N(\zeta_\star+(2+\mathrm i)L\rho)\ge1\bigr)=0.
 \end{gathered}
\end{equation}
Here we use $g(\zeta_\star+\mathrm i\rho)=1+O(\sqrt\rho)$ and
$\Re g(\zeta_\star+(2+\mathrm i)L\rho)\le1-c\sqrt{L\rho}$.

We now establish the saddle location estimate. On the event $\{\mu_+\le \zeta_\star+L\rho\}$, let
$u_i=\zeta_\star+2L\rho-\mu_i\ge L\rho$. Since $(L\rho)^2\le a_iL\rho$ for any $i$, 
\[
 r_N(\zeta_\star+3L\rho)=\frac1N\sum_i\frac1{u_i+L\rho}
 \le\frac1N\sum_i\frac{a_i}{a_i^2+(L\rho)^2}
 =\Re r_N(\zeta_\star+(2+\mathrm i)L\eta ).
\]
As $L\to\infty$, the right end is smaller than $1$ with probability tending to $1$ by \eqref{rewrite:e:variance-edge-resolvent}.
Monotonicity therefore gives
$|\zeta-\zeta_\star|=O_{\bbP}(\rho)$, and hence $\zeta-\mu_+=O_{\bbP}(\rho)$.

We next show $(\zeta-\mu_+)^{-1}=O_\bbP(\rho^{-1})$. 
For any fixed $K$, on the event $\{|\zeta-\zeta_\star|\le K\rho\}$, we have
$$|\zeta_\star+i\rho-\mu_i|\le|\zeta+i\rho-\mu_i|+K\rho\le(1+K)|\zeta+i\rho-\mu_i|.$$
On this event, direct computation gives
\[
 \begin{aligned}
 N|r_N(\zeta_\star+i\rho)-r_N(\zeta+i\rho)|
 &\le|\zeta-\zeta_\star|\sum_i\frac1{|\zeta_\star+i\rho-\mu_i|\,|\zeta+i\rho-\mu_i|}\\
 &\le|\zeta-\zeta_\star|\sum_i(1+K)|\zeta_\star+i\rho-\mu_i|^{-2}\\
 &=(1+K)K\rho\cdot (-N\rho^{-1}\Im r_N(\zeta_\star+i\rho))=O_\bbP(\rho^{-1})\,.
 \end{aligned}
\]
Taking $K\to\infty$ and combining with \eqref{rewrite:e:variance-edge-resolvent} yields that $|Nr_N(\zeta+i\rho)-N|=O_\bbP(\rho^{-1})$, which implies
\[
 O_{\bbP}(\rho^{-1})
 =N-N\Re r_N(\zeta+\mathrm i\rho)
 =\sum_i\frac{\rho^2}{(\zeta-\mu_i)((\zeta-\mu_i)^2+\rho^2)}
 \ge\frac{\rho^2}{d(d^2+\rho^2)}.
\]
This proves $(\zeta-\mu_+)^{-1}=O_{\bbP}(\rho^{-1})$.

To estimate $V_N$ we use
\[
 \sum_{j=2}^N(\mu_+-\mu_j)^{-2}=O_{\bbP}(\rho^{-2}).
\]
For $\beta>0$, this is
\cite[Lemma~4.1, $l=2$, and Corollary~2.7]{collins2025fluctuations};
any zero eigenvalues contribute only $O(N)$.
For $\beta<0$, criticality implies $\alpha=(1-s)^{-2}>1$,
and the estimate follows from \cite[Theorem~3]{LiValkoXu2026LowerEdge}.
Finally, monotonicity gives
\[
 (\zeta-\mu_+)^{-2}\le V_N\le (\zeta-\mu_+)^{-2}
       +\sum_{j=2}^N(\mu_+-\mu_j)^{-2}=(\zeta-\mu_+)^{-2}+O_{\bbP}(\rho^{-2}).
\]
Since $\zeta-\mu_+=O_{\bbP}(\rho)$ and $(\zeta-\mu_+)^{-1}=O_{\bbP}(\rho^{-1})$, this gives $V_N=O_\bbP(\rho^{-2})$ and
$V_N^{-1}=O_{\bbP}(\rho^2)$.
\end{proof}

\begin{proof}[Proof of Lemma~\ref{rewrite:l:density}]
Condition on $\bG$ and set $a_i=(\zeta-\mu_i)^{-1}$. Then
$\sum_i a_i=N$ and $\sum_i a_i^2=V_N$.
At criticality, Lemma~\ref{l:weighted-chi-square-density} gives
$C^{-1}V_N^{-1/2}\le f_\zeta(N)\le CV_N^{-1/2}$, so
Lemma~\ref{rewrite:l:variance} proves both asserted stochastic bounds.

In the high-temperature regime, the weights $a_i$ are uniformly bounded and $V_N/N \to (1-\theta)^{-1}$ with probability tending to one. In particular, $\max\{ a_i/\sqrt V_N \}=O_{\bbP}(N^{-1/2})$. Apply the estimate in \eqref{eq:Fourier-est-local} and \eqref{eq:Fourier-est-far} to the normalized weights $a_i/\sqrt V_N$, with $\delta$ replaced by $\max\{ a_i/\sqrt V_N \}$. After rescaling the density, they give
\[
 f_\zeta(N)=(4\pi V_N)^{-1/2}(1+o_{\bbP}(1))
 =(1+o_{\bbP}(1))\sqrt{\frac{1-\theta}{4\pi N}}.
\]
The logarithmic conclusions follow immediately.
\end{proof}

\subsection{Proof of Lemma~\ref{rewrite:l:comparison}}
\label{rewrite:ss:comparison-proof}

\begin{proof}
Since $\Phi_N'(\zeta)=0$ and
\[
 \Phi_N''(z)=\frac12\sum_i(z-\mu_i)^{-2},\qquad z>\mu_+,
\]
we estimate the difference by Taylor's formula.

In the high-temperature regime, $z_N=\zeta_\star$. As argued before, $\mu_+=\mu_\star+o(1)$ with overwhelming probability.
Fact~\ref{f:zeta-above-spectrum} and Lemma~\ref{rewrite:l:variance} imply that $\zeta_\star$ and $\zeta$ lie in a
fixed neighborhood separated from the spectrum, with probability
tending to one. On the interval between them, $\Phi_N''(z)=O_{\bbP}(N)$.
Consequently,
\[
 |\Phi_N(\zeta)-\Phi_N(\zeta_\star)|
 =O_{\bbP}\bigl(N|\zeta-\zeta_\star|^2\bigr)
 =O_{\bbP}(N^{-1})=o_{\bbP}(1).
\]

At criticality, $z_N-\zeta_\star=\omega n^{-2/3}\asymp\omega N^{-2/3}$ and $\omega\to\infty$.
Lemma~\ref{rewrite:l:variance} therefore gives $z_N>\zeta$ with
probability tending to one and
$z_N-\zeta=O_{\bbP}(\omega N^{-2/3})$. For $z\ge\zeta$,
$0<\Phi_N''(z)\le V_N/2$, so
\[
 0\le\Phi_N(z_N)-\Phi_N(\zeta)
 \le\frac14(z_N-\zeta)^2V_N
 =O_{\bbP}(\omega^2)
 =o_{\bbP}(\sqrt{\log N}). \qedhere
\]
\end{proof}

\subsection{Proof of Lemma~\ref{rewrite:l:action-clt}}
\label{rewrite:ss:clt-proof}

For the centering corrections we use the identity
\begin{equation}{\label{eq:aspect-ratio-derivative}}
    \partial_\gamma\int\log(z-\beta x)\,\nu_{\MP,\gamma}(\de x) =\log(1-\beta g_\gamma(z)),
\end{equation}
where $g_\gamma(z)=\int (z-\beta x)^{-1} \nu_{\MP,\gamma}(\de x)$. 
This follows by differentiating the MP potential formula
\eqref{e:wishart-saddle-potential}; the terms containing
$\partial_\gamma g_\gamma$ cancel by the MP equation.

Write $\lambda=\min\{\alpha,\alpha^{-1}\}=\lim n/m\le 1$, so
$n/m=\lambda+O(n^{-1})$, and let $\mathrm{MP}_\lambda$ denote the normalized
MP law supported on $[(1-\sqrt\lambda)^2,(1+\sqrt\lambda)^2]$.
Note that $\nu_{\MP,\lambda}=(1-\lambda)\delta_0+\lambda\,\mathrm{MP}_\lambda$.

\begin{proof}
\noindent\emph{High temperature.}
Define the statistic centered at the finite aspect ratio,
\[
 \widetilde L_N^{\mathrm{ht}}
 =N\int\log(\zeta_\star-\beta x)\,\nu_{\MP,M/N}(\de x)
  -\sum_{i=1}^N\log|\zeta_\star-\mu_i|.
\]
The analytic linear statistics CLT of
\cite[Theorem~1.1]{BaiSilverstein2004}, with the mean and variance
evaluated by \cite[\S3.3, (3.6)--(3.7) and (3.27)]{baik2016fluctuations},
gives
\[
 \widetilde L_N^{\mathrm{ht}}\Longrightarrow
 \cN\left(-\frac12\log(1-\theta),-2\log(1-\theta)\right).
\]
By \eqref{rewrite:e:app-deterministic-saddle} and \eqref{eq:aspect-ratio-derivative},
\begin{align*}
    &\frac N2(\zeta_\star-1)-\frac N2\int\log(\zeta_\star-\beta x)\,\nu_{\MP,M/N}(\de x)  \\ 
    =\ & \overline{\Phi}_N(\zeta_\star) + \frac{N}{2} (\alpha-M/N) \left( \log(1-\beta g_\alpha(\zeta_\star))+O(N^{-1}) \right)  \\ 
    =\ & \frac{\alpha N}{2}\log s + \frac{N}{2} (\alpha-M/N) \left( -\log s+O(N^{-1}) \right) = \frac M2\log s+O(N^{-1}).
\end{align*}
Consequently, since $z_N=\zeta_\star$,
\[
 \Phi_N(z_N)-\frac M2\log s
 =\frac12\widetilde L_N^{\mathrm{ht}}+o_{\bbP}(1),
\]
which proves the high-temperature assertion.

\noindent\emph{Criticality.}
Let $x_1,\ldots,x_n$ be the eigenvalues of $m^{-1}\bX\bX^\top$,
where $\bX=\bG$ if $M\le N$ and $\bX=\bG^\top$ otherwise. Define
\[
 v_N=\frac{Nz_N}{\beta m}.
\]
Since $m/N=\max\{1,\alpha\}+O(n^{-1})$, at criticality there is a
deterministic multiplier
\[
 \widehat{\omega}
 =\frac{\omega}{|\beta|\max\{1,\alpha\}}+O(n^{-1/3})
 \asymp\omega
\]
such that
\[
 v_N=\begin{cases}
 (1+\sqrt\lambda)^2+\widehat{\omega}n^{-2/3},&\beta>0,\\[2mm]
 (1-\sqrt\lambda)^2-\widehat{\omega}n^{-2/3},&\beta<0.
 \end{cases}
\]
Thus the $O(n^{-1})$ correction to the evaluation point is absorbed
into an $O(n^{-1/3})$ change of $\widehat{\omega}$.
The edge CLTs used below apply throughout the range
$\widehat{\omega}\asymp(\log\log n)^3$. The statistic
\[
 \widetilde L_N=n\int\log|v_N-x|\,\mathrm{MP}_{\lambda}(\de x)
                 -\sum_{i=1}^n\log|v_N-x_i|
\]
satisfies the CLT of \cite[Theorem~1.1 and Lemma~3.1]{CollinsWoodfinLe2025Edge}
when $\beta>0$, and of Lemma~\ref{rewrite:l:app-lower-edge-clt}
when $\beta<0$. 
In the lower-edge case, criticality
implies $\alpha=(1-s)^{-2}>1$, so $\lambda$ stays away from $0$ and $1$. Thus
\begin{equation}\label{rewrite:e:critical-logdet}
 \frac{\widetilde L_N-\frac16\log N}{\sqrt{\frac23\log N}}
 \Longrightarrow\cN(0,1).
\end{equation}

It remains to compare this statistic with $\Phi_N(z_N)$. The exact
coordinate change $z_N=\beta m v_N/N$ and
$\lambda-n/m=O(n^{-1})$, together with \eqref{eq:aspect-ratio-derivative}, give
\[
 \Phi_N(z_N)=\overline\Phi_N(z_N)+\frac12\widetilde L_N
              +O_{\bbP}(1).
\]
Moreover,
the square-root edge expansion $1-g(\zeta_\star+v)=O(\sqrt v)$ gives
\begin{align}
    \overline\Phi_N(z_N)= \overline{\Phi}_N(\zeta_\star) + \frac{N}{2} \int_0^{z_N-\zeta_\star} \left(1-g(\zeta_\star+v)\right) \de v = \frac M2\log s+O(\omega^{3/2}).
\end{align}
As $\omega^{3/2}=(\log\log n)^{9/2}=o(\sqrt{\log N})$, we conclude that
\[
 \Phi_N(z_N)-\frac M2\log s
 =\frac12\widetilde L_N+o_{\bbP}(\sqrt{\log N}).
\]
Equation~\eqref{rewrite:e:critical-logdet} proves the assertion.
\end{proof}

We now prove the lower-edge counterpart used above.

\begin{lem}\label{rewrite:l:app-lower-edge-clt}
Fix $\lambda\in(0,1)$ and let $n\le m$ satisfy
$n/m=\lambda+O(n^{-1})$,
and $\bX\in\bbR^{n\times m}$ have independent standard Gaussian entries.
Let $\widehat{\omega}\asymp(\log\log n)^3$, and set
$\bS=m^{-1}\bX\bX^\top$ and
$u=(1-\sqrt\lambda)^2-\widehat{\omega}n^{-2/3}$. Then
\begin{equation}{\label{eq:lower-edge-clt}}
    \frac{\log|\det(\bS-uI)|-n\int\log(x-u)\,\mathrm{MP}_\lambda(\de x)
       +\frac16\log n}{\sqrt{\frac23\log n}}
 \ \Longrightarrow\ \mathcal N(0,1).
\end{equation}
\end{lem}
\begin{proof}
We adapt the proof of \cite{CollinsWoodfinLe2025Edge}, recording the necessary changes.
By Lemma~\ref{l:golub-kahan}, $\bX\bX^\top$ has the same eigenvalue law as $\bL\bL^\top$, where $\bL$
is lower bidiagonal with independent diagonal entries $a_i$ and subdiagonal entries $b_i$, where $a_i^2\sim\chi^2_{m-n+i}$ and $b_i^2\sim\chi^2_i$. 
Its $i$'th leading principal minor $D_i=\det((\bT-muI)_{[1,i]})$ satisfy $D_0=1$, $D_1=a_1^2-mu$, and
\begin{equation}{\label{eq:recur-D-i}}
    D_i=(a_i^2+b_{i-1}^2-mu)D_{i-1}-a_{i-1}^2b_{i-1}^2D_{i-2}.
\end{equation}
The deterministic analog of this recursion is given by
\[
    D_i' = A_i D_{i-1}' - p_i D_{i-2}', \qquad A_i=m-n+2i-1-mu, \qquad p_i=(m-n+i-1)(i-1),
\]
which has characteristic roots
\begin{equation}{\label{eq:characteristic-roots}}
    \rho_i^{\pm} =(A_i\pm\sqrt{\Delta_i}\,)/2, \qquad \Delta_i=A_i^2-4p_i.
\end{equation}
This motivates us to consider the normalized recursion by putting $P_i=D_i/\prod_{j\le i}\rho_j^+$. The root-product calculation in
\cite[Lemma~3.1]{CollinsWoodfinLe2025Edge} applies with $\rho_i^+$, leading to 
\[
    \sum_{i=1}^n\log\frac{\rho_i^+}{m} =n\int\log(x-u)\,\mathrm{MP}_{n/m}(\de x)+O(1)=n\int\log(x-u)\,\mathrm{MP}_{\lambda}(\de x)+O(1).
\]
Combined with 
\[
    \log|\det(\bS-uI)|\overset d=\log|P_n|+\sum_i\log(\rho_i^+/m),
\]
it suffices to prove
\begin{equation}{\label{eq:reduced-clt}}
    \frac{ \log|P_n|+\frac{1}{6}\log n }{ \sqrt{\frac{2}{3}\log n} } \Longrightarrow \cN(0,1).
\end{equation}
For simplicity, set
\[
    \tau_i=\frac{m-n+i}{\rho_i^+},\quad \delta_i=\frac{i-1}{\rho_i^+},\quad u_i=\frac{a_i^2-(m-n+i)}{\rho_i^+},\quad v_i=\frac{b_{i-1}^2-(i-1)}{\rho_i^+},\quad v_1=0.
\]
The recursion then can be written as $P_0=1,P_1=1+u_1$, and
\begin{equation}{\label{eq:recur-P-i}}
    P_i = \left( \tau_i+\delta_i-\frac{mu}{\rho_i^+}+u_i+v_i \right) P_{i-1} - (\tau_{i-1}+u_{i-1})(\delta_i+v_i) P_{i-2}.
\end{equation}
Our general approach, modeled after the methods in \cite{CollinsWoodfinLe2025Edge}, is to approximate the recursion for $P_i$ by a linear recursion. The observation is that $P_i/P_{i-1}$ is close to $1$, so we introduce the quantity $R_i=1-P_i/P_{i-1}$. 
Writing $X_i=(1-\tau_{i-1})(\delta_i u_{i-1}-v_i)$ and
$Z_i=u_{i-1}v_i+\delta_i u_{i-1}+\tau_{i-1}v_i$, the exact recurrence for $R_i$ is (below we write $\omega_i=\rho_i^-/\rho_{i-1}^+=\tau_{i-1}\delta_i$ and $d_i=\rho_i^-/\rho_i^+-\omega_i$)
\begin{equation}{\label{eq:recur-R-i}}
    (R_i+u_i) = X_i+\omega_i(R_{i-1}+u_{i-1})-d_i+u_{i-1}v_i +\frac{Z_iR_{i-1}+\omega_iR_{i-1}^2}{1-R_{i-1}},\qquad R_1=-u_1.
\end{equation}
We now define a linear recursion on a new sequence of variables $Y_i'$, which we will show are a good approximation of $R_i+u_i$. Let
\begin{equation}{\label{eq:recur-Y-i}}
    Y_i'=(X_i-d_i)+\omega_iY'_{i-1}, \qquad Y_1'=0.
\end{equation}
Put $G_{n+1}=1$ and $G_i=1+\omega_iG_{i+1}$. For this recursion, we have
\begin{equation}{\label{eq:Y-sums-of-id-variable}}
    \sum_{i=1}^{n} Y_i' = \sum_{i=1}^{n} Y_i - \sum_{i=2}^{n} G_{i+1} d_i,  \qquad Y_i = \sum_{ 2 \le j \le i } W_{j,i} X_j, \qquad W_{i,i}=1,\ W_{j,i}=\prod_{k=j+1}^i\omega_k.  
\end{equation}
Before applying the arguments in \cite[Section~4]{CollinsWoodfinLe2025Edge}, we first collect some deterministic estimates for the parameters.
\begin{fac}\label{f:deterministic-parameter-bound}
  Define $\eta_i=((n-i+1)/n+\widehat{\omega} n^{-2/3})^{1/2}$. 
  \begin{itemize}
      \item By \eqref{eq:characteristic-roots}, $\Delta_i \asymp n^2\eta_i^2$, $A_i\asymp n$, $\rho_i^+\asymp n$, and $\rho_i^-\ge0$. Also, by the identities
      \[
          1-\omega_i=\frac{\sqrt{\Delta_{i-1}}+\sqrt{\Delta_i}-2}{2\rho_{i-1}^+}, \qquad \rho_{i-1}^+-\rho_i^+=-1+\frac{2(mu-1)}{\sqrt{\Delta_{i-1}}+\sqrt{\Delta_i}},
      \] 
      we have $0<\omega_i<1$, $1-\omega_i\asymp\eta_i$, $|d_i|\le C/(n\eta_i)$, and $|\omega_{i+1}-\omega_i|\le C/(n\eta_i)$.
      Then
      \begin{equation}\label{rewrite:eq:a:lower-green}
         W_{j,i}\le e^{-c(i-j)\eta_i},\qquad \sum_{j=2}^iW_{j,i}\le C\eta_i^{-1},\qquad G_i,G_{i+1}\le C\eta_i^{-1}.
      \end{equation}
      The last bound follows by backward induction: since $\eta_{i+1}^{-1}-\eta_i^{-1}\le C/(n\eta_i^3)$ and $n\eta_i^3\ge\widehat{\omega}^{3/2}\to\infty$, $\eta_i^{-1}-\omega_i\eta_{i+1}^{-1}\ge c/2$ for large $n$.
      \item The $X_i$ are independent and centered, and the chi-square generating function gives 
      \[
          u_i,v_i,X_i\in\mathrm{SG}(C/n,C/n)
      \]
      in the notation of \cite[Definition~2.4]{CollinsWoodfinLe2025Edge}. Moreover,
      \[
          \bbE X_i^2=2\delta_i(1-\tau_{i-1})^2 \left(\frac{\omega_i}{\rho_{i-1}^+}+\frac1{\rho_i^+}\right), \qquad \bbE X_i^4\le Cn^{-2}.
      \]
  \end{itemize}
\end{fac}
The stopping and truncation arguments in \cite[Lemmas~4.1--4.4 and 4.7--4.8]{CollinsWoodfinLe2025Edge} apply to the recurrence \eqref{eq:recur-Y-i}: their sub-gamma estimates are unchanged, and their propagation bounds follow from the deterministic estimates \eqref{rewrite:eq:a:lower-green} and $|d_i|\le C/(n\eta_i)$. Consequently their decompositions (4.6) and (4.21) give
\begin{align}
    \log|P_n| 
    = -\sum_{i=1}^nY_i+\sum_{i=2}^nG_{i+1}d_i-Q_n+O_{\bbP}(1), \qquad Q_n=\sum_{i=2}^n(G_i-1)Y_{i-1}^2. \label{rewrite:eq:a:lower-linearization}
\end{align}
Write $y_i=(n-i)/n$ and restrict first to
$\mathcal I_n=\{i:(\log\log n)\eta_n^2\le y_i\le(\log\log n)^{-2}\}$.
With $e_i=\sqrt{y_i}+\eta_n^2/y_i$, the root identities yield
\[
 \begin{split}
 1-\omega_i&=2(1-\sqrt\lambda)\sqrt{y_i}\{1+O(e_i)\},\qquad
 \bbE X_i^2=\frac{4(1-\sqrt\lambda)^2}{n}\{1+O(e_i)\},\\
 G_{i+1}^2\bbE X_i^2&=\frac1{ny_i}\{1+O(e_i)\},\qquad
 G_{i+1}d_i=\frac1{4ny_i}\{1+O(e_i)\}.
 \end{split}
\]
Summing the window estimates and using
$\sum_{i\notin\mathcal I_n}(n\eta_i^2)^{-1}=O(\log\log n)$ gives
\[
 \operatorname{Var}\left(\sum_iY_i\right)
 =\sum_{i=2}^nG_{i+1}^2\bbE X_i^2
 =\frac23\log n+O(\log\log n),\qquad
 \sum_{i=2}^nG_{i+1}d_i=\frac16\log n+O(\log\log n).
\]
For $i\le j$, the exact identities
$\operatorname{Cov}(Y_i,Y_j)=W_{i,j}\bbE Y_i^2$ and
$\operatorname{Cov}(Y_i^2,Y_j^2)=W_{i,j}^2\operatorname{Var}(Y_i^2)$ give
\[
 \bbE Q_n=\frac12\left\{\operatorname{Var}\left(\sum_iY_i\right)
                   -\sum_i\bbE Y_i^2\right\}
 =\frac13\log n+O(\log\log n),\qquad
 \operatorname{Var}(Q_n)\le\frac{C}{n\eta_n^3}=o(1),
\]
using \eqref{rewrite:eq:a:lower-green} and
$\bbE Y_i^4\le C/(n^2\eta_i^2)$.
Thus \eqref{rewrite:eq:a:lower-linearization} becomes
\[
    \log|P_n|=-\sum_iY_i-\frac16\log n+O_{\bbP}(\log\log n)
\]
Using Fact~\ref{f:deterministic-parameter-bound}, it is easy to check that $Y_i$ in the form of \eqref{eq:Y-sums-of-id-variable} satisfies the Lyapunov condition. So \eqref{eq:reduced-clt} holds.
\end{proof}

\bibliographystyle{alpha}
\bibliography{bib}

\end{document}